\documentclass[a4paper,10pt]{amsart}

\include{amssymb}
\usepackage{amssymb}

\newcommand{\R}{\mathbb{R}^{n}}

\newcommand{\B}{\mathcal{B}}
\newcommand{\DL}{\mathcal{D}_{L^{1}}'}

\newcommand{\ram}{r_{T}'}
\newcommand{\vsp}{\vspace*{.5cm}}
\newcommand{\ra}{\mathfrak{A}}
\newcommand{\Ss}{\text{ sing supp }}
\newcommand{\hyp}{\mathcal{B}_{K}}
\newcommand{\Bm}{\mathcal{B}_{m}}
\newcommand{\FB}{\mathcal{F}}

\newcommand{\M}{\mathcal{M}}

\newcommand{\Top}{\mathcal{T}}
\begin{document}

\title{SOME REMARKS ON TR\`EVES' CONJECTURE}
\author{T. Dahn}
\thanks{Tove Dahn, Lund University, Sweden}

\begin{abstract}
 I will give a discussion of the conditions involved in Treves' conjecture on analytic hypoellipticity. I will discuss some microlocally characteristic sets and introduce a topology of monotropic functionals as suitable for solving the conjecture. The pseudodifferential operator representation is inspired by Cousin (cf. \cite {Cousin95})
\end{abstract}

\maketitle

\section{Introduction}
Treves' conjecture is existence of an involutive stratification equivalent with hypoellipticity. The concept of hypoellipticity is very sensitive to change of topology but there are geometric sets that are characteristic. We will discuss lineality and a set that relates to orthogonality. We will also consider three sets that occur in literature and that we consider as not characteristic. The first is relative representation of spectral function to a hypoelliptic operator (\ref{Nilsson}). The third (section \ref{Ex}) relates to hypoelliptic operators as limits of operators dependent on a parameter . In the second we consider continuation of the contact transform to ($\Top$), which is considered as a B\"acklund transform. For this continuation, algebraicity is considered to be characteristic for hypoellipticity (\ref{Iversen}). It is necessary for hypoellipticity that the singularities have measure zero and in this study we assume parabolic singularities. The regular approximations are transversals and we only briefly discuss some possible generalizations..

\vsp

The set of lineality is defined for a polynomial over a real (or
complex) vector space $E_{\mathbb{R}}$ (or $E_{\mathbb{C}})$ is
$$ \Delta(P)=\{ \eta \in E_{\mathbb{R}} \quad P(\xi + i t
\eta)-P(\xi)=0 \quad \forall \xi \in E_{\mathbb{R}} \quad \forall t
\in \mathbb{R} \}$$ It can be proved that $\Delta$ is standard
complexified in the topology for $Exp$ (cf. \cite{Martineau} ), why
it is sufficient to consider purely imaginary translations as above.
The set can be generalized to symbol classes where $\Delta$ has a
locally algebraic definition or where the definition is locally
algebraic modulo monotropy. The pseudo differential operators are
realized from the symbol ideals using a representation derived from Cousin.

\vsp

For constant coefficients polynomial differential operators, we note that the class of operators hypoelliptic in $\mathcal{D}'$ is not radical. We can prove for the radical to the class, that the lineality is decreasing for iteration. For variable coefficients polynomial differential operators, we consider formally hypoelliptic operators, that is where the symbol is equivalent in strength with a constant coefficients polynomial operators, as the variable varies. We also assume that the real part of the symbol is unbounded and does not change sign in the infinity, as the variable varies on a connected set.

\vsp

The generalization to more general symbols will be using a lifting operator acting on a dynamical system, that maps into analytic symbols $f(\zeta) \in (I)(\Omega)$, where $(I)$ is an ideal over a pseudo convex domain. We will mainly discuss operators $A_{\lambda}$ on the form $A_{\lambda}=P_{\lambda}+H_{\lambda}$, where $P_{\lambda}$ is a polynomial for finite parameter values and $H_{\lambda}$ is regularizing.
 
 \newtheorem{mainresult}{ Proposition }[section]
 \begin{mainresult} \label{conjecture}
 Assume $S$ a pseudo differential operator, self-adjoint and of
 exponentially finite type. Assume the symbol in $(I)(\Omega)$,
 where $(I)$ is finitely generated and $\Omega$ is pseudo convex.
 Assume the lineality to $S$, $\Omega_{0}$ is decreasing for
 iteration. Assume singular points are mapped on to singular points
 in the dynamical system, with tangent determined (global pseudo
 base). Then, for $u \in \mathcal{D}'$
 $$ WF_{a}(Su)=\tilde{\Omega}_{0} \cup WF_{a}(u)$$
 \end{mainresult}
Here $\tilde{\Omega}$ is a set only dependent on $\Omega$ and the symbol.
\vsp

If $\Omega_{0}=\lim_{j}\Omega_{j}$, where $\Omega_{j} \subset$ a
pseudo convex set (and $\Omega_{j}$ algebraic), then $\Omega_{0}$
must be an analytic set.
Given that the level surfaces are of order $1$, $\Omega_{0}$ has a
locally algebraic definition through transversality. Conversely, if
$\Omega_{0}$ has an algebraic definition and if we have a global
pseudo base for $(I)$, then regular approximations are transversals
and $\Omega_{0}$ is given by regular approximations.
If $\Omega_{j}$ are given by the lineality locally to $A_{\lambda}=P_{\lambda}+H_{\lambda}$
and $P_{\lambda} \sim P^{\lambda}$, then $\Omega_{0}$ is a set of lineality for the limit of $P_{\lambda}$.

\vsp

We will use monotropic functionals to study both the symbols to hypoelliptic operators
and the equations in operator space. For the representations we consider, monotropy is microlocally indifferent, that is does not influence the geometry in a microlocally significant way. We will use the notation $f \sim_{m}0$  explained as follows.
Between the spaces $\dot{\B}(\mathbb{R}^{n})$ and
$\B(\mathbb{R}^{n})$, we consider over an $\epsilon-$ neighborhood
of the real space, the space $\Bm$ of $C^{\infty}-$ functions
bounded in the real infinity by a small constant with all
derivatives. Thus, consider $D^{\alpha}\phi - \mu_{\alpha}
\rightarrow 0$ in the real infinity, for all $\alpha$ and
$\mu_{\alpha}$ constants. Obviously, the space of monotropic functionals $\Bm' \subset
\mathcal{D}_{L^{1}}'$, why $T \in \Bm'$ has representation
$\sum_{\mid \alpha \mid \geq k} D^{\alpha}f_{\alpha}$ with
$f_{\alpha} \in L^{1}$. If $T \in \mathcal{D}_{L^{1}}'$ and $\phi
\in \Bm$, there is a $S \in \Bm'$ such that $S=T$ over $\Bm$. We
have that $\mathbb{R}^{n}=\cup^{\infty}_{j=0} K_{j}$, for compact
sets $K_{j}$. Let $\Phi_{j,1}=(S-T) |_{K_{j}} \in \mathcal{E}'
\subset \Bm'$ and $\Phi_{j,2}=\Phi-\Phi_{j,1}$. We chose $S$ such
that $\Phi_{j,1}=0$ for all $j$ and $\Phi_{j,2} \in \Bm'$. This
gives existence of a functional $S$ such that $S(\phi)=\sum_{\alpha}
f_{\alpha}(x)dx=\lim_{j \rightarrow \infty} T_{j}(\phi)$, where the
limit is taken in $\mathcal{D}_{L^{1}}'$

\vsp

Assume $\Omega_{0} \subset U \subset V$, where $U$ is an open set.
Assume $U$ quasi-porteur for $S \in H'(V)$, that is
$T={}^{t}i(u_{S})$ for $u_{S} \in H'(U)$, where $i$ is the
restriction homomorphism. Assume $r_{T}'$ the transposed ramifier. 
Algebraicity for $r_{T}'$ means that we
can prove that the wave front-set is defined by $b_{\Gamma}$ (\cite{Sjostrand82}) in
$H'$. Assume for the vorticity to the dynamical system $\widehat{w}_{0}$ changes sign finitely many times locally. On
regions where $\widehat{w}_{0}$ has constant sign, we have isolated
singularities in a sufficiently small neighborhood. The lift
function $F$ in $f(\zeta)=F(\gamma)(\zeta)$, can be represented by
$\prod_{p} F_{p}$ relative a division in contingent regions.

\newtheorem{Kohn}[mainresult]{ Proposition }
\begin{Kohn}
Assume $F$ reduced and $F_{T}$ algebraically dependent on $T$. Then
$F_{T}$ is not regularizing.
\end{Kohn}

Proof:\\
We have assumed conditions on the ramifier $\ram$ such that we have
existence of constants $c_{1},c_{2}$ such that $c_{1} \mid \gamma
\mid \leq \mid \ram \gamma \mid \leq c_{2} \mid \gamma \mid$ as
$\mid \gamma \mid \rightarrow \infty$, that is the type $\mid F
\mid_{I}=\limsup_{r \rightarrow \infty} \frac{1}{r} \log \mid
F(\gamma) \mid$ is not dependent on $T$ in the $\mid \gamma \mid-$
infinity and $\mid F_{T} \mid_{I}=\mid F \mid_{I}$. If for this
reason $F$ is not of type $-\infty$, then the same holds for
$F_{T}$. $\Box$.

\newtheorem{rad}[mainresult]{ Proposition }
\begin{rad}
The condition that $F(\gamma_{T})$ is analytically hypoelliptic does
not imply that $\mbox{ Re } F(\gamma_{T})$ or $\mbox{ Im } F(\gamma_{T})$ is
analytically hypoelliptic.
\end{rad}

Assume $P_{T}$ the pseudo-differential operator
that corresponds to $F_{T}$ and that $P_{T}u=f_{T}$ in $H'(V)$, for
an open set $V$, where we are assuming $f_{T}$ holomorphic, that
$\lim_{T \rightarrow 0}f_{T}=f$ in $H'(V)$ and $\lim_{T \rightarrow
0}<u,P_{T} \varphi>=<f,\varphi>$, for $\varphi \in H$. We are
assuming that $P_{T}$ maps $H \rightarrow H$ and that $D(P_{T})$,
the domain for $P_{T}$, has $D(P_{T}) \subset H(V)$. 

\subsection{Paradoxal arguments}

 First note that
among parametrices to partially hypoelliptic differential operators,
considered as Fredholm operators on $L^{2}$, there are examples of
operators with non-trivial kernels. These can be proved to be
hypoelliptic outside the kernel. If they are defined as regularizing
on the kernel, they will not be hypoelliptic there. The class of
partially hypoelliptic differential operators can be shown to be
different from the class of hypoelliptic differential operators on
$L^{2}$. The following argument for $C^{\infty}$-hypoellipticity is
based on two fairly trivial observations,
\begin{itemize}
\item[i)] The Dirac measure $\delta_{0}$ is not a ($C^{\infty}-$)hypoelliptic
operator.
\item[ii)] If $E$ is a parametrix to a differential operator $P$
such that $PE-\delta_{x} \equiv 0$ in $V$, an open set in the real
space (a neighborhood of $x$), then $P$ is not a
($C^{\infty}-$)hypoelliptic operator.
\end{itemize}

Proof of the observations:\\
For the first proposition, define a convolution operator on
$\mathcal{E}'$, $H(\varphi)=E_{0}*\varphi$, where $E_{0}$ is a
fundamental solution with singularities in $0$ to $P(D)$ and where
$P(D)$ is a ($C^{\infty}-$)hypoelliptic differential operator with
constant coefficients. If $\delta_{0}$ were ($C^{\infty}-$)
hypoelliptic over $\mathcal{E}'$, then $\Ss \delta_{0}*\varphi=\Ss
\varphi$ and also $\Ss H({}^{t}P(D) \varphi)= \Ss \varphi$, but
since $E_{0}$ is regularizing outside the origin, $\varphi$ can have
singularities that $H(\varphi)$ does not have.

For the second proposition, we use the notation $I_{E}(\varphi)=\int
E(x,y) \varphi(y) dy$ and $I$ denotes the identity operator, that is
an operator such that $\Ss{ Iu }=\Ss{ u }$ for all $u \in
\mathcal{D}'$. If $P(D)$ were $(C^{\infty}-)$ hypoelliptic, then
$I_{E}-I$ would be locally regularizing. If locally $I_{PE}=I$, we
also have that locally $u-Pu \in C^{\infty}$ for $u \in \mathcal{D}$. But if
$P(D)$ is $(C^{\infty}-)$ hypoelliptic, then the same must hold for
$P-1$ and we have a contradiction. $\Box$

\vsp

 The first observation can immediately be adapted to analytic
hypoellipticity. For the second we note that if $P$ is a
differential operator, then $P-1$ can not be regularizing and the
proof is conclusive also for analytic hypoellipticity. As a
consequence of this, the pseudo differential operators that we are
studying will be assumed locally not regularizing.

\vsp

 The conclusions are as follows. If $f$ is the symbol to an
analytically hypoelliptic pseudo differential operator in the class
that we are studying, we have that all approximations $f_{T}$ can be
chosen regular.  The condition that the dependence of $T$ is
algebraic for $f_{T}$ is necessary to avoid a paradox in the
analogue to Weyl's lemma. It is necessary to have symplecticity on each stratum.
The involution is used to guarantee existence of an inverse lifting
function, since in this case $F_{T}$ can be chosen regular in $T$.

\vsp

For a symbol in $\Bm$ over the real space (modulo regularizing operators), we again
consider (modulo monotropy) locally algebraic symbols. For an constant coefficients polynomial operator a condition equivalent with hypoellipticity is that every distributional parametrix
is very regular. These parametrices map $\mathcal{D}' \rightarrow \mathcal{D}^{'F}$, why it 
is necessary for the pseudo differential operator to be hypoelliptic, that the symbol is of real dominant type (orthogonal real and imaginary parts). The analysis is focused on the microlocal contribution from the lineality. The singular support is considered as a formal
support in a ball of $\epsilon-$ radius..

\vspace*{.5cm}

 Assume temporarily
that the operator is not self-adjoint. Consider $E$, a parametrix to a homogeneously hypoelliptic, constant coefficients operator $P(D)$.
Then, $I_E - I$ is regularizing in $\mathcal{D'}^F$ and thus is represented by a kernel in
$C^{\infty}$, which has a regularizing action in $\mathcal{E'}$. However, it is not trivial to
extend this action to $\mathcal{D'}$. Consider instead $C_{I_E}=I_E \varphi - \varphi I_E$, for some suitable
real function in $C^{\infty}_0$ acting on $\mathcal{D'}$. Since
$$ C_{I_E}f=\varphi f - \varphi I_E f + I_E \varphi f - \varphi f \  f \in \mathcal{D'}$$
This operator will be regularizing in $\mathcal{D'}$.

\section{Lineality}

\subsection{The lineality and the wavefront set}
The lineality $\Omega_{0}$ can be considered as the "boundary" to
the frequency component. More precisely, assume $\Gamma$ a simple
cone in $\Omega_{0}$ and $B_{\Gamma}=\lim_{t \rightarrow
0}A_{\Gamma}$, where $A_{\Gamma}=\FB^{-1} \tau_{\Gamma} \FB \mbox{:}
H(E_{\mathbb{R}}) \rightarrow H'(E_{\mathbb{R}})$, for a real vector
space $E_{\mathbb{R}}$. Assume $h_{F}$ the growth indicator to
$B_{\Gamma}$ and that $g$ is growth indicator for the frequency
component to $WF_{a}(u)$ (cf. \cite{Martineau} Ch.2, Theorem 4.3).
As $h_{F}=g$ on $\Delta_{0}=\Omega_{0} \backslash 0$, we see that
cones in $\Delta_{0}$ have indicator $\geq0$. Let $W$ be the convex
closure of the real support to $B_{\Gamma}$, that is $W=\{y \quad
<y,\eta> \leq h_{F}(\eta) \quad \mid \eta \mid=1 \}$. Let $W_{+}=\{
y \in W \quad <y,\eta> \geq 0 \quad \mid \eta \mid=1 \}$ and let
$W_{-}$ be the complementary set. Let $V_{+}=\{ \eta \quad <y,\eta>
\geq 0 \quad y \in W_{+} \}$, then $\Delta_{0} \cap V_{+}=\{ \eta
\quad <y,\eta>=0 \quad y \in W_{+} \}$. Further, since $g=h_{F}=0$ on
$\Delta_{0}$, we must have $\Sigma \cap \Delta_{0} \cap
V_{-}=\emptyset$.

\subsection{ The lineality is standard complexified }
 We can show that the lineality to a polynomial, is standard complexified in $Exp$,
why it is sufficient to consider completely imaginary translations
of the real space. We shall now see that if we have lineality and if
the lineality is locally algebraic, there is lineality in a complete
disk (cf. \cite{Julia23}). Assume $0$ an essential singularity and
$\Delta$ simply connected and closed $\ni 0$ with boundary $\Gamma'
\cup \Gamma''$. If for a holomorphic function $f$, $\mid f \mid$ is
bounded on $\Delta$ and $f(z) \rightarrow w$ as $\Gamma' \ni z
\rightarrow 0$ and $\Gamma'' \ni z \rightarrow 0$, then $f(z)
\rightarrow w$ uniformly as $z \rightarrow 0$ in $\Delta$.
Conversely, if the limits on $\Gamma',\Gamma''$ are different, then
$f$ can not be bounded on $\Delta$. Assume $\Delta$ with a algebraic
definition locally, then given a sector $A0B$ where $f$ is assumed
holomorphic, if $f \rightarrow w$ as $z \rightarrow 0$ on a line
$0L$ in this sector, the same holds for any sector inner to $A0B$.
Thus, if we have lineality on a line $OL$, we have lineality on the
disk. The conclusion also holds for the several dimensional set of lineality,
but since hypoellipticity can be derived from one dimensional translations, we do not prove this here.

\subsection{Remarks on hypoellipticity and symmetry}
 An operator is considered as
hypoelliptic, if its symbol is reduced in a neighborhood of the
infinity, but for a holomorphic symbol it is not simultaneously
reduced in a neighborhood of the origin. Note also that if $f(z)$ is
reduced as $z \rightarrow \infty$, then $\overline{f(z)}$ is not necessarily
reduced as $z \rightarrow \infty$. If
$f(\overline{z})=\overline{f(z)}$, we have that $f(\overline{z})$ is
not necessarily reduced, as $z \rightarrow \infty$. This property is
consequently not symmetric with respect to the real axes. A necessary condition on a mapping 
$c$ to preserve reducedness, when$f(c(z))=c(f(z))$, is that it is bijective.

\vsp

In this context we consider the property $(P)$ for a continuous function $d$,
that is $d(\frac{1}{T})=\frac{1}{d(T)}$ as $T \rightarrow \infty$. For
instance if $d$ is the distance function to the boundary, if there
is no essential singularity in the infinity and if all singularities
are isolated in the finite plane, then $d$ is globally reduced and
$d$ has the property $(P)$.

\subsection{The property (P)}

Assume again that $f=e^{\varphi}$ with $\varphi=e^{\alpha}$ and $L(e^{\varphi})=\widehat{L}(\varphi)=e^{\tilde{L}(\varphi)}$ and if $\widetilde{\widetilde{L}}(-\alpha)=-\widetilde{\widetilde{L}}(\alpha)$, we have $L(e^{\frac{1}{\varphi}})=e^{\frac{1}{\tilde{L}(\varphi)}}$, which we denote property $\log (P)$. If $L$ is algebraic, we have that it has the property $\log (P)$. The property (P) means that $\tilde{L}(\varphi)+\tilde{L}(\frac{1}{\varphi}) \sim 0$ and if we have both the properties, we get $\tilde{L}^{-1}=- \tilde{L}$. If we assume $\tilde{L} \sim_{m} \tilde{W}$, where $W$ is algebraic in the infinity, then $e^{\tilde{W}(-\varphi)}=W^{-1}(e^{\varphi})$. We will consider $\tilde{L} \rightarrow 0$, such that we have existence of an algebraic morphism $\tilde{W}$ such that $\tilde{L} \sim_{m} \tilde{W}$ where $\tilde{W}$ has the property (P). We assume existence of $L^{-1}$ over an involutive set where we have a regular approximation. If $\varphi$ is a holomorphic function with $\varphi(\zeta_{T})=\varphi_{T}(\zeta)$ and $\zeta_{T} \rightarrow \zeta_{0}$, as $T \rightarrow \infty$, then using Weierstrass theorem, we have existence of $s$ continuous, such that $s(\varphi_{T}+a)=\zeta_{T}$, for a constant $a$ and $s(a)=\zeta_{0}$. Further, $s$ can be approximated by polynomials of $1/(\varphi_{T}+a)$.

\subsection{Lineality and the characteristic set}

Treves' conjecture is given for the characteristic set $\Sigma$ and our argument is given for the set of lineality. We will now argue that the conjecture can be derived from our result.
Assume $\Sigma=\{ \zeta \quad f(\zeta)=0 \}$ and $\Sigma=\Sigma_{1} \cap \Sigma_{2}$, where $\Sigma_{1}=\{ \mbox{ Re }f(\zeta)=0 \}$. Thus, if $\mbox{ Im }f$ is algebraic, we have that $\Sigma_{2}$ is removable. The condition $\mbox{ Re }f \bot \mbox{ Im }f$ is considered as necessary for hypoellipticity. We note in this connection the
well-known Weyl's lemma (cf. \cite{AhlforsSario60}), if $w \in
L^{2}(\mid z \mid <1)$ and for all $V \in C^{2}_{0}(\mid z \mid<1)$,
we have $(w,dV)=(w,dV^{\diamondsuit})=0$, (harmonic conjugation) then $w$ is equivalent
with a $C^{1}$ form. 

\vsp

Assume $(I)=(\mbox{ ker }h)$, where $h$ is a homomorphism and assume existence of a homomorphism $g$, such that $d h(f)=g(f) dz$. If $g$ is algebraic and $g^{-1}(0)=const$, we can define $\Delta$ as semi-algebraic. Note however that existence for a global base for $g$, does not imply existence of a global base for $h$. Let $C_{1}=\{ f=c \}$ and assume $\Delta=V_{1} \backslash C_{1}$, where $V_{1}=\{ f_{1}=0 \}$ and $f_{1}=\tau f -f$. Let $\Delta_{0}$ be $\Delta \backslash \{ x_{0} \}$, where $x_{0}$ is the intersection point. We can choose $g(f_{1})=0$ on $C_{1}$ and $g(f_{1}) \neq 0$ on $\Delta$. Note that if $\Delta \cup C_{1}=V_{1}$, then $I(V_{1})=I(\Delta)I(C_{1})$. Assume $\Delta_{0} \cap C_{1}= \emptyset$, then $g \in I(V_{1})$ implies $g=pq$, where $p \in I(\Delta)$ and $q \in I(C_{1})$.  Assuming $C_{1}$ oriented, we can choose $\Delta$ as locally algebraic $p_{+}p_{-}q=g$, where $p_{\pm}$, have one-sided zero-sets. If we assume $C_{2}=\{ f=c \quad \frac{d}{d T}f=c' \}$
and $I_{2}$ the ideal of non-constant functions and $NI_{2}=V_{1} \cup V_{2}$, then $INI_{2} \sim \mbox{ rad }I_{2}$. If $V_{1} \cap V_{2}=\emptyset$, we can write $g=g_{1}g_{2} \in \mbox{ rad }I_{2}$.  Assume $\rho$ a measure such that $\rho_{1}(I(\Sigma))=\rho(I(\Sigma_{1}))$ and correspondingly for $\rho_{2}$.If $\Sigma_{1} \cap \Sigma_{2}=\emptyset$, then the measures can not be absolute continuous with respect to each other. If we instead consider two ideals of analytic functions
$I_{1}=\{f \quad dh(f)=0 \}$ and $I_{2}=\{f \quad f=const. \}$ and the corresponding measures $\rho_{1}(I)=\rho(I_{1})$,$\rho_{2}(I)=\rho(I_{2})$. Then if
$0=\rho(I_{1}^{c})$ implies $\rho(I_{2}^{c})=0$, we have $\rho_{2}$ is absolute continuous with respect to $\rho_{1}$. Thus, we have existence of $f_{0}$ Baire (cf. \cite{Garding55}), such that $\rho_{1}(f_{0}f)=\rho_{2}(f)$ and $f \in L^{1}(\rho_{1})$.

\newtheorem{flip}{ Proposition }[section]
\begin{flip}
Given an analytic symbol with first surfaces $C$, the lineality can be studied locally as transversals. Conversely, given the lineality and a normal model, the lineality approximates the first surfaces to the symbol.
\end{flip}

Existence of lineality can be seen as a proposition of possibility to continue the symbol on a set of infinite order, that is the symbol is not reduced with respect to analytic continuation.  Assume for a measure $\rho$, $\rho(\Top \varphi)=\rho(\varphi^{*})$, for $\varphi \in L^{1}$, on an algebraic set and $\rho(\Top \varphi)=0$ implies $\rho(\varphi^{*})=0$, then we have existence of $\varphi_{0}$ Baire such that $\rho(\Top \varphi)=\rho(\varphi_{0} \varphi^{*})$

\vsp

 We know that (cf. \cite{Poincare87}) every form
$\sum_{j} B_{j} dx_{j}$ invariant relative closed contours, has the
representation $\int \sum_{j} B_{j} d x_{j}=\int dW+ \int \sum_{j}
c_{j} dx_{j}$, where $dW$ is exact and the last integral is an
absolute invariant. The argument can be repeated for our ramifier
and $\int dV=\int B(dx_{T},dy)-B(dx,dy)$ with $V(x,y)=W(\ram
x,y)-W(x,y)$ and $dW$ exact. We have assumed that the ramifier is close to translation and we have the following
explanation of this. Assume $\sum_{j} F_{j}d x_{j}$ invariant in the sense that $\sum_{j}
\int F_{j}(r_{T}' x)dx_{j}=\int F_{j}(x)dx_{j}$ and assume that $\tau_{\Gamma}$ is translation. Let
$$ dK_{T}=\sum_{j} \big[ F_{j}(x)(r_{T} \zeta)-F_{j}(x)(\zeta)
\big]dx_{j}$$
$$ dL_{T}=\sum_{j} \big[F_{j}(x)(\tau_{\Gamma}
\zeta)-F_{j}(x)(\zeta) \big] dx_{j}$$ We can prove
 that over regular approximations, we have that $\int dK_{T} +\int
 \sum_{j} C_{T,j} dx_{j} \sim \int dL_{T} + \int C_{T,j}' dx_{j}$,
 for constants $C_{T,j},C_{T,j}'$.

\section{Involution}

\subsection{Introduction}

Given a multivalued surface, a canonical approximation is the spiral Puiseux approximation,
but some results require a tangent determined, why we prefer transversal approximations. Sufficient conditions for existence of transversals are discussed in connection with the lifting principle.

We note that assuming polynomial right hand sides, for the associated dynamical system, monotropy is microlocally invariant. That is since a bounded set can not contribute as lineality, obviously $\epsilon$ translation does not affect this proposition.
In this case monotropy (cf. \cite{Cousin95}) correponds exactly to adding a small constant (the value in the origin to a polynomial) to the symbol in the infinity. For analytic right hand sides, the two monotropy concepts are no longer equivalent, but the microlocal invariance can be proved for both separately.

\subsection{Tr\`{e}ves curves}
Assume $<,>_{1}=\mbox{ Re } < y_{T}^{*},y_{T}>-1$ for $\gamma_{T} \in
\Gamma$. If $\Gamma_{T}$ describes a line, we
have that $\gamma_{T}^{*}$ describes a line. Let $\Sigma=\{ \gamma_{T}
\quad \frac{d}{dT}F_{T}=F_{T} \}$ and $\Sigma_{0}=\{ \gamma_{T}
\quad \frac{d}{dT}F_{T}=F_{T}=0 \}$. Let $F_{T}'$ be the transposed operator
to $F_{T}$ with respect to $\frac{d}{dT}$, that is $F_{T}'
\frac{d}{dT}=\frac{d}{dT}F_{T}$, why on $\Sigma$,
$F(\frac{d}{dT}\gamma_{T}-\gamma_{T})=0 \Leftrightarrow
\frac{d}{dT}(F_{T}'-F_{T})=0$.  Let $A=T \Sigma=\{ \gamma_{T} \quad
\frac{d^{2}}{dT^{2}}F_{T}=\frac{d}{dT}F_{T} \}$.  If for every $\theta_{T} \in T \Sigma$, we have $<
\frac{d}{dT} \gamma_{T}, \theta_{T}>_{1}=0$, we have that
$\frac{d}{dT} \gamma_{T} \in \mbox{ bd }A$. Further, $\frac{d
F_{T}}{d \gamma_{T}} \frac{d \gamma_{T}}{d T}=F_{T}' \frac{d
\gamma_{T}}{dT}$ over $A$, why $F_{T}'=d F_{T}/d \gamma_{T}$. Thus,
for instance $F_{T}'(\gamma_{0})=0$, where $F_{T}(\gamma_{0})$ is
constant. If $F_{T}'$ maps $A^{\bot} \rightarrow A^{\bot}$ over
$\Sigma$, we have that $< \frac{d}{dT}F_{T}(\gamma),
\theta_{T}>_{1}=0$, so $\frac{d}{dT}F_{T}(\gamma) \in A^{\bot}$ and
$< F_{T}(\gamma),\theta_{T}>_{1}=0$. A sufficient condition for
$F_{T}$ to map $A^{\bot} \rightarrow A^{\bot}$ is that
$<F_{T}(\gamma),\theta_{T}>_{1}=\rho_{T}<\gamma_{T},\theta_{T}>_{1}$,
where $\rho_{T}$ is a function, not involving any differentials (a
multiplier).  The 
proposition, is that $\gamma_{T}^{*} \bot (\mbox{ bd }A) \Rightarrow
\gamma_{T}^{*} \bot A$, which can be fulfilled if $A$ is on one side
 locally of a hyperplane. If the symbol ideal is symmetric and
finitely generated over a pseudo convex domain, this can be assumed.

\vsp

Assume $\Phi \bot \mbox{ bd } \Sigma=\{ F_{T}(\eta)=c \}$, for a
constant $c$, implies $\Phi \bot \{ \eta \geq c_{1} \}$ (a semi
algebraic characteristic set $= \Sigma$). If
$<F_{T}(\eta),\phi>=C_{T}<\eta,\phi>=0$. Assume further that
$F_{T}(\eta)=c \Leftrightarrow \eta=c_{1}$, for constants $c,c_{1}$,
why $F_{T}$ maps $\Sigma \rightarrow \Sigma$. We know that if
$\eta_{T}^{*}=y_{T}^{*}/x^{*}$ with $x^{*},y_{T}^{*}$ polynomials
and $\int_{\Sigma} \eta_{T}^{*} dx dx^{*}=0$, then $\Sigma$ has
measure zero. In the same manner for $\eta_{T}$. Assume $\mbox{ bd
}\Sigma=$ $\{$ the set where $\eta$ changes sign $\}$ and where
$\Sigma$ has measure $>0$. If we have existence of $\gamma \bot
\Sigma$ holomorphic, we must $\gamma \equiv 0$ on $\Sigma$, through
Hurwitz theorem and we conclude that there can not exist an
algebraic $\gamma \bot \Sigma$ with these conditions. The conclusion
is that if $\gamma \bot \mbox{ bd }\Sigma$, we can not have, for
$\gamma$ algebraic, that $\Sigma$ stays locally on one side of a
hyperplane. More precisely, if there are $2m$ characteristics
through a singular point (cf. \cite{Bendixsson01}), where $m$ is
referring to the order of $X,Y$ in the associated dynamical system,
and if the sign is changed passing the characteristics, then the set
of for instance positive sign is not separated by a hyperplane. By
giving the characteristics a direction however, the problem can be
handled. Assume $\Sigma_{+}=\Sigma_{1} \cup \Sigma_{3}$, the domain
for positive sign and that $\eta$ is an algebraic characteristic
with $\int_{\Sigma_{+}} \eta dx dy=0$, then $\eta=0$ either on
$\Sigma_{1}$ or $\Sigma_{3}$, depending on which direction $\eta$
has. We can thus have half-characteristics $\eta \bot \Sigma_{1}
\cup \Sigma_{3}$ with algebraic definition.

\vsp

$F_{T}$ is said to be reduced for involution, if given existence of
a regular approximation $G_{T}$ in $(I)$ with $(I)=\mbox{ ker
}H_{V}$, we have $\{F_{T},G_{T} \}=0$ on $S_{T}$ implies $T=0$. In
this case there are no level surfaces to $F_{T}$ on $S_{T}$. Over reduced $x$
, we have that $\ram x=x$ implies $\ram-id$ is only
locally algebraic. If $\ram$ is algebraic in $T$ with minimally
defined singularities, then $\ram -id \sim_{0}$ a polynomial.

\newtheorem{bdcond}[mainresult]{ Boundary condition}
\begin{bdcond} \label{bdcond}
The boundary is characterized by the condition that $F_{T}$ is
holomorphic in $T$, for $T \notin \Sigma$ or $\frac{d F_{T}}{d T}$
holomorphic in $T$, for $T \notin \Sigma$, where $\Sigma$ is given
by locally isolated points and the regularity is close to the
boundary.
\end{bdcond}
 More precisely, let $\Sigma=\{ \zeta_{T} \quad F_{T}=const.\quad
\frac{d F_{T}}{d T}=const. \}$ and as previously $(I_{1})=\{
\gamma_{T} \quad F(\gamma_{T}) \text{ is not constant } \}$, where
$F_{T}$ is assumed holomorphically dependent on the one dimensional
parameter $T$. Let $(I_{2})=((I_{1}) \cap (\frac{d}{dT}(I_{1}))$ and
$N(I_{2})=V_{1} \cup V_{2}$, where $V_{1}=\{ \zeta_{T} \quad F_{T}
\text{ is not constant } \}$ and $V_{2}=\{ \zeta_{T} \quad
\frac{d}{dT}F_{T} \text{ is not constant }\}$. Using the
Nullstellensatz, we can form $IN(I_{2}) \sim rad(I_{2})$ and we
claim that $(I_{2})$ is radical. The condition can be generalized to
higher order derivatives.
\newtheorem{Treves}[mainresult]{ Lemma }
\begin{Treves}
The condition that $F_{T}$ is not reduced for involution means that
there exist Tr\`{e}ves curves in $S_{T}$.
\end{Treves}

Proof:\\
 Assume for this reason that $T \neq 0$ and that there
exist $\gamma_{T} \subset S_{T}$ such that $\{ G_{T},F_{T} \}=0$
over $\gamma_{T}$, where $F_{T}$ is a lifting function and $G_{T}$
is a regular approximation of a singular point. Assume
$G_{T}(\gamma)=G(\gamma_{T})$ and
$\frac{d}{dT}G(\gamma_{T})=G_{1}(\frac{d \gamma_{T}}{dT})$ with
$G_{1}$ invertible over $\frac{d \gamma_{T}}{dT}$. Assume existence
of $v_{T}$, a regular approximation of $T \Sigma$ with $\frac{d
v_{T}}{dT}=G_{1}\frac{d \gamma_{T}}{dT}$, then $< \frac{d
v_{T}}{dT},\theta_{T}>_{1}=0$ and if $G_{1} \text{:} A^{\bot}
\rightarrow A^{\bot}$ (independence of $T$ at the boundary), we see
that there exist Tr\`{e}ves curves for $F_{T}$ in $S_{T}$.$\Box$

\subsection{ First surfaces }

Consider the system $\frac{dx}{X}=\frac{dy}{Y}=dt$ and the
corresponding variation equations
$\frac{dx^{*}}{dt}=\frac{dX}{dx}x^{*}+\frac{dX}{dy}y^{*}$ and
$\frac{dy^{*}}{dt}=\frac{dY}{dx}x^{*}+\frac{dY}{dy}y^{*}$. Assume
$F_{T}(x,y,x^{*},y^{*})$ a first integral to the variation
equations, algebraic in $x,y$ and homogeneous of order $1$ in
$x^{*},y^{*}$. It is well known that $\int F_{T}(dx,dy)$ is
invariant integral to the given system. Conversely, if $\int
F_{T}(dx,dy)$ is invariant integral to the system, then
$F_{T}(x^{*},y^{*})$ is integral to the variation equation. Assume
$V$ the Hamilton function to the system, that is
$\frac{dx}{dt}=\frac{dV}{dx^{*}}$,$\frac{dx^{*}}{dt}=\frac{-dV}{dx}$
and
$\frac{dy}{dt}=\frac{dV}{dy^{*}}$,$\frac{dy^{*}}{dt}=\frac{-dV}{dy}$.
Then given that $\{ V,V_{1} \}=0$, also $V_{1}$ is a Hamilton
function. If $V_{2}$ is a Hamilton function, we have that $\{
V,V_{2} \}=0$ and $\{ V_{1},V_{2} \}=0$. We will consider an
involutive set $S_{T}$ such that for $F_{T}$ a lifting function and
$V$ a Hamilton function, $\{ V,F_{T} \}=0$ over $S_{T}$. One of the
most important problems in this approach seems to be existence of an
inverse for $F_{T}$. A sufficient condition is reducedness, but this
is not suitable in connection with invariant integrals.

\vsp

 Assume $V$ a Hamilton function and $F_{T}$ a lifting
function to the system $\{ \gamma_{T} \}$ corresponding to the
symbol. Further that $G_{T}$ is a regular approximation (with respect to $\frac{d}{dT})$ to the
singularity in $\{ \gamma_{T} \}$, not necessarily a lifting
function. As $\{V,\cdot \}=H_{V}$ defines an ideal $(I)$, we note
that if $F_{T} \in (I)=(I)(S_{T})$ with
\begin{equation} \label{Hamilton}
\frac{d F_{T}}{d x^{*}}=\frac{d V}{d x^{*}} \text{  and  }\frac{d
F_{T}}{d y^{*}}=\frac{d V}{d y^{*}}
\end{equation}
and if $G_{T} \in (I)$, we have that $\frac{d F_{T}}{dT}=\frac{d
G_{T}}{dT}$. We see that $F_{T}$ is regular under these conditions.
The proposition is that existence of $G_{T}$ regular in $(I)(S_{T})$
and $S_{T}$ involutive means that the lifting function with
(\ref{Hamilton}) is regular.

\subsection{Continuation of the representation}
 Assume
$W \subset V \subset V'$ and $\varLambda$ complex varieties and consider
the mapping $r^{\bot} \text{: }V^{\bot} \rightarrow \varLambda^{\bot}$
with $\mbox{ ker } r^{\bot}=W^{\bot}$. We then have, given $T \in
H'(\varLambda)$, there is a $U \in H'(V')$ with $\FB(U)=\FB(T)$ if and
only if $\FB(U)$ is constant on $W^{\bot}$. Particularly, if
$\FB(T)$ has isolated singularities in the infinity, there is a
continuation principle through the projection method.

\vsp

Given a finitely generated
system with polynomial right hand sides $P,Q$. If the
constant surface corresponding to $P/Q$ is $=\{ 0 \}$, then $L_{T}$
is reduced with respect to contraction, that is $$
-\frac{dL_{T}}{dy} / \frac{dL_{T}}{dx}=\frac{dx}{dy} \quad
\Rightarrow T=0$$. Particularly, consider $F$ over $\gamma$ with right hand sides
$P,Q$ and $df_{T}=dF(\zeta+T)-dF(\zeta)$ for $\zeta \in \Omega$.
Then, over the lineality for $Q/P$, $\frac{d
f_{T}}{dx}/\frac{df_{T}}{dy}=dy/dx=Q/P$. In the same manner, if $G$
is a different form to the same system and $dg_{T}$ as above, if $\{
f_{T},g_{T} \}=0$ over $V$, then $df_{T}=dg_{T}$ over $V$. If
further $\{ f_{T},g_{T} \} \in \overline{I(V)}$, we have
$df_{T}=dg_{T}+c_{T}$ on $V$, for a constant $c_{T}$. Note that
there is an ideal $J$ such that $rad J \sim_{m} \overline{I(V)}$. If
$f_{\Gamma}=F(x+i\Gamma,y)-F(x,y+i\Gamma)$ and
$g_{\Gamma}=f_{\Gamma}^{*}$, we have for an involutive set, that
$df_{\Gamma}=dg_{\Gamma}$. If we have $f_{\Gamma}=g_{\Gamma}$ in
$L^{2} \cap H$ we know that $g_{\Gamma}={}^{t}f_{\Gamma}$.
Obviously, $H_{P}$ defines a functional in $H'$. If $\int_{\Sigma}
H_{P}(f_{T}) dx dy=0$, we have either $H_{P}(f_{T}) \equiv 0$ on
$\Sigma$ or $\int_{\Sigma} dx dy=0$.

 \subsection{ Continuous ramification }

We are assuming the ramifier defines a regular covering (\cite{Grauert60}), that is we
are assuming $\Psi : (I)(\Omega) \rightarrow (I_{T})(\Omega)$, where
the first is a Hausdorff space, $\Psi$ is continuous, proper and
almost injective (singular points are mapped on to a discrete set (subset of transversals) in
$\Omega$). We write $r_{T}'$ for $\Psi$, $N : (I)(\Omega)
\rightarrow \Omega$ and the ramifier is the lift $\Omega \rightarrow
(I_{T})(\Omega)$, such that $r_{T}'I(\Omega)=I(r_{T} \Omega)$.
Denote the critical points to $r_{T}$ with $A$ and assume that they
are parabolic. We assume $\Psi$ such that $I(r_{T} A)$ is nowhere
dense in $(I)(\Omega)$ and so that $\Psi$ is locally a homeomorphism
outside critical points. Finally, we are assuming that for all
$\gamma \in (I)(\Omega)$, there is a small neighborhood
$U_{\gamma}$, open and arc-wise connected, such that the
$U_{\gamma}-I(r_{T}A)$  is arc-wise connected. Wherever $\Psi$ is
holomorphically dependent on the parameter, the inverse $r_{T}$ will
be assumed continuous outside a discrete set. If for instance $A=\{
\zeta \quad d_{\zeta}f(\zeta)=f(\zeta)=0 \}$, we are studying points
$\zeta_{T}$ that can be used to reach $A$ from $\{ f=c \}$, for a
constant $c$.

\vsp
Assume $\Omega_{0} \subset \Omega$, where $\Omega$ is assumed a
pseudo convex domain. Assume $U$ an open set such that $\Omega_{0}
\subset U \subset \Omega$. Assume $T(=B_{\Gamma})$ an analytic
functional, $T \in H'(\Omega)$, quasi portable by $\Omega_{0}$, that
is we have existence of $u \in H'(U)$ with $T=i_{U,V}(u)$
(restriction homomorphism). Let $\Omega_{0} \subset \cup^{N}_{j=1}
U_{j}$, for open sets $U_{j} \subset \Omega$ and $T=\sum_{j=1}^{N}
T_{j}$, that is we can write $T_{j}=i_{U_{j},\Omega}(u_{j})$ with
$u_{j} \in H'(U_{j})$. Assume now the restriction homomorphism
algebraic, then we have if $\Omega_{0}$ is complex analytic in a
real analytic vector space, that $T$ is portable by $\Omega_{0}$
(cf. \cite{Martineau} Ch.2, Section 2).

\vsp

Assume $h$ algebraic and let $v_{T}(x)=h(r_{T}'x)-h(x)$, where
$r_{T}'$ is a continuous linear mapping $d(r_{T}' f)/dx=d f/dx$ and
we write $f(r_{T} \zeta)=r_{T}' f(\zeta)$. Let $\Delta=\{ T \quad
f(r_{T} \zeta)=f(\zeta) \quad \forall \zeta \}$. Over an ideal,
finitely generated and of Schwartz type topology with (weakly)
compact translation (cf. \cite{Martineau}), there are given $T_{0}
\in \Delta$, $T_{j}$ regular such that $T_{j} \rightarrow T_{0}$ and
$f(r_{T_{j}}) \zeta)=f(\zeta)+C_{T_{j}}$, for constants $0 \neq
C_{T_{j}} \rightarrow 0$ as $j \rightarrow \infty$. The sets $\{
v_{T}=0 \}$ will not contribute micro locally, however the sets $\{
v_{T}=const \}$ contributes to invariance in the tangent space and
gives a micro local contribution.

\vsp

Assume $L$ an analytic line, transversal in a first surface $S_{0}$
through $p_{0}$ and consider a neighborhood $\Gamma$ of $p_{0}$ on
$L$. Denote $\Sigma_{\Gamma}$ the set of points that can be joined
with a point in $\Gamma$, through a first surface to $f$. We assume
$L$ transversal to every first surface through $\Gamma$ of order
$1$. Transversality means existence of regular approximations. We
will in this approach not assume minimally defined singularities. If
for a first surface $S$, we have $S \cap \overline{\Sigma}_{\Gamma}
\neq \emptyset$, we have $S' \subset \overline{\Sigma}_{\Gamma}$,
for all $S' \sim S$ (conjugated in the sense of \cite{Nishino68}).
Thus for a generalization of the inhomogeneous Hange's result, it is
sufficient to consider the normal tube. That is if $\Gamma$ gives a
micro local contribution in $p_{0}$, then if $S^{\bot}$
(transversal) has $S^{\bot} \cap \Sigma_{1}(u) = \emptyset$, we have
$S \cap \overline{\Sigma}_{\Gamma} \neq \emptyset$ implies $S
\subset \overline{\Sigma}_{\Gamma}$, so that $\Gamma \in S^{\bot}$.
Note however that it is necessary for micro local contribution, that
the set is not bounded globally.

\subsection{The condition on involution}

First a few notes on the lifting principle. Assume $\gamma \in
\mathcal{P}$, an analytic polyeder. It is not true that the lifting
principle holds over every $\mathcal{P}$, but by constructing a
normal model $\Sigma$ (\cite{Nishino96}) to $\mathcal{P}$, we have
always (modulo monotropy) a lifting function. Let $\Omega=\{ \zeta
\quad r_{T}' \gamma(\zeta) \in \Sigma \quad \gamma \in \mathcal{P}
\}$, by the definition of the ramifier $r_{T}$, $\Omega=\{ r_ {T}
\zeta \quad \gamma(r_{T} \zeta) \in \Sigma \quad \gamma \in
\mathcal{P} \}$. We assume $\gamma_{T}$ (real-) analytic on $V
\times \Omega \ni (T,\zeta)$. For $\zeta$ fix in a neighborhood
defined by $T$, $F$ can be chosen holomorphic. Let $\mathcal{P}=\{
\gamma(\zeta) \quad \gamma \text{ holomorphic in } \zeta \in \Omega
\}$. Then, if we assume $\mathcal{P}$ finitely generated over
$\Omega$ and $r_{T}' \mathcal{P}=\Sigma$, we get a corresponding
$\tilde{\Omega}=\{ r_{T} \Omega \}$ and $f(r_{T}
\zeta)=F(\gamma_{T})(\zeta)$ for $\zeta$ fix, can be extended to the
domain for $f$, in a neighborhood of a first surface. Thus, the
construction is such that $\Omega$ is a neighborhood of $\{ T=0 \}$
and $\zeta$ in a first surface, why we have existence locally of a
lifting function for a normal model.

 \vsp

 The condition on
involution gives existence of the inverse lifting function
$G_{T}=F_{T}^{-1}$. We are now interested in determining the domain
where $G_{T}$ is constant, algebraic, holomorphic etc. Note that if
$G_{T}(f)$ is algebraic in $f$ and $f$ is the symbol to a
hypoelliptic operator, then in the real space, $G_{T}(f) \in \Bm$.
Assume existence of $G'$, derivative with respect to argument, then
from the regularity conditions for the dynamical system, $G_{T}(f)$
has isolated singularities and if $G'$ holomorphic or constant, we
must have isolated singularities for the symbol $f_{T}$. Consider
$(I)_{+}=\{ \gamma_{T} \quad F(\gamma_{T})=F(\gamma) \quad \mbox{ Im } T >0
\}$ and correspondingly $(I)_{-}$. Assume $F_{T}$ algebraic in $T$,
then the signs will give an orientation to the first surfaces. Thus,
$(I)_{+}$ will correspond to conjugate classes of first surfaces
(\cite{Nishino62}). For instance in case
$F(\gamma_{T})=F(\gamma_{\overline{T}})$, we have the same first
surface in $(I)_{\pm}$ but different orientations. We will assume
the number of classes constant, when $\mbox{ Im } T$ is small (compare with
the regularity conditions (\cite{Bendixsson01})). The regularity for
$G_{T}$ will now determine the character of the first surfaces.
Regular first surfaces, for instance have only trivial conjugates,
which will be the case if $G_{T}$ is reduced. We have noted that all
normal approximations can be chosen regular. 

\vsp

Consider the symbol $F=PF_{0}$ with $P$ a polynomial,
$F_{0}=\widehat{f_{0}}$ and $F_{0}$ holomorphic or monotropic with a
holomorphic function. Let
$$ (I)_{\Omega}=\{ f_{0} \in (\Bm)' \quad P \widehat{f_{0}}=0
\text{ on } \Omega \}$$ If $\Lambda=Z_{P}$ (zero-set), we have that
$F_{0} \bot \Lambda$ implies $f_{0} \in (I)_{\Lambda}$. Conversely, if
the polynomial $P$ is reduced and $\mid PF_{0} \mid < \epsilon$ at
the boundary for a small number $\epsilon$, then the Nullstellensatz
(\cite{Oka60}) gives that $F_{0}$ is bounded by a small number at
the boundary.

\subsection{ The lifting principle }
Assume the right hand sides to the associated dynamical system $X,Y$
are polynomials in $\zeta$, then according to the lifting principle
(cf. \cite{Oka60}), we have on $\mid X \mid \leq 1$,$\mid Y \mid
\leq 1$, existence of a function $F$ holomorphic in $x,y$, such that
$f(\zeta)=F(x,y)(\zeta)$. If $\zeta$ is in a polynomially convex and
compact set, $f$ can be represented as a polynomial. Assume
$\vartheta=Y/X$ and $\eta=y/x$, for polynomials $X,Y$. Further, for
constants, $c,c'$, $\mid \vartheta -\eta \mid >c$ and $\mid \eta
\mid<c' \mid \eta \vartheta -1 \mid$ locally. We can determine $w$
algebraic and locally maximal, such that $\mid w \vartheta -1 \mid <
1$. For $\eta \vartheta \sim_{m} w \vartheta$, we have existence of
a holomorphic function $F$, such that $F(\eta)(\zeta) \sim_{m}
f(\zeta)$ and $F(\eta)=const.$ $\Leftrightarrow \eta=const.$ If $F$
is invariant for monotropy, the result $F(\eta)=f$ follows directly
from the lifting principle. Assume $\mathcal{P}$ an analytic
polyeder with separation condition (cf. \cite{Oka60}),
$\mathcal{P}=(x,h(x))$. Assume $\overline{\Delta}=\{ \mid z_{j} \mid
\leq 1 \}$ and $\Delta^{\epsilon}=\{ \mid z_{j} \mid \leq 1+
\epsilon \}$, close to $\overline{\Delta}$, for $j=1,2$. Assume
$\Sigma=\Phi(\mathcal{P})$, such that $\Phi(\delta \mathcal{P})
\subset \delta \mathcal{P}$ (conformal) and that $\Sigma$ is an
analytic set with continuation  in $\Delta^{\epsilon}$. Then,
$\Sigma$ is a (normal) model for $\mathcal{P}$. Assume $f$ analytic
on $\mathcal{P}$, then we have existence of $F$ holomorphic on
$\overline{\Delta}$, such that
$f(\zeta)=F(\Phi(\mathcal{P}))(\zeta)$. Note that we are assuming
$\zeta$ in a symmetric neighborhood of $\{T=0 \}$. We can, according
to Rouch\'{e}'s principle assume, $\mid z_{j}-w_{j} \mid < \epsilon
\mid z_{j} \mid$ and $\mid w_{j} \mid \leq 1$, for $z_{j} \sim_{m}
w_{j}$. For $w_{j}$ polynomials, this is a proposition on $F$ being
invariant for monotropy.

\vsp

\subsection{ Exactness and involution }

We will use the Poisson bracket ${U,V}=\sum_{i} \frac{\delta
U}{\delta x_{i}} \frac{\delta V}{\delta y_{i}}-\frac{\delta
U}{\delta y_{i}} \frac{\delta V}{\delta x_{i}}$. Assume $V$ defined
through $\frac{\delta V}{\delta y_{i}}=\frac{d x_{i}}{d
t}$,$\frac{\delta V}{\delta x_{i}}=-\frac{d y_{i}}{d t}$. Concerning
the two possibilities for $(x_{i},y_{i})$ where $i=1,2$,
A) $(x,y,x^{*},y^{*})$
B) $(x,x^{*},y,y^{*})$
,it does not appear to be important what
representation we use.

\vsp

Consider the sets $\Phi_{\vartheta}=\{ e^{\vartheta} M=W \}$ and
analogously for $\Phi^{*}$. Thus,
$F^{\diamondsuit}(M)=F(e^{\vartheta} M)$. Assume over an involutive
set that $\exists F^{-1}$ and let $G=F^{-1} F^{\diamondsuit}$ over
$M$. Then, $G(M)/M=e^{\pm \varphi}$. We will study the parabolic
sets $\pm \varphi < 0$, so that $G(M)=const. M$. The spectrum is $\{
e^{\vartheta} M=W \}$, then for a lifting operator $F$, invertible
and over $\vartheta < 0$, we have $F^{\diamondsuit}(M)=const. F(M)$,
if the constant is real, we have real eigenvectors. There will be a
boundary in this approach, given by the set where $\vartheta$
changes sign. Finally, we consider the sets where $\vartheta > 0$
(real and holomorphic). If the underlying sets in $\Omega$ are
simply connected, these sets constitute neighborhoods of the
constant surfaces. If we consider $F$ as an analytic functional, we
have that $F$ has the closure of $\{ \vartheta < \alpha \}$ as
semi-porteur if and only if the type for $\widehat{F}$ is $\leq
\alpha$, which particularly means that it is portable by any convex
neighborhood of the semi-porteur.

\subsection{ Dependence of parameter }
Given a closed trajectory, that does not end in a singular point
$P$, that is the point $P$ stays inner to the trajectory. The point
$P$ is called a center, if there are infinitely many closed
trajectories, arbitrarily close to $P$, that circumscribes the
point. We could say that the trajectory $\gamma_{T} \rightarrow
\gamma_{0}=P$, but does not reach it.  We will assuming the boundary
not $C^{1}$, but holomorphic and with only parabolic singularities,
consider the problem of removing the center point as a Dirichlet
problem.

\vsp

There are certain conclusions on the singularities in
$\Omega_{\zeta}$, given the dependence of the parameter in $L^{1}$.
We have the following weak form of minimally defined singularities.
For $F=w \in \Bm'$, if $x,y \in \Bm$ and $\int_{I}w_{T}(x,y) d \sigma
\rightarrow \int_{I}w(x,y) d \sigma$, through a normal and regular
approximation. Assume that the dependence of $T$ is holomorphic and
$w$ algebraic in $(x,y)$. We have that $\{ (x,y) \quad
w_{T}(x,y)=w(x,y) \}$ has $\sigma-$ measure zero. Assume that $\int
\mid w \mid^{2} d \sigma < \infty$ and $w w^{*}=w^{*}w$ and that
$(x,y)$ is in the normal tube. Then we have normal and regular
approximations, say $g_{T}$ of $\{ x=const.,y=const. \}$. Assume
$f_{T} \rightarrow 0$, as $T \rightarrow 0$ normally and regularly,
such that $\frac{d f_{T}}{d T}$ is holomorphic in $T$ (that is not a
non-zero constant). If $\frac{d f_{T}}{d T}=C \frac{d g_{T}}{d T}$,
for a constant $C$ on a domain of positive measure, we still have a
regular approximation. If $\gamma_{T}=(f_{T},h(f_{T})$ is the
regular approximation and if the dependence of $T$ is algebraic in
$\frac{d \gamma_{T}}{d T}$, then according to Hurwitz theorem, since
polynomials never have zero-sets of infinite order, then the
zero-set must have measure zero. Thus, given existence of regular
and normal approximations, where we assume algebraic dependence of
the parameter $T$, in the tangent space, then all normal
approximations, algebraically dependent on the parameter in the
tangent space, can be assumed regular (at least after adding a
regular approximation).

\newtheorem{Yodovich}[mainresult]{Proposition }
\begin{Yodovich} \label{Yodovich}
Assume $F_{T}$ with $L^{1}-$ dependence in the parameter and
existence of a normal and regular approximation algebraically
dependent of the parameter in the tangent space, then all normal
approximations, algebraically dependent of the parameter in the
tangent space, can be chosen as regular.
\end{Yodovich}

Note that when the parameter is with respect to the ramifier, we assume algebraic dependence
over transversals and tangents. There are numerous examples where $(dI)$ has a global (pseudo-)base,
but not $(I)$.
Finally, note that of $\Omega_{1}(dI)=\{ T  \quad r_{T}'F_{x}=F_{x}
\quad F_{x} \in (dI) \}$ and $\Omega_{2}(dI)=\{ T  \quad
r_{T}'F_{x}^{2}=F_{x}^{2} \quad F_{x} \in (dI) \}$, where $F \neq 0$
we have that $T \in \Omega_{2} \Rightarrow T \in \Omega_{1}$ iff
$r_{T}'$ is algebraic in the sense that it is geometrically
equivalent with a polynomial. Assume all approximations of a
parabolic singular point are on the form
$\eta_{T}=\alpha_{T}e^{\varphi_{T}}$, since we know that all normal
approximations are regular, we can assume the singularities for
$\alpha_{T}$ simple. Assume $\eta_{T}(x^{j}) \sim_{m}
\eta_{T}^{j}(x)$, then it is sufficient to consider the case where
$\eta_{j,T}=\frac{d^{j}}{d x^{j}} \eta_{T}$ has isolated
singularities. Since $\mid e^{-\varphi^{*}} \eta_{j,T}^{*} \mid < M$
as $\mid x^{*} \mid \rightarrow \infty$ implies $\mid e^{-
\varphi^{*}} \eta_{T}^{*} \mid < M$, as $\mid x^{*} \mid \rightarrow
\infty$. We will see that monotropy is a micro local invariant, this means that
it is sufficient to consider parabolic approximations for
$\eta^{*}$.

\vsp

Note that presence of lineality for the symbol, may result
in $\mbox{ Im }F$ in the space of hyperfunctions. We now note that
if $F$ is symmetric, entire and of finite type in $Exp$, then the
condition that $f$ represents a hypoelliptic operator, means that
for some $\lambda$, $(\mbox{ Im})^{\lambda} =\sum A_{j} F_{j}$ on a
domain of holomorphy, for constant coefficients and a global
pseudo-base representing the ideal of hypoelliptic operators. Thus,
symbols to hypoelliptic operators do not have imaginary part outside
the space of distributions and if hyperfunction representation is
necessary, we must have contribution of lineality in the infinity.

\subsection{ A generalized Cousin integral }

We denote with $\tilde{M}=-Ydx+Xdy$ and correspondingly for $\tilde{W}$
Assume $\tilde{M}$ exact and $\tilde{W}$ closed, then the form
corresponding to $\widehat{M}$ is exact after analytic continuation
and in the same manner for $\widehat{W}$. Note however that the forms
corresponding to $\widehat{M}$ and $\widehat{W}$ are not locally
holomorphic, that is we do not have locally isolated singularities
and the center case could appear.

\vsp

 Assume $\mu$ is a positively definite measure and consider
$$ \Phi_{\mu}(d \gamma)=\int_{P_{\mu}=0}d \mu(\gamma)$$
where $P_{\mu}$ is a polynomial and gives a local definition of
$\Delta$. Approximating a singular point through $d \gamma
\rightarrow 0$, then either $\Phi_{\mu}(d \gamma) \rightarrow 0$ or
we have existence of a point support measure $\mu'$ such that $\big[
\Phi_{\mu}+\Phi_{\mu'} \big](d \gamma) \rightarrow 0$. Thus, for the
measure corresponding to a hypoelliptic operator, we can choose
$\mu$ with point support. Assume
$\Phi_{\mu}(d(\gamma_{T}-\gamma_{0}))=\int_{\gamma_{T}-\gamma_{0}}d\mu$.
If $d \mu$ is a reduced measure, we must have
$\gamma_{T}=\gamma_{0}$. We know that if $d \mu$ is holomorphic
(that is holomorphic coefficients), then $d \mu$ will be reduced,
for $T$ close to $0$. Assume $d \mu$ continuous and locally bounded,
for all $T$ and that $\widetilde{d \mu}=d \mu+d \mu_{0}$, where $d
\mu_{0}$ is assumed with point support and $\widetilde{d \mu}$ is
holomorphic. Then $\int_{\gamma_{T}-\gamma_{0}} \widetilde{d \mu}=0$
implies $\gamma_{T}=\gamma_{0}$. Assume $\gamma_{T}$ a closed
contour and $\gamma_{0}$ a point, then for $T$ not close to $0$, we
have $\int_{\gamma_{T}-\gamma_{0}} d \mu=0$, implies $\gamma_{T}
\neq \gamma_{0}$. This case includes the case with a center (cf.
\cite{Bendixsson01}, Theorem 4).

\section{Stratification}

\subsection{Introduction}
If we consider a hypoelliptic analytic symbol $f$ as locally reduced, it is naturally necessary to use a stratification to define a globally hypoelliptic symbol. The model is centered around the set of lineality and we are always assuming the lineality locally is a subset of a domain of holomorphy,
which means that its local complement set is analytic, We consider it to be necessary for the concept of hypoellipticity to have an approximation property for $\log f$. We will discuss an interpolation property. Further, it is necessary to have a concept of orthogonality between the real and imaginary parts of the symbol.  

\subsection{ The arithmetic mean }

For the arithmetic mean, we have that $$\lim_{\epsilon \rightarrow
0} \int_{C_{\epsilon}} MV dz(T)=MV(0)$$ given that $MV$ is
holomorphic, regularly that is without a porteur (cf. \cite{Martineau}). If for all closed
contours $\int_{C_{\epsilon}} MV dz(T)=MV(0)$ implies
$C_{\epsilon}=\{ 0 \}$, then $MV$ is reduced for analytic
continuation. If $\int_{C_{\epsilon}-0} MV
dz(T)=0$ for all closed contours in a leaf $\mathcal{L}$, then the
form $MVdz(T)$ is closed in $\mathcal{L}$ and we have a mean value
property above for the arithmetic mean in $\mathcal{L}$. Further, the
closed contour $C_{\epsilon} \sim 0$ on $\mathcal{L}$.

\subsection{ The concept of stratification  }
Assume $X \subset Y$ are separable topologiclal vectorspaces. 
We say that $Y$ is a stratifiable space if it has the property that
to any open set $U$ we associate a sequence
$\{  U_{j} \}^{\infty}_{j=1}$ of open sets in $X$, such that
\begin{itemize}
\item[i)] $\overline{U_{n}} \subset U$ for all $n$
\item[ii)] $U=\cup^{\infty}_{j=1} U_{j}$
\item[iii)] $U \subset V$ implies $U_{n} \subset V_{n}$ for all $n$
\end{itemize}
Further, (cf. \cite{Cauty73}) given a topological vector space $X$
and with $Y$ as above, we can associate a topological vector space
$Z(X)$, such that $X$ is closed in $Z(X)$. We say that $X$ is
locally RA (retractible), if $X$ has a local extension property with respect to
the stratification. Particularly, if $\Gamma$ is closed in $Y$ and
$f$ is a continuous mapping $\Gamma \rightarrow X$, we have
existence of $\tilde{f}$ that maps $Y \rightarrow Z(X)$.

\subsection{ A stratification using averages }
A topological vector space $X$ is stratifiable, if for any open set
$U$, there is a continuous mapping $f_{U}$ $X \rightarrow nbhd 0$,
such that $f_{U}^{-1}(0)=X-U$ and if $U,V$ are open sets with $U
\subset V$, we have $f_{U} \leq f_{V}$. We will for this reason
study the averages $\M_{1} \geq \M_{2} \geq \ldots \geq \varphi$,
where the boundary $\M_{j}=\varphi$ is common for all the averages
and where $\M_{j} \rightarrow 0$, as $\varphi \rightarrow 0$. Let
$F_{1}=\{ \M_{1} \geq \varphi \}$ and let $f_{1}$ be a continuous
function such that $\mbox{ ker }f_{1}=\mbox{ bd }F_{1}$ and
$f_{1}=\M(\varphi)-\varphi$. If $\M$ is holomorphic and
$\M(x^{0})=\lim_{\rho \rightarrow 0} \int_{C_{\rho}} \M(x)dx$ and if
$C \subset \mbox{ bd }F_{1}$ and $\varphi(x^{0})=\lim_{\rho \rightarrow
0}\int_{C_{\rho}}\varphi(x)dx$. Further,
$f_{2}(\varphi)=\M_{2}(\varphi)-\varphi$ with $\mbox{ ker }f_{2}=\mbox{ bd
}F_{2}$, where $F_{2}=\{ f_{2} \geq 0 \}$ why $F_{2} \subset F_{1}$
and $f_{2} \leq f_{1}$, and so on.

\vsp

A stratification of $\dot{\mathcal{B}}$ can be mapped into a
stratification of $\Bm$, through $i_{a} : \dot{B} \rightarrow \Bm$
and $i_{a}(\varphi)=\varphi+a$, for a constant $a$. This is a compact mapping
with $i_{a}(\varphi_{j}-a)=i_{a}(\varphi_{j})-a$. $(\Bm)$ stratified in this
manner with topology induced of Schwartz type is $\mathcal{FS}$, why
the dual space $(\Bm)'$ is $\mathcal{DS}$ (cf. \cite{Martineau})

\subsection{ The arithmetic mean and duality in $L^{1}$ }
Assume $F_{1}=\{ \M(\phi) \geq \phi \}$ and $f_{1}=\M(\phi)-\phi$ and
$\Gamma=\{ f_{1}=0 \}$. If we assume $\M(\phi)$ holomorphic, we must
have that $Flux(\M(\phi))=0$. Note that
$\M(\phi^{\diamondsuit})=\M(\phi)^{\diamondsuit}$ in $L^{1}$ and
$\M(d \phi)=d \M(\phi)$. Thus, given a dynamical system with right hand sides
harmonic conjugates, satisfying the regularity conditions, we see
that the arithmetic mean satisfies a condition on vanishing flux. If $\Gamma$ is always reduced to a
dynamical system considered in $L^{1}$, the boundary problem is
solvable in $L^{1}$.

\vsp

Assume now that the boundary value problem is solvable for
$MV(\phi)$, that is assume $\Delta MV(\phi)=0$ on an open set
$\Omega$. Using duality with respect to the scalar product in
$L^{1}$, we consider
$$ 0 \rightarrow \phi \rightarrow MV(\phi) \rightarrow \Delta MV(\phi)
\rightarrow 0$$
$$ 0 \leftarrow MV_{-1}(\Delta \sigma) \leftarrow \Delta \sigma
\leftarrow \sigma \leftarrow 0$$ We are thus assuming $\Delta \sigma
\in L^{1}$ with $\sigma \in L^{1}$. Let $E=\Delta L^{1}$ and $X=\{
\phi \quad {}^{t}MV(\phi) \in E \}$, that is for a $f \in
L^{1}$, we have ${}^{t}MV(\phi)=\Delta f$ in $L^{1}$. More
precisely, we can describe $\Delta L^{1}=E$ through the closure of
$(M,W)$ with respect to $L^{1}$. Assume $\Phi \bot
MV(W)^{\diamondsuit}$ and $\Phi \bot MV(W)$ in $L^{1}$. Assume
$\Phi$ with support in a bounded neighborhood of the boundary
(restriction to strata). The relations will then also hold in
$L^{2}$ and we can apply Weyl's lemma to conclude $\Phi \in C^{1}$
locally. Assume in a neighborhood of the boundary that $0 = \mid
(\varphi,MV(W)) \mid \geq \mid (\varphi, W) \mid$. Thus, if the
problem is solvable for $MV$, it is solvable for $(M,W)$, given the
inequality above. We now
have $\varphi \in C^{1}$. The parametrix to the problem
then has a trivial kernel and the problem is solvable.

\newtheorem{strat}{ Proposition }[section]
\begin{strat}

The arithmetic means applied to $f$ (and $\log f$) form a stratification over $(\Bm)'$
associated to $f$ in a finitely generated symmetric ideal of analytic functions over a pseudoconvex domain with transversals given by a locally algebraic ramifier.
We have assumed parabolic singularities and no essential singularity in the infinity.
\end{strat}

 \subsection{Reduction to tangent space}
 Assume $F \sim V_{1}+i V_{2}$, and consider the condition

 \begin{equation} \label{nec_he}
 \frac{d}{d x} \log V_{1} \text{  reduced and  }
 \frac{d}{d y} \log V_{1} \text{ reduced  }
 \end{equation}
Given that $V_{1}+iV_{2}$ is hypoelliptic with $\vartheta=\log
V_{1}$, we have that if the property (\ref{nec_he}) holds for
$\M_{N}(\vartheta)$, then it also holds for $\vartheta$. Note also
that if $M \bot W$ with $\Top W=0$, then we can not conclude that
$\Top V_{2}$ has vanishing flux. However, if the condition
(\ref{nec_he}) is satisfied for $V_{1}$ and $M \bot W$, we can
conclude that $V_{1} \bot V_{2}$. Let $\tilde{M}=X d x+ Y d y$ and
$\tilde{W}=d F$. Then we can consider $F^{\diamondsuit}$ defined as
$\frac{d F}{d x}=\frac{d F^{\diamondsuit}}{d y}$ and $\frac{d F}{d
y}=-\frac{d F^{\diamondsuit}}{d y}$, so that $\tilde{M}=d
F^{\diamondsuit}$. If the involution is taken over $F,V,G$, where
$V$ is the Hamilton function, $F$ is the lifting function and $G$ is
a regular approximation, then we can relate the involutive set to a
condition $\int_{C_{0}} d F=0$. Assume $F_{N}$ corresponds to
$\M_{N}(F)$ and $C_{N}$ is the corresponding contour, such that
$\int_{C_{0}} d F \sim \int_{C_{N}} d F_{N}$ and $C_{N} \subset
C_{0} \subset C_{-N}$. Then, the conclusion is that the
stratification of negative order is a covering of the involutive
set.

\subsection{Example}
Assume for instance that $V=V_{1}+iV_{2}+\Delta$, such that
$\frac{\delta}{\delta x} \Delta=\frac{\delta}{\delta y} \Delta=0$
and where $\Delta$ is defined through involution and through the
conditions $\frac{\delta}{\delta x}(V_{1} \bot
V_{2})=\frac{\delta}{\delta y}(V_{1} \bot V_{2})$. Hypoellipticity
means that $\mbox{ supp } \Delta=\{ 0 \}$.

\subsection{The lineality as closed contours}

The lineality has a pre-image in the contour $C_{T}$ in the
following manner. Let $f_{T}=e^{\vartheta_{T}}$ and assume
$\vartheta_{T}-\vartheta \equiv 0$ and $\zeta_{T} \in \Delta$ locally (lineality).
Assume $C_{T}=F^{-1}\{ \vartheta_{T}-\vartheta \}$ describes a simple contour with an analytic parametrization,
then on $C_{T}$, $\M_{-N}(\vartheta_{T}-\vartheta) \leq \vartheta_{T}-\vartheta$.
Assume $\Delta=\{ \vartheta_{T}-\vartheta \equiv 0 \}$ locally analytic,
then we have locally $I(\Delta)=\{ \vartheta_{T} \quad
\vartheta_{T}-\vartheta \equiv 0 \quad \zeta_{T} \in \Delta \}$ and
$NI(\Delta)=\Delta$. This means that $C_{T} \mbox{:} 0\equiv
\M_{-N}(\vartheta_{T}-\vartheta)$ has a point in common with
$I(\Delta)$, that is $\{ \zeta_{T} \quad F^{-1}(\vartheta_{T}-\vartheta) \in C_{T} \}$ and $\Delta$ have points
in common. The contour $C_{T}$ gives a micro-local contribution, if
$\M_{-N}(\vartheta_{T}-\vartheta) \equiv 0$, that is for a given $F$, $\M_{-N}$ maps
locally the geometric ideal $I(\Delta)$ on to the closed contour $C_{T}$. 

\vsp

Consider again the problem if the zero set has points in common with $C_{T}$. If $\M_{-N} \in
\mathcal{D}_{L^{1}}^{(k')}$, we can assume that the restriction of a
complex operator $P(\delta_{T})$ to the real space, is such that
$P(\delta_{T})\sigma \mid_{\mathbb{R}} \sim \M_{-N}$ in
$\mathcal{D}_{L^{1}}^{(k')}$, where $\sigma \in L^{1}(\mathbb{R})$.
Extend the definition of $\sigma$ (standard complexify) to
$L^{1}(\mathbb{R}^{2})$. We can then, in a neighborhood of the
boundary corresponding to the symbol, assume that the parametrix to
$P(\delta_{T})$ is injective,
$E(P(\delta_{T})\sigma(\phi))=\sigma(\phi)-r_{T}$, where $r_{T}$ is
regularizing and $r_{T} \rightarrow 0$ as $T \rightarrow 0$. We must
assume that $\sigma$ is not identically $0$, but that $\sigma_{T}
\equiv 0$ on $C_{T}$. The proposition is now that $\phi_{T}$ has a
zero on $C_{T}$.
Assume for this reason that $\{ \sigma_{s} \}$ is a family of
measures, depending on a parameter as above, such that $\sigma_{s}
\rightarrow \delta_{\gamma_{0}}$ as $s \rightarrow 0$ continuously. This
is motivated by the condition that $C_{T}$ has an analytic
parametrization as a closed and simple contour. If $\sigma_{s}(\phi)
\equiv 0$ on $C_{T}$ as $s \rightarrow 0$, then we have existence of
$\gamma_{0}$ on $C_{T}$ such that $\phi(\gamma_{0})=0$.

\subsection{Further remarks on the stratification}

Assume a global pseudo base in the tangent space and that $F(d z)=f(z) d z$, where $f$ is
given by a locally reduced function. We are assuming $F$ has no lineality in the tangent space and that $\Delta$ can be given as a semi-algebraic set. If $F \in L^{1}$ in the parameter,
then $\frac{d F}{d z}=f(z)$ a.e. A sufficient (and necessary) condition for equality, is that $F$ is absolute continuous. For example, if $d f= f_{0} d g$, where $f_{0}$ corresponds to continuation. If $f$ is reduced with respect to analytic continuation (over strata) then $f_{0}$ is locally reduced. If $df,dg$ are of type $0$, then the same must hold for $f_{0}$.
If $g=hf$, then over $d h=0$, for $f_{0}$ to be reduced, we must have that $h$ is minimally defined. The first relation particularly means that $F$ preserves order of zero's if $f$
is regular, particularly $F$ maps exponentials onto exponentials. If $f$ is absolute continuous, then zerosets are mapped onto zerosets. When $F(e^{x})=e^{\tilde{F}(x)}$, if we assume $\tilde{F}(\overline{x})=\overline{\tilde{F(x)}}$. If we only have $F(e^{x})=\beta e^{\phi}$, where $\phi(\overline{x})=\overline{\phi(x)}$, then $e^{-\phi(x)} f(\overline{x})=\big[ \beta(\overline{x})-\overline{\beta(x)} \big] e^{- 2 i \mbox{ Im } \phi}$. Note that reducedness for $\beta$ is not necessarily symmetric.

\subsection{ Condition $(M_{1})$ relative the stratification }
We are assuming $\Top
 \mbox{:} x \frac{d \eta}{d x} \rightarrow x^{*} \frac{d \eta^{\vartriangle}}{d
 x^{*}}$ and $x \frac{d \eta}{d x} \neq 0$
 and that the systems $(M,W)$... are regular. In particular, we assume $\ram
 \frac{d F}{dx}=\frac{d F}{d x} \ram$ and in the same manner for
 $y$. Further, we are assuming that $F(X,Y) \sim R(M,W)$
 and $F \mbox{:} 0 \rightarrow \Delta$ and $F^{-1} \mbox{:} \Delta
 \rightarrow 0$ locally. We are thus assuming that the lineality is defined as
 "independent" of the system. The stratification is formed over $(\Bm)'$ and is relative
 $\Top W=e^{\phi^{\vartriangle}}$, where $\pm \phi^{\vartriangle} > 0$ on each strata. We use the
 following concept of condition $(M_{1})$. Let $L_{N}(\omega_{T})=\frac{1}{\mid \rho
 \mid} \int_{C_{N}} \M_{N}(\omega_{T})dz(T)$, where $C_{N}$ is a closed
 contour, parameterized through $T$, such that $T \rightarrow C_{N}
 \downarrow \{ 0 \}$ as $N \uparrow \infty$ and $\frac{1}{T}
 \rightarrow C_{-N} \uparrow \mbox{ bd } \mathbf{C}$, as $-N
 \uparrow \infty$ and where $\rho$ is the radius for $C_{N}$, where
 $\rho=\rho(N,T)$. The condition $(M_{1})$ is that $\lim_{N \rightarrow
 \infty} L_{N}(\omega_{T})$ is regular, that is that the functional
 corresponding to $\M_{N}(\omega_{T})$ is of real type. Note that if
 $\sigma_{N}^{\vartriangle} \in L^{1}$, with $\parallel \sigma_{N}^{\vartriangle}
 \parallel=1$, we have the same argument for $\M_{-N}$ and $L_{-N}$.
 We will now argue that if $\omega_{T}^{\vartriangle} \equiv 0$, for $T$
 small, then $\M_{-N}(\omega_{T}^{\vartriangle}) \equiv 0$ on the contour
 $C_{-N}$. Let $X_{N}=\{ \omega_{T} \leq \M_{N}(\omega_{T}) \}$, for $N
 \geq N_{0}$, where $\M_{N}(\omega_{T}) \in L^{1}$ and algebraic in
 $T$, for $T$ small. If we extend the definition of
 $\sigma_{N}^{\vartriangle}$, such that $\sigma_{N}^{\vartriangle}$ is the evaluation
 functional on the boundary $\{ \omega_{T}=\M_{N}(\omega_{T}) \}$, with
 $\omega_{T} \in L^{1}$, then if
 $<\M_{-N}(\omega_{T}^{\vartriangle}),\M_{N}(\omega_{T})>_{1} \equiv 0$, we can chose
 $\M_{-N}(\omega_{T}^{\vartriangle})=\sigma_{N}^{\vartriangle}$ on the inner of $X_{N}$, why
 $\omega_{T}^{\vartriangle} \bot \M_{2N}(\omega_{T})$, for $T$ small. If $\M_{2N}-I$
 is locally algebraic, then $\omega_{T} \in \Gamma = \Gamma^{\bot}$.
 We are assuming the Lagrange condition
 $\Gamma=I(\Delta)=I(\Delta)^{\bot}=\Gamma^{\bot}$

\vsp

 Assume $\omega_{T}^{\vartriangle} \equiv \omega_{\infty}^{\vartriangle}$ ($\sim \omega_{1/T}$),
 for $T$ small, on a set with complex dimension, then we must have
 existence of $\M_{-N}$ as described above, such that for $N$ large,  $<
 \M_{-N}(\omega_{T}^{\vartriangle}),\M_{N}(\omega_{T})>_{1} \equiv 0$ or equivalently
 $< \omega_{T}^{\vartriangle},\M_{2N}(\omega_{T})>_{1} \equiv 0$, according to the
 conditions, we have that $\M_{-N}(\omega_{T}^{\vartriangle}) \equiv 0$ on
 $C_{-N}$, for $N$ large. We are assuming that $C_{-N}$ includes the
 real infinity, as $N \rightarrow \infty$. Conversely, if
 $L_{-N}(\omega_{T}^{\vartriangle}) \neq 0$, as $N > N_{0}$ implies $\omega_{T}$ is
 not $\equiv \omega_{\infty}$, in the real infinity. The conclusion is
 that if the stratification has condition $(M_{1})$ in the infinity, it is not possible to
 have lineality.

 \vsp

 Consider the limit $L_{N}(\omega_{T})=\int_{C_{\delta}}
 \M_{N}(\omega_{T}-\omega)d z$, where $z(T) \in C_{\delta}$ a closed
 contour of radius $\delta$ and let $A_{N}$ be the porteur set to
 this limit considered in $H'$. Obviously, we have $A_{N} \subset
 \Delta$, for $N \geq 0$. Consider the stratification of $(\Bm)$
 with $\{ X_{N}^{\vartriangle} \}$, that is a stratification using the means
 $\M_{N}$. If $L_{N}$ are not regular, that is $A_{N} \neq \{ 0 \}$,
 then we have on a connected set that $\omega_{T}-\omega \equiv 0$. (We are assuming
 Schwartz type topology for the symbol space). Conversely, consider
 the stratification of $(\Bm)'$ and the contours $\{ C_{T} \}$ that
 contribute to $\Delta$ through common points. In this case, if
 $\M_{-N}$ are of real type, there is no possibility of lineality.
 Thus, given an operator with lineality, we do not have condition $(M_{1})$
 for $(\Bm)'$ in the stratification using $\M_{-N}$.

 \newtheorem{condition (M)}{ Proposition }[section]
 \begin{condition (M)}
 If the stratification that we are considering has condition $(M_{1})$, that is if all the means are
 of real type, then the symbol ideal is locally reduced and conversely.
 \end{condition (M)}

We will discuss two other similar topological conditions in a later section. Since it is topological, we prefer the set of lineality to characterize hypoellipticity. The condition $(M_{1})$ at the boundary, means that the boundary behavior
does not influence the microlocal behaviour in the infinity. A globally hypoelliptic operator
is in this context a globally defined operator that is hypoelliptic in the infinity and for
which the topology for the symbol space has condition $(M_{1})$ (or a similar topological condition) at the boundary.

\subsection{Reduction to real type}
Assume $F$ holomorphic and of finite exponential type. Further that $F$ has finitely many zero's on
$X \backslash U_{0}$, where $X$ is assumed a bounded domain and $U_{0}$ is a neighborhood of the infinity. Further, we assume that the zero's
$P_{1}, \ldots ,P_{\nu}$ are isolated and of finite order. Assume $U_{1}$ is a neighborhood of
$P_{1}$ that does not contain any other zero's. Then we have on $X$ a holomorphic function $F_{1}$,
such that $F-F_{1}$ is of type $0$ on $X$ and $F_{1}$ is of type $- \infty$ on $U_{1}$. The
remaining $P_{j}'s$ are treated in the same way. Thus, $F-\sum_{j} F_{j}$ is of type $0$ on $X$ and
each $F_{j}$ is of type $- \infty$ on the corresponding $U_{j}$.

\subsection{ Remarks on a spectral mapping problem}

The definition of the mapping $\Top$ starts with $-Y dx+X dy
\rightarrow -\widehat{Y} dx^{*}+ \widehat{X} dy^{*}$ and we are
requiring $\{ W=0 \} \rightarrow \{ \widehat{W}=0 \} \rightarrow \{
\widehat{\widehat{W}}=0 \}$. We consider the multipliers $\chi X=Y$,
$\chi^{\vartriangle} \widehat{X}=\widehat{Y}$, $\lambda H=G$ and $\lambda^{\vartriangle}
\widehat{H}=\widehat{G}$. We assume $\Top \mbox{:} \chi \rightarrow
\chi^{\vartriangle}$ and $\{ \eta=\chi \} \rightarrow \{ \eta^{\vartriangle}=\chi^{\vartriangle} \}$.
We have that $\Top$ preserves parabolic points, but is usually
not a contact transform. If $\Top$ has the property that it maps
constants on constants and exponentials on exponentials, we know
that $\Top$ preserves parabolic approximations.
Through the condition on vanishing flux, we can assume $(w,\Top w)$
pure and that $\Top$ preserves analyticity.

\vsp

Consider $(J)=$ $\{ f$ $\int_{(I)} f d \sigma(t)$ $=0$ $\tilde{V}
\}$, where $\tilde{V}$ is a geometric set.
One of the more difficult problems in our approach is to see that
the spectral mapping result we use respects the stratification, that is if
starting with a stratification of $(\Bm)'$ and $\widehat{W}$, $\{
X_{j}^{*} \}$, we have that the sets $\{ X_{j}^{\vartriangle} \}$ where $\Top X_{j}=X_{j}^{\vartriangle}$
constitute a stratification. Consider $\Phi^{\vartriangle}=\{\vartheta^{\vartriangle}
\quad  e^{-\vartheta^{\vartriangle}}\chi^{\vartriangle}=const \quad \exists \vartheta^{\vartriangle} \}$
and $\Phi=\{ \vartheta \quad e^{-\vartheta}\chi=const \quad \exists
\vartheta \}$. Consider the Legendre transform $R$, according to
$\mbox{RE} < R(\chi),\chi >-1$. Let
$R(e^{\vartheta})=\widehat{R}(\vartheta)=\widehat{I}R(\vartheta)=e^{\vartheta^{*}}$
and we note that $\big[ \widehat{R},I \big]=\big[ \widehat{I},R
\big]$ implies that $R$ is algebraic in $H'$ over $\Phi^{\vartriangle}
\rightarrow \Phi^{*} \rightarrow \Phi$ and over a regular parabolic approximation, we can
argue as in the spectral mapping theorem. For a hypoelliptic system, the continuation to $\Top$ is algebraic
and the stratification of $X^{*}$ gives a stratification of $X^{\vartriangle}$. We can conversely argue that if these stratifications are equivalent, the system has no lineality.

\section{Topology}
\subsection{Introduction}
The concept of hypoellipticity is dependent on topology and we will use the monotropic functionals both for limits in the symbol space and for the equations in the operator space. 
The topological arguments are comparative and we compare with the more familiar hyperfunctions. However there are geometric sets that are characteristic for hypoellipticity, such as lineality and the set of orthogonality, for all topologies that we consider. Several parameters are necessary to define the class of hypoelliptic symbols.
We give the approximation of the operator using operators dependent on a parameter and a second parameter is used to trace the transversal in determining microlocal contribution. Since this is an analytical study and not a geometrical, we do not attempt to minimize the number of parameters.

\subsection{Topological fundamentals}
 The space $H(V)$,
where $V$ is a complex analytic variety, countable in the infinity, is the space of holomorphic
functions with topology of uniform convergence on compact sets. This is a Frechet type of space (FS)
and the dual space is denoted (DS). Given F-spaces $\{ E_{i} \}$, if $i \mbox{:} E_{i+1} \rightarrow E_{i}$
the projective space is (FS). If ${}^{t} i \mbox{:} E_{i} \rightarrow E_{i+1}$ compact, the
inductive limit is compact. We start with a topology of Schwartz type, that is given a separated space $E$, if $V$ is a convex disc neighborhood of
the origin in $E$, then we have existence of a convex disc in $E$ that is a neighborhood of the origin
such that $U \subset V$ and such that $E_{\tilde{U}} \rightarrow E_{\tilde{V}}$ is compact,
where $E_{\tilde{U}}$ is the completion of the normalized set $E_{U}$. The topological arguments in this study are comparative. The symbols modulo regularizing action are considered in i neighborhood of the real space where we compare with monotropic functionals and the $\mathcal{D}_{L{p}}'$ spaces ($p=1,2)$. We also give a brief comparison with the hyperfunctions.

\newtheorem{singularities}{ Proposition }[section]
\begin{singularities}
If $(I)$ is an ideal of holomorphy with topology of Schwartz type
and a compact translation, consider $rad(I)$ with Schwartz type
topology and a weakly compact translation, if $\psi \sim_{m} 0$ in
the $\mid \zeta \mid-$ infinity, is in $rad(I)$, then $\{ d_{\zeta}
\psi=\psi=0 \}$ is nowhere dense in $N(I)$.
\end{singularities}

\newtheorem{Bounded}[singularities]{ Proposition }
\begin{Bounded}
Assume $(J)=\mbox{ ker }h$ a finitely generated ideal with topology
of Schwartz type and $r_{T}'$ weakly compact. Assume for all $\psi
\in (J)$, we have $h(r_{T}' \psi)/\psi \sim_{m} 0$ in the $\zeta-$
infinity. If $\eta \neq const.$, we have that $\psi$ is in a bounded
set with respect to the origin.
\end{Bounded}

Proof:\\
We can prove an estimate
$$ \mid \eta -c_{R} \mid < 1/\mid \zeta \mid \qquad \mid \zeta \mid
> R $$
for a constant $c_{R}$ and $R$ sufficiently large. Thus, for $\phi
\in (J)$
$$ \mid \phi \mid < c_{1}/\mid \zeta \mid +c_{2} \mid \eta(\phi)
\mid \qquad \mid \zeta \mid > R$$ for constants $c_{1},c_{2}$.
Symmetry follows from the conditions on $r_{T}'$. $\Box$

\subsection{ Monotropic functionals }

Assume $\Bm$ test functions, that is $C^{\infty}-$ functions,
bounded by a small constant in the infinity, such that $\dot{\B}
\subset \Bm \subset \mathcal{E}$ and $(\Bm)' \subset \DL$. The
Fourier transform over the real space is $Pf_{0}$, where $P$ is a
polynomial and $f_{0}$ is a continuous function. We will modify
$f_{0}$ to an $\epsilon-$ neighborhood of the real space as follows

\begin{itemize}
\item[i)] $F_{0}$ is continuous on the real space and locally
bounded on an $\epsilon-$ neighborhood of the real space.
\item[ii)] we have existence of $\lim_{\Gamma \rightarrow 0}F_{0}(\xi +
\Gamma)$, for any line $\Gamma$
\item[iii)]any line $\Gamma \subset \Delta(F_{0})$ is such that
$\Gamma \subset \Omega$, where $\Omega$ is a domain of holomorphy.

\end{itemize}
Note that the difference $\tau_{\Gamma}F_{0}-F_{0}$, even when it is
not holomorphic, will preserve constant value over the lineality
corresponding to $F_{0}$. Finally, assume

\begin{itemize}
\item[iv)]$F_{0} \sim_{m} W_{0}$, where $W_{0}$ is holomorphic and
in $Exp$ of finite type.
\end{itemize}
We then have existence of $B_{\Gamma}$ (modulo monotropy). Assume
further that the translation is algebraic over $PF_{0}$ and for
$W_{0}$, that the lineality is quasi-porteur (cf. \cite{Martineau}) for $B_{\Gamma}$.

\vsp

The first observation is that if $f_{T} \in L^{1}$ then $\M_{N}(\frac{d^{N}}{ T^{N}} f_{T})=$ $\sigma_{T} \in L^{1}$.
Particularly, we have $f \in \DL^{m}$ implies $\M_{m}(f) \in L^{1}$. As $d f=\alpha d x$, we have
if $d \M(f)=\beta d x $, then we must have $\beta=\M(\alpha)$, thus in $L^{1}$, $\frac{d}{d T} \M(\alpha)=\M(\frac{d}{d T} \alpha)$
We now argue that $\M_{-N}$ is surjective in $\Bm'$. Consider for this reason $\M_{N}$ in $\Bm$ and assume that
$(I)$ is defined by $ f \in (I)$ $\Leftrightarrow f \in \Bm$ and $f \leq \M(f)$.
We then have $\M(f-f_{0})=0$ implies $f-f_{0}=0$. Thus, we have that $\M_{-N}$ is surjective over
$(I)'$. As $(I) \subset \Bm$, we must have $\Bm' \subset (I)'$, why the surjectivity follows for $\Bm'$.
Note that $f=e^{\varphi}$ with $\varphi$ subharmonic, if we let $\tilde{\M}(e^{\varphi})=e^{\M(\varphi)}$
and $\M(\varphi-\varphi_{0})=0$, then $\varphi-\varphi_{0}=0$ implies $f=f_{0}$.

\vsp

We have studied regular approximations according to
$F(\zeta+T_{j})=F(\zeta)+c_{j}$ as $0 \neq c_{j} \rightarrow 0$ as
$T_{j} \rightarrow T_{0}$. Note that more generally, for $dF(\zeta +
T)-dF(\zeta)=dL_{T}(\zeta)$ with for instance $dL_{T} \sim_{m}0$ as
$\mid \zeta \mid \rightarrow \infty$.
The projection method gives that $f \sim_{m}0$ as $\mid \zeta \mid
\rightarrow \infty$, means existence of $g$ holomorphic such that $g
\rightarrow 0$ as $\mid \zeta \mid \rightarrow \infty$. Thus we have
existence of a polynomial $P$ such that $\mid
g(\zeta)-P(\frac{1}{\zeta}) \mid < \epsilon$ as $\mid \zeta \mid
\rightarrow \infty$, where we have assumed $f=\tau_{\epsilon}g$.
Note that $1/\zeta=(1/\zeta_{1},\ldots,1/\zeta_{n})$ and we can assume the condition
in some variabels and assume the others fixed and in the finite plane.

\subsection{ Algebraicity for exponential representations }

Consider the following problem, when for a continuous homomorphism $L$ and $\mathcal{L}'=\{ \exists \mbox{!  } \eta \quad L(e^{\psi})=e^{< \eta,\psi>} \}$, do we have $L \in \mathcal{L}'$. Let $\mathcal{L}^{0}=\{ \exists \mbox{!  } \eta \quad L(\psi)=< \eta,\psi > \}$, where in this case $L$ is assumed continuous and linear. If $X=H(\Omega)$, for an open set $\Omega$, we assume $\widehat{X}=\{ e^{\psi} \quad \psi \in H(\Omega) \}$. If we have existence of $\eta_{x}$, for $x$ fix, such that $L \in \mathcal{L}^{0}(\widehat{X})$, we have $L \in (\widehat{X})'$, . Assume for $M \in (X)'$, that $M(\psi)=< \widehat{L},\psi>=<L,\widehat{\psi}>$, then $L=\mathcal{F}^{-1}M$. When $L \in \mathcal{L}'$ and if $L$ is algebraic in $e^{\psi}$, that is linear in $\psi$, we have $L(e^{\psi_{1}+\psi_{2}})=L(e^{\psi_{1}})L(e^{\psi_{2}})=e^{<\eta,\psi_{1}>}e^{<\eta,\psi_{2}>}$. Further if $L_{1},L_{2} \in \mathcal{L}'$, we have $L_{1}L_{2}(e^{\psi})=e^{<\eta_{1}+\eta_{2},\psi>}=e^{<\eta_{1},\psi>}e^{<\eta_{2},\psi>}$. If $\big[ \widehat{I},N \big](\psi)=e^{<\eta,\psi>}$ and $\big[ \widehat{N},I \big](\psi)=<\eta,\psi>$, we assume the commutator $C$ such that $C \big[ \widehat{I},N \big]=\big[ \widehat{N},I \big]$, then $\mathcal{F}^{-1} C \mathcal{F} \big[ I,N \big]=\big[ N,I \big]$.
If $\big[ \widehat{I},N \big]=\big[ \widehat{N},I \big]$, we say that $N$ is algebraic. Let $<N(\psi),\theta>=<\psi,{}^{t}N(\theta)>$. Then $N(\psi) \in (\widehat{X})'$ implies ${}^{t}N(e^{\theta}) \in X'$, for $x$ fix and ${}^{t}N \widehat{I}=\big[ \mathcal{F} N \big]^{t}$. We have $<\widehat{I}(\psi),\theta>=<\psi,I {}^{t} \mathcal{F} \theta >$ iff $< e^{\psi},\theta>=<\psi,{}^{t} \mathcal{F} \theta>$.

 \vsp

Assume $\Sigma \ni \gamma \rightarrow \mid \gamma \mid \in
\mathcal{R}$ is on the form $(e^{\varphi},h(e^{\varphi}))$. Assume
$L$ within a constant is algebraic over $\mathcal{R}$ (does not
imply algebraic over $\gamma$). Note
that if $\pi$ is the projection $\Sigma \rightarrow \mathcal{R}$ and
$\pi^{-1} \mathcal{R}=\tilde{\Sigma}$, then $\tilde{\Sigma}$ may
have points in the edge, even when $\Sigma$ does not. Thus. there
may be points in common for $L_{T} \in (I_{\tilde{\Sigma}})$, that
are not present for $L_{T} \in (I_{\Sigma})$.

\subsection{Some generalizations}

Assume $V_{1} \subset V \subset V_{2}$, where $V_{1},V_{2}$ are
semi-algebraic and $V$ analytic. Assume $V_{1} \subset
\Omega_{1}$ a domain of holomorphy, such that the limit $B_{\Gamma}$
is independent of starting-point, then $V_{1}$ is quasi-porteur to
$B_{\Gamma}$ and the same follows for $V$. Further, since $V$ is
analytic, $V$ is porteur to $B_{\Gamma}$. Assume $\tilde{V}$ the
extension to full lines. Assume $g_{1}$ algebraic, such that if
$V_{1}=\{ p_{1}(\gamma) \geq 0 \}$, $g_{1}p_{1}(\gamma)=\gamma$
locally. Then $g(0)=0$ and $V_{1} \rightarrow \gamma \rightarrow
\tilde{\gamma}$, where the last mapping is into the wave front set,
but regular approximations are assumed in $V_{1} \subset V_{2}$. If
$V_{1}$ is porteur to a functional $T \in H'$, we can chose $V_{1}$
as a cone which through the topology can be assumed compact. Note
that $V_{1}^{0}=\{ \gamma_{T}^{*} \quad \mbox{ Re }
<\gamma_{T},\gamma_{T}^{*}> \leq 1  \quad \gamma_{T} \in V_{1} \}$
and $bdV_{1}^{0}$ corresponds to the orthogonal complement to
$bdV_{1}$. If $B_{\Gamma}$ has indicator $h_{0}$, we have for
$\Gamma$ a compact, convex cone $\ni 0$ with inner points, that
$B_{\Gamma}$ is portable by $\Gamma$ $\Leftrightarrow h_{0} \leq 1$
on $bd \Gamma^{0}$.

\vsp

If $\gamma$ is defined by a homomorphism $h$ such that $h^{N}=1$, we
have that all regular approximations of a singular point $P$,
$\gamma \rightarrow P$ as $t \rightarrow \infty$, can be seen as on
one side of a hyperplane $\{ (x,h(x)) \quad dh(x) \geq \mu dx \}$,
for a constant $\mu$. More precisely, for a curve that reaches
$\Sigma$ as $t \rightarrow \infty$, if the part of the curve that is
situated outside $\Sigma$ is finitely generated, we claim that
$\gamma^{\bot}$ can be chosen locally on side of a hyperplane (cf.
section on paradoxal arguments). We are assuming in the following
that $x$ is reduced. We have $\{ p_{1}(x,y)=0 \} \sim \{
\tilde{p}(\eta)=0 \}$, for a polynomial $\tilde{p}$ and $\eta=y/x$.
Further, there are polynomials in $x,y$, $r_{1},s_{1}$ such that
$p_{1}(\gamma_{T})=\frac{d}{dT}r_{1}(\gamma_{T})=s_{T}(\frac{d}{dT}\gamma_{T})$.
If $p_{1}$ is reduced, there is a polynomial $q_{1}$ such that
$p_{1}=q_{1}^{2}$. Assume $F_{T} \rightarrow 0$ and
$\frac{d}{dT}F_{T}>0$ over $\{ p_{1}(\gamma_{T})>0 \}$. If we let
$\frac{d}{dT}F_{T} \sim p_{1}$, we are assuming $F_{T}$ conformal in
$T$ and algebraic in $x,y$. Assume $\tilde{V}_{1},\tilde{V}_{2}$ are
extensions to full lines with indicators $h_{1}$ and $h_{2}$. It is
then sufficient to prove that $h_{1} \leq c h_{2}$, for a constant
$c$, to have $\tilde{V}_{1}
 \subset \tilde{V}_{2}$. Further, $c_{1}h_{1} \leq h_{2} \leq
 c_{2}h_{1}$, for constants $c_{1},c_{2}$, gives $\tilde{V}_{1} \sim
 \tilde{V}_{2}$. If $V_{1}$ is a semi-algebraic quasi-porteur to a
 functional for instance defining the first surface to the symbol
 with $V_{1} \subset V$, where $V$ is analytic, then $V$ is porteur.
 Assume that $p_{1},p_{2}$ have the same micro-local properties, in
 the sense that their sets of constant sign coincides. Assume
 $V_{i}=\{ p_{i}(\gamma_{T}) \geq 0 \}$ and $V_{1} \subset V \subset
 V_{2}$ and $V=\{ g=0 \}$ with $g$ analytic, then $\{ p_{2}<0 \}
 \subset \text{ supp }g \subset \{p_{1}<0 \}$, Thus, if singular
 points are in $\{ g=0 \}$, then regular points will be in an
 "octant". We could say that the micro-local contribution from the
 symbol, is given by this "octant".

\subsection{A comparison of hyperfunctions and monotropic functionals}
If the symbol $F_{T}(\gamma)$ preserves
a constant value in the $\gamma$- infinity, then $F_{T} \in \Bm$,
that is it is $C^{\infty}$ and bounded by a small constant in the
infinity. In this case the Cauchy inequalities can be satisfied for
a monotropic function, that is there is a $\varphi_{T} \in \ra$
(real-analytic functions) such that $F_{T} \sim_{m} \varphi_{T}$. If
for $F_{T}(\gamma)=\sum_{\alpha}F_{\alpha,T}/\alpha !
\gamma^{\alpha}$ there is a number $\rho$ with $\rho < A$, for a
constant $A$, such that $\rho^{\alpha} \sup \mid F_{\alpha,T} \mid
\rightarrow 0$ as $\mid \alpha \mid \rightarrow \infty$, then we
have that $F_{T}$ is entire in $\gamma$ and of exponential type $A$
(\cite{Martineau}). Note that a sufficient condition for existence
of a global pseudo-base for the symbol ideal, is that it has an
induced topology with Oka's property. 

\vsp

If $\gamma_{T}$ is in $\ra(\Omega_{\zeta})$, then $F_{T} \in
\hyp(\Omega_{\zeta})$ that is hyperfunctions with compact support
(\cite{Komatsu71}). In the case where $F$ is real-analytic, then so
is $f$. If we assume instead that
$G_{T,k}=\frac{d^{k}}{dT^{k}}F_{T}$ has isolated singularities at
the boundary and preserves constant value in the infinity, then
intuitively we would have at worst algebraic singularities in the
infinity. Assume $f \in \Bm$ in $x,y$, then we have
$f(x,y)=g(x,y)+P(\frac{1}{x},\frac{1}{y})$, where $g$ is radial and
bounded by a small number in the infinity and $P$ is polynomial.
Further, if $f\in \Bm$ there is a $\varphi \in \ra$ such that $\mid
f - \varphi \mid < \epsilon$ at the boundary, for a small number
$\epsilon$.

\vsp

 An important difference between $\Bm$ and $\ra$ is the algebraic
properties. A function $f$ is in $\ra$ if both its real and
imaginary parts satisfy the Cauchy's inequalities. More precisely,
assume $L \in \B(\Omega)$ where $L^{2}$ is defined by composition,
such that $\phi \in \mathfrak{A}(\Omega)$ and $L^{2}(\phi)=L(J(D)
\varphi)$, where $J(D)$ is a local elliptic operator (cf.
\cite{Kaneko}) and where $L \in \mathcal{E}'(\Omega)$ such that
$J(D)\varphi \in \mathcal{E}$. Then $\varphi \in \mathfrak{A}$ and
also $J(D)\varphi \in \mathfrak{A}$, from the properties of $J(D)$.
Any element in $\mathfrak{A}$ has a representation through $J(D)
\varphi$ as above, why $L$ is defined on $\mathfrak{A}$ and $L \in
\B_{K}$.
 However we can
have $1/f^{N} \rightarrow 0$ in the infinity, for some iterate $N$
without having $1/f \in \Bm$, for instance if $f$ is the symbol to a
self-adjoint operator partially hypoelliptic in $\mathcal{D}$ with $\mbox{ Re } f
\prec \prec \mbox{ Im } f$. Thus we do not expect a radical behavior in the
case of monotropic functionals.

\vsp

 Another important difference is the global property of the
hyperfunctions (\cite{Komatsu71}), which is not present with the
monotropic functionals. However, we can give the following argument.
Let $F(T,\gamma)=F(\gamma_{T})$, where $\gamma=(x,y)$ and $T \in
V=\cup^{N}_{j=1}V_{j}$, the parameter space. Assume $F$ algebraic in
the parameter in the sense that
$$ F(T_{1}.T_{2},\gamma)=F(T_{1},\gamma).F(T_{2},\gamma)$$
(the dot signifies concatenation of curve segments). We are assuming
$\gamma_{T}$ real-analytic as $T \in V$, but we are not assuming $T
\rightarrow \gamma_{T}$ algebraic in $T$. With these conditions we
do not necessarily have that $F_{T}(\gamma)$ is algebraic in $x$ and
$y$. For the symbol now holds that if $f(r_{T} \zeta)=0$ for all
$T_{j} \in V_{j}$ and $\zeta$ fixed, then we have $f(r_{T} \zeta)=0$
for $T=T_{1} \ldots T_{N} \in V$.

\subsection{The operator space}

Given an open set $\Omega \subset \R$ and a differential operator $P$ with constant coefficients,
we have that $P(D) \mathcal{A}(\Omega)=\mathcal{A}(\Omega)$. If $P$ a differential operator
hypoelliptic (in the sense of $\mathcal{A}$), we have $P(D) \mathcal{B}(\Omega)=\mathcal{B}(\Omega)$

\vsp

The global property is particularly
interesting in connection with the solvability problem. It is known
for an elliptic d.o $P(D)$, that $P \ra(\Omega)=\ra(\Omega)$
(\cite{Kawai71}). For an elliptic d.o $P(D)$, (hypoelliptic in the
sense of $\ra$), we have that $P(D) \B(\Omega)=\B(\Omega)$ (Harvey
\cite{Komatsu71}), where $\Omega$ is an open set in $\R$. Note also
that through majorants an operator with coefficients $\leq \mid
c_{\alpha} \mid$ and such that $\limsup \mid c_{\alpha}
\mid^{\frac{1}{\mid \alpha \mid}} \rightarrow 0$, as $\mid \alpha
\mid \rightarrow \infty$ (compare exponentially finite type) maps
$\B(\Omega) \rightarrow \B(\Omega)$ (cf. \cite{Kaneko})

\vsp

 We have not discussed
hypoellipticity in the sense of monotropic functionals, but we can
relate to $\ra$-hypoellipticity using monotropy. We assume the
proposition (\ref{conjecture}). We can now use that monotropy is
micro-locally invariant in the
symbol-space. Assume for this reason $P$ a ps.d.o hypoelliptic in
$\ra-$ sense. Further that $Pu \in \Bm$ and $\phi \in \ra$ such that
$\phi= Pv$ in $\ra$ and $\mid Pu - \phi \mid < \epsilon$ at the
boundary $\Gamma$. The problem is now to prove existence of a symbol
$P_{1}$ such that $\Delta(P)=\Delta(P_{1})$ and $P_{1}u \in \ra$
with $\mid P_{1}u- \phi \mid < \epsilon$ at the boundary. Define
$P_{1}$ such that $\widehat{P_{1}} \sim_{m} \widehat{P}$, that is
$\tau_{\epsilon} \widehat{P_{1}} \sim_{0} \widehat{P}$. The
conclusion is that given a ps.d.o $P$, hypoelliptic in $\ra-$ sense,
there is a ps.d.o $P_{1}$ with the same lineality ($=\{ 0 \}$), such
that $WF_{a}(P_{1}u)=WF_{a}(u)$, why using the claim (\ref{conjecture}), we
have that $P$ is hypoelliptic in $\Bm$.

\vsp

For the discussion of the symbol space, we will use a topological
argument. Assume the topology of Schwartz type and with (weakly)
compact translation.
 Let $\Omega_{0}^{(j)}=\{ \Gamma \quad
F(\gamma_{j})(\zeta + \Gamma)=F(\gamma_{j})(\zeta) \quad \exists \gamma_{j} \}$, where $F_{\Gamma}$ are assumed to satisfy the regularity conditions for the dynamical system. 
 We can now prove for a sequence of $\gamma_{j}$, that approximate a singularity, $\Omega_{0}^{(j)}
\downarrow \{ 0 \}$ as $j \uparrow \infty$. Let $J$ be defined by
$N(J)=\Omega_{0}^{(j)}$, for some $j$ and let $\overline{(I_{0})}=\{
\gamma \quad F(\gamma)(\zeta + \Gamma)=F(\gamma)(\zeta)+C_{\Gamma}
\text{  for some constant }C_{\Gamma} \}$. Then we have existence of
$J$ as above with $rad(J) \sim \overline{I_{0}}$. Further, there is
a regular sequence of $\Gamma_{j}$, such that $F(\zeta +
\Gamma_{j})=F(\zeta)+C_{\Gamma_{j}}$ where $0 \neq C_{\Gamma_{j}}
\rightarrow 0$ as $\Gamma_{j} \rightarrow \Omega_{0}$. In the same
manner, if the dependence of $T$ is holomorphic for $\ram$, we have
existence of a regular sequence of $T$ outside $\Omega_{0}'=\{ T
\quad F(r_{T} \zeta)=F(\zeta)+C_{T} \quad C_{T}=0 \}$, such that $0
\neq C_{T} \rightarrow 0$ as $T \rightarrow \Omega_{0}'$. This
motivates why there is no loss of generality in assuming that
$r_{T}$ behaves locally as translation, in a regular approximation
of $\Omega_{0}'$.

\section{The mapping $\Top$}
\subsection{Introduction}
We have seen that certain trace sets (clustersets) are characteristic for hypoellipticity,
more precisely the absence of these sets is necessary. The mapping $\Top$ which is derived from dynamical systems theory (\cite{Bendixsson01}) will in this study be used to define and describe these sets. Characteristic for hypoellipticity, assuming the real and imaginary parts of the symbol are orthogonal, is that $\Top$, given as a continuation of the contact transform (Legendre), is (topologically) algebraic. 

\subsection{ Systems of multipliers } Consider the system with right
hand sides $(X,Y)$ and $\eta^{*} \widehat{X}=\widehat{Y}$. In the
same manner, to the system $(M,W)$, $\gamma H=G$. This can be seen
as a multiplier problem. Note that if $\eta$ is a polynomial and the
corresponding convolution equation is seen over $\mathcal{E}'$, then
$\eta \delta_{0}\widehat{*}X=\eta^{*} \widehat{X}$. Thus, if we
assume $X,Y$ are holomorphic and of type $0$, then
$\eta^{*}=\widehat{\eta}$. Assume using the Fourier-Borel transform,
that $M=X_{x}+Y_{y} \rightarrow H$ and $W=Y_{x}-X_{y} \rightarrow
G$, then the condition $\widehat{W}=0$ is the condition that $x^{*}
\frac{d \eta^{*}}{d x^{*}}=0$, that is $\eta^{*} \widehat{X}=
\widehat{Y}$. Further, $M_{1}=H_{x}+G_{y}$ and $W_{1}=G_{x}-H_{y}$,
why the condition $\widehat{W}_{1}=0$ is a condition $\gamma^{*}
\widehat{H}=\widehat{G}$, that is $x^{*} \frac{d \gamma^{*}}{d
x^{*}}=0$. For $M_{2}=(\widehat{X})_{x^{*}}+(\widehat{Y})_{y^{*}}$
and $W_{2}=(\widehat{Y})_{x^{*}}-(\widehat{X})_{y^{*}}$, if we
assume $\widehat{X},\widehat{Y}$ are holomorphic of type $0$, then
the condition $\widehat{W}_{2}=0$ is a condition $x \frac{d
\eta}{dx}=0$, that is $\eta X=Y$. Finally,
$M_{0}=(\widehat{H})_{x^{*}}+(\widehat{G})_{y^{*}} \rightarrow C
\frac{d x'}{dt}$ and
$W_{0}=(\widehat{G})_{x^{*}}-(\widehat{H})_{y^{*}} \rightarrow -C
\frac{d y'}{dt}$, give the condition $\widehat{W}_{0}=0$ is a
condition $x \frac{d \gamma}{d x}=0$, assuming
$\widehat{H},\widehat{G}$ are holomorphic. Let $S$ denote the
mapping $(X,Y) \rightarrow (M,W)$ and $T$ the mapping $(M,W)
\rightarrow (M_{2},W_{2})$. If we assume that all elements are
holomorphic over regular approximations, we can prove that the
mapping $T$ preserves order of zero.

\vsp

Consider the following scheme
\begin{eqnarray}
\nonumber
(M,W) \rightarrow & (H,G) \rightarrow & (M_{1},W_{1}) \\
\nonumber
\downarrow_{T} & \downarrow & \downarrow_{T} \\
\nonumber
 (M_{2},W_{2}) \rightarrow & (\widehat{H},\widehat{G})
\rightarrow & (M_{0},W_{0})
\end{eqnarray}
and the corresponding characteristic sets $\widehat{\Sigma}=\{
(H,G)=0 \}$, $\widehat{\Sigma}_{1}=\{
(\widehat{M}_{1},\widehat{W}_{1})=0 \}$, $\widehat{\Sigma}_{0}=\{
(\widehat{M}_{0},\widehat{W}_{0})=0 \}$, $\widehat{\Sigma}_{2}=\{
(\widehat{M}_{2},\widehat{W}_{2})=0 \}$. Let $\Sigma \rightarrow_{T}
\Sigma_{2}$,$\Sigma_{1} \rightarrow_{T} \Sigma_{0}$. Let
$\widetilde{S}(M_{2},W_{2})=(M_{1},W_{1})$, The sets $\Sigma_{1}
\rightarrow \Sigma_{2}$ are connected through $S \FB S=\widetilde{S}
S \FB$ over $(X,Y)$.

\vsp

Assume $g=\widehat{Y}/\widehat{X}-\eta^{*} >0$, then we know that
there exists a measure $v$, non-negative and slowly growing such
that $g=\widehat{v}$. The condition on positivity implies exactness
over the tangent space (global pseudo-base). We can say that $\Top$
maps contingent regions on contingent regions, in the sense that the
order of the regions are preserved, that is the number of defining
functions is preserved. Let $\vartheta_{2}=W_{2}/M_{2}$ and
$\vartheta_{1}=W_{1}/M_{1}$. Then, we have that $\vartheta_{1}$
changes sign as $(\vartheta_{2}+\eta)/(1+\eta \vartheta_{2})$.

\subsection{ Degenerated points for the method }
The problem of determining $x \in L^{2} \cap H$, so that $f=c_{1}x +
c_{2}h(x)$ is trivial over $\{ \eta=const. \}$. Consider a
neighborhood where $\eta$ is quasi conformal, $\mid \eta(x)-x \mid
<c$, locally for a constant $c$. Assume
$\frac{-\widehat{M}}{\widehat{W}}=\frac{xX+yY}{xY-yX}$, where $M=0$
is the equation for mass conservation and $W$ is the vorticity. The
Poincar\'{e} index counts the number of changes of sign $-\infty$ to
$\infty$ and back, for this quotient, why if for instance $\{ \eta +
\vartheta \geq \mu \}$, for a positive constant $\mu$, the index is
zero. For $\eta$ conformal we have $0< \mid \vartheta - \eta \mid$.
Assume $\eta$ algebraic over $\vartheta=\eta$. In a neighborhood of
$\widehat{W}=0$, if $\eta$ is constant, then
$\frac{-\widehat{M}}{\widehat{W}}$ changes sign as $1/(\vartheta +
\eta)$. In the same manner in a neighborhood of $\widehat{M}=0$, the
corresponding test is for $1/(\eta^{-1}+\vartheta)$. Note that conformal mappings do not preserve continuum or reducedness, unless they are bijective.

\vsp

Consider the mapping $dx \rightarrow dv(x) \rightarrow dh(x)$, then
$F(x,v(x))$ have isolated singularities. To describe the
singularities to $F(x,h(x))$, we start with $(M,W)$ and consider
$W=const.$ (not necessarily non-zero) or $W$ regular. We will use
$\widehat{W}=G=Pf$, where $P$ is a polynomial and $f$ is a
continuous function and uniformly bounded (in the real space). When
have $\widehat{W}=0$, equivalently $x^{*} \frac{d \eta^{*}}{d
x^{*}}=0$ $\Leftrightarrow \eta^{*} \widehat{X}=\widehat{Y}$. In the
same manner $x^{*} \frac{d \gamma^{*}}{d x^{*}}=0$ $\Leftrightarrow
\gamma^{*} \widehat{H}=\widehat{G}$ and we can map $\eta^{*}
\rightarrow \gamma^{*}$ by a contact transform (in the case where $\Top$ is
a contact transform). In the case where
the multipliers are polynomials, there is a simple connection $\eta
X=Y$, $\gamma H=G$ and $\eta^{*}=\widehat{\eta}$,
$\gamma^{*}=\widehat{\gamma}$. Locally, where $H \neq 0$, we have
$\mid \gamma \mid \leq C \mid P \mid$, for a constant $C$. Consider
the set $\Sigma_{\gamma^{*}}=\{ \gamma^{*}= d \gamma^{*} \}$ and
outside this set, $G$ is regular if and only if $H$ is regular
$(\neq 0)$. In the same manner for $\Sigma_{\eta^{*}}$.
Consider for $G=Pf$, the set $\Sigma_{f}=\{ f=0 \}$ and $\Sigma_{f}
\subset \{ x^{*} \frac{d \gamma^{*}}{d x^{*}}=0 \}$. Assume $\Omega=\{
\gamma=P=0 \}$ algebraic. Further, either $\Sigma_{\gamma^{*}}$ has
isolated singularities locally or $\Sigma_{\gamma^{*}}=\{ f=0 \}
\backslash \Omega$. We are assuming $H \neq 0$ why $\gamma^{*}=0$
implies $f=0$ while $\gamma^{*}=H=0$ implies $(x,y) \in \Omega$.
Where $H \neq 0$, we have $f=0$ outside $\Omega$. Further, that
$H=0$ corresponds to the proposition on
lineality in the tangent space and if $H \neq 0$, we can only have
lineality in the symbol space.

\vsp

A neighborhood of a singular point $P$, is divided into contingent
regions where $\widehat{W}$ has constant sign. Let $(J_{p}')=\{
(x,y) \quad Y \geq \mu_{p} X \}$, for a constant $\mu_{p}$. Then we
have $(J_{p}')=\{ (x,y) \quad \gamma^{*} \geq \mu_{p} \}$ which, under the conditions above is a semi-algebraic set and on this set, every boundary curve has a local
maximum. In the case where $\gamma^{*}$ is not a polynomial, we can
find semi-algebraic sets $V_{1},V_{2}$, such that $V_{1} \subset
(J_{p}') \subset V_{2}$ ($H \neq 0$). Assume the boundary is given
by a quasi conformal mapping $c$ such that we have existence of a
conformal mapping $c_{1}$ on the boundary with $c \sim_{m} c_{1}$,
$c_{1}(x_{0},y_{0})=0$ and either $\frac{d c_{1}}{dy} \neq 0$ or
$\frac{d c_{1}}{d x} \neq 0$ in $(x_{0},y_{0})$. We can assume that
$c_{1}$ has isolated singularities on the boundary. If we assume $c$
algebraic modulo monotropy at the boundary with $w^{-1}c$ conformal
and $w^{-1}c(x_{0})=0$ and $w$ algebraic. Thus, $c$ has isolated
singularities at the boundary. Assume $c$ satisfies Pfaff's
equation, $dc-pdx-qdy=0$ with $p=\frac{dc}{dx}$ and
$q=\frac{dc}{dy}$. If $p \neq 0$ at the boundary $dy/dx=-q/p$ and
the condition $\mid p \mid \leq k \mid q \mid$ gives in the context
$\mid dy \mid \geq c \mid dx \mid$ at the boundary.

\subsection{The spectrum to $\Top$}

Assuming $y=h(x)$ and $\frac{d \ram x}{d x}=1$ (we can assume $x$ reduced) and $Y_{T}=\chi_{T} X_{T}$ and if $\frac{d h}{d x}$
is locally bounded, the system is well defined. Let $\Top \mbox{:  }(X,Y) \rightarrow (\widehat{X},\widehat{Y})$
and $\Top(X d x + Y d y)=\widehat{X} d x^{*}+\widehat{Y} d y^{*}$, where $d x \rightarrow d x^{*}$
is the Legendre. The discussion has to do about when $\Top$ is analytic. A sufficient condition
for this is that $\Top$ is pure. 

\vsp

Given a formally self-adjoint and hypoelliptic differential operator $L$, the spectrum to a self-adjoint realization in $L^{1}(\R)$,
$\sigma(\mathcal{A})=\{ L(\xi) \quad \xi \in \R \}$, is left semi-bounded.
If we let $\omega=\chi_{(-\infty,\lambda)} \circ L$ and $\mathcal{F}(E_{\lambda} u )=\omega \mathcal{F} u$, for the spectral projection $E_{\lambda}$
and $\mathcal{F}^{-1} w \in \DL$. For Baire functions we have the spectral mapping theorem, that is if $A$ is harmonic conjugation, $\sigma(\Top(A)) \subset \overline{\Top(\sigma(A))}$
and when $\Top$ is algebraic, we have equality.

\subsection{ The spectrum to multipliers }
The spectral theory is set over $\chi X=Y$ and $\chi^{\vartriangle}
\widehat{X}=\widehat{Y}$, that is $\frac{d h(x)}{dt}=\chi \frac{d
x}{dt}$ and $\lambda^{\vartriangle} H=G$, $\lambda\widehat{H}=\widehat{G}$. Let
$\Phi=\{ \lambda=const. \}$ and $\Phi^{\vartriangle}=\{ \lambda^{\vartriangle}=const \}$.
Let $\Psi$ and $\Psi^{\vartriangle}$ be the constant sets to $\chi$ and
$\chi^{\vartriangle}$. Let $V$ be the Hamilton function corresponding to
$\Delta V=-W$. Assume that $\Top \mbox{:} \lambda \rightarrow
\lambda^{\vartriangle}$ is analytic on $\Phi$ and $\Psi$. Assume $A$ the linear
operator corresponding to $A(M)=M^{\diamondsuit}$, that is harmonic conjugation. Thus, if $Y=\eta
X$, then $A(X)=\chi X$ and $A(\widehat{X})=\chi^{\vartriangle} \widehat{X}$. If
$Sp(\Top(A)) \subset \mathbb{R}$, then $Sp(\Top(A))$ is constant
$\subset \Phi^{\vartriangle}$.

\vsp

We will also consider the sets $\{ e^{-\vartheta^{\vartriangle}}
\chi^{\vartriangle}=const. \}$ and $\{ e^{-\vartheta}\chi=const. \}$ as
parabolic Riemann surfaces.  Over these sets $\Top$, when it acts as a
Legendre transform, we can consider it as algebraic and we can apply the spectral
mapping theorem. Note that through symplecticity we have that if
$\chi^{\vartriangle}/\chi=const.$, then this holds in a point. If $\chi$ is
algebraic in $T$, then $\chi^{\vartriangle}$ is algebraic $1/T$. We add the
condition $\Top(0)=0$ to the definition of $\Top$, corresponding to the
condition that $P_{0}=P_{0}^{\vartriangle}$ exists, which is necessary for
analytic continuation.

\subsection{ Conclusions concerning the trace formula }
\label{Nilsson}
In the representation $e_{\lambda}(x,y)=\int_{\Delta(\xi) < \lambda}
e^{i (x-y) \cdot \xi} d \xi$,(cf. \cite{Nilsson63}) a trace in $(x-y)$ corresponds to a
trace in $\xi$, considered as a functional in $H'$. Through
Iversen's result, the correspondent to $\ram$ in $\xi$, is a
function, however multi-valued in most cases. If we reduce the
situation to a real parameter, the continuity of $r_{T}$ means that
there can be no trace in $\xi$ corresponding to the leafs, since the
operator is elliptic. The only possibility is through change of
leafs, but through the conditions that the covering is regular, we
know that these trace sets do not have a measure .

\subsection{Conjugation}
Assume $\chi^{\vartriangle}H=G$ and consider the two FBI-transforms $F_{G}$ and
$F_{H}$, where the kernels are harmonically conjugated. We are assuming
$H=e^{\Phi}$, where $\Phi \geq 0$ over a parabolic surface and where
$\Phi < 0$ implies $\Phi=const.$ If $\chi^{\vartriangle}=1$, we have
$F_{G}=F_{H}$ and we assume $F_{H/G}=F_{1}=\delta_{0}$, the
evaluation functional. Consider the problem of geometric
equivalence. The mapping $\Top : \Psi= x \frac{d \chi}{ dx }
\rightarrow x^{*} \frac{d \chi^{\vartriangle}}{d x^{*}}$, maps the set
$\chi^{\vartriangle}=0$ (or $\chi=const$) on to $\Top \Psi=0$. Further
$\frac{-1}{\chi^{\vartriangle}}=\frac{-H}{G}=\frac{\Top \Psi}{\overline{\Top
\Psi}}$ ( quotient of polynomials ) and we have
${\Top}^{2}\Psi=-\Psi$. We can show that $\Top$ is not algebraic but
if $\Top$ maps $0$ onto $0$, it has the property that $\Top
\Psi=0$ implies $\Psi=0$. In the case where $\chi^{\vartriangle}$ is constant,
we have not excluded the case where $G=0$ ( and $H=0$). These points
are singular and parabolic. Thus, if $\chi^{\vartriangle}=1$, we have $\Top
\Psi + \overline{\Top \Psi}=0$. In the case where $\chi^{\vartriangle}$ is
algebraic, we have that $-H/G$ changes sign as $\Top \Psi/
\overline{\Top \Psi}$, a quotient of polynomials.

\vsp

We are considering two types of conjugates (still referring to kernels),
$F_{H}^{\diamondsuit}=F_{G}$ and $F_{\Psi}^{\vartriangle}=F_{\Top
\Psi}$. We assume as above $F_{1}=F_{\Top \Psi/\overline{\Top
\Psi}}$ corresponds to the proposition $F_{T \Psi}=F_{\overline{\Top
\Psi}}$, where as before $\Top
\Psi=\vartheta^{\vartriangle}-\chi^{\vartriangle}$ and
$\overline{\Top \Psi}=\vartheta^{\vartriangle} +
\chi^{\vartriangle}$. Particularly, if $\Top \Psi=\rho \Psi$ and
$\Top (\rho \Psi)=\rho_{1} \Top \Psi$, we have $\Top \Psi=
-1/\rho_{1} \Psi$. Further, $\frac{x^{*}}{x} \frac{d
\chi^{\vartriangle}}{ d \eta}=\rho \frac{d x^{*}}{dx}$. Now let
$\chi^{\vartriangle}=i \rho^{\vartriangle}$, for a real
$\rho^{\vartriangle}$ and assume
$\overline{\vartheta^{\vartriangle}}=\vartheta^{\vartriangle}$.
Thus, $\overline{\Top \Psi}=S \Psi$, why $S^{2}=\overline{\Top}
\Top$. A differential $\Psi=w+iw^{\diamondsuit}$ is said to be pure,
if $\Psi^{\diamondsuit}=-i\Psi$. We have $\Top(w+i\Top w)=-i(w+i\Top
w)$, $\Top(w+iSw)=-i(-w+i\Top w)$ and $S(w+i\Top w)=-i(-w+iSw)$.
Finally, $S(w+iSw)=-i(-w +iSw)$. Thus, $w+i\Top w$ is pure and if
$w$ is symmetric with respect to the origin, also $w+iSw$ is pure.
If $\Top \Psi=g$ and $g=\alpha \overline{g}$, we have
$\overline{g}=\overline{\alpha} g$ which is pure, if $\alpha=i$. Let
$p=(w,w^{\diamondsuit})$ and consider $r=(w,\Top w)$. Through the
definition of $\Top$, we have $\Top w^{\diamondsuit}=(\Top
w)^{\diamondsuit}$. Thus, if $(w,w^{\diamondsuit})$ is pure, the
same holds for $(w + iw^{\diamondsuit}) + i(w^{\diamondsuit} + iT
w^{\diamondsuit})$ and it follows from the condition $T^{2}=-I$ over
$(w,w^{\diamondsuit})$ that $r$ is pure.

\vsp

Assuming $w$ is pure, we have that $(w,w^{\diamondsuit \diamondsuit})$ is pure. This
follows since $(w+iw^{\diamondsuit \diamondsuit})=w-w^{\diamondsuit}$, which is pure if $w \bot
w^{\diamondsuit}$, where we have used that $i w^{\diamondsuit \diamondsuit}=-w^{\diamondsuit}$ and if
$w^{\diamondsuit}=-iw$, then $(w+iw^{\diamondsuit \diamondsuit})=w+iw$. If $(iw^{\diamondsuit})^{\diamondsuit}=-iw$, then the
form $(w+iw^{\diamondsuit \diamondsuit})$ is pure.

\subsection{ Symplecticity and forms }

Assume $u=a dx+b dy$ and let $F(u)=F(a)dx+F(b)dy$. Then
$F(u^{\diamondsuit})=-F(b)dx+F(a)dy=F(u)^{\diamondsuit}$, assuming
that $F(-b)=-F(b)$. Further, if
$\frac{d}{dx}F_{1}^{\diamondsuit}=F^{\diamondsuit}\frac{d}{dy}$,
then $\Delta F_{1}^{\diamondsuit}=F^{\diamondsuit} \Delta$. Assume
$F^{\diamondsuit}$ a homomorphism and that it maps $0 \rightarrow
0$. We can as before write, $< F(M),
\theta>=\rho_{T}(x,y)<M,\theta>$ and as
$F(M)=-F^{\diamondsuit}(W)=\Delta F_{1}^{\diamondsuit}(V)$, such
that if $\Delta V=0$, then $F^{\diamondsuit}(W)=0$ and $F(M)=0$.

\vsp

 Consider the mapping $\Top \mbox{:}
X dx+Y dy \rightarrow \widehat{X} dx^{*} + \widehat{Y} dy^{*}$. We
are assuming $\gamma$ and $\gamma^{\vartriangle}$ in duality and write $\Top
\tilde{M}(d \gamma)=\tilde{M_{2}}(d \gamma^{\vartriangle})$. $\Top$ is first assumed
an extension of a contact transform in the sense that $\Top \tilde{M}=\rho
\mathcal{L} \tilde{M}$, where $\rho$ is at least a Baire function. Assume $\Top f_{\tau}/f_{\tau}=const.
\Rightarrow \tau=0$ and $d \Top f_{\tau}/d f_{\tau}=const. \Rightarrow
\tau=0$. Assume $\alpha(w, w^{\diamondsuit})$ a symplectic form over $E=V \times V^{\diamondsuit}$ and consider $E \times \Top
E$. Assume $S$ an involutive set with respect to the bracket
$\frac{d V}{dx} \frac{d}{dy}-\frac{d V}{dy} \frac{d}{d
x}+\frac{dV}{d x^{*}}\frac{d}{dy^{*}}-\frac{dV}{dy^{*}}\frac{d}{d
x^{*}}$. Assume $\Top$ equivalent with a form symplectic with
respect to $\alpha$, in the sense that the sets $\{ \chi= \frac{d}{d
x} \chi=0 \}$ and $\{ \chi^{\vartriangle}= \frac{d}{d x^{*}} \chi^{\vartriangle}=0 \}$ are
both minimally defined and equivalent. That is formally,
$\alpha(\Top w,\Top w^{\diamondsuit})=\varrho
\alpha(w,w^{\diamondsuit})$ and $\varrho \sim constant$.

\vsp

Consider the form $(p,q)_{\sigma}=-(q,p)_{\sigma}$, where
$p=(w,w^{\diamondsuit})$ and $q=\Top p$, where we are assuming
$(p,q)_{\sigma}=(\Top p,\Top^{2}p)_{\sigma}=-(q,p)_{\sigma}$, that
is skew-symmetric and bilinear, assuming the double transform in
$Exp_{0}$ is equivalent with $-I$ (after analytic continuation).
Through the conditions
$(w_{\tau})^{\diamondsuit}=(w^{\diamondsuit})_{\tau}$ and the
quotient $\Top (w_{\tau})/(\Top w)_{\tau}$ is never algebraic. We
conclude that $\Top$ under these conditions is symplectic for
$()_{\sigma}$ and that the involutive set $S$ has a corresponding
extended involutive set with respect to $()_{\sigma}$.

\newtheorem{Backlund}{ Proposition}[section]
\begin{Backlund}
 The mapping $\Top$ when planar (pure) preserves analyticity. 
\end{Backlund}

\subsection{ The reflection principle }
 Consider $\Phi^{*}=\{\vartheta^{*}
\quad  e^{-\vartheta^{*}}\chi^{*}=const \quad \exists \vartheta^{*}
\}$ and $\Phi=\{ \vartheta \quad e^{-\vartheta}\chi=const \quad
\exists \vartheta \}$. Consider the Legendre transform $R$,
according to $\mbox{RE} < R(\chi),\chi >-1$. 
\vsp

Further, note that if $\chi^{*}/\chi=f_{0}(\chi)$, where $f_{0}$ has
slow growth like $e^{-\varphi}$ as $0 \leq \varphi \rightarrow
\infty$. Thus, $R(e^{\vartheta})=e^{-\varphi + \vartheta}$ and
$\vartheta^{*}=-\varphi+\vartheta$. Note that if $\Top$ is
considered as a continuous morphism on a Banach algebra $A$ with
$\Top (e^{\vartheta})=e^{\vartheta^{\vartriangle}}$, then $\vartheta^{\vartriangle} \in
A$. For instance if $\vartheta_{T}$ is algebraic in $T$ or a quotient of algebraic functions,
for $T$ close to $0$, then the same holds for $\vartheta_{T}^{\vartriangle}$.The spectrum for $R(f)/f$ contains under the conditions on
$f_{0}$ both $0$ and $\infty$. If $R$ preserves first surfaces, then
we can extend the definition of $\Top$ to $\Top (0)=0$, that is a
part of the boundary. Since $\Top$ is pure if we assume the
corresponding form closed, we have an analytic mapping.

\vsp

Assume that the segments $\gamma,\gamma^{*}$
have a point in common $P_{0}=P_{0}^{*}$. As
$\gamma_{0}=\int_{\Gamma} d \gamma= \int_{\Gamma^{\bot}} d
\gamma^{*}=\gamma_{0}^{*}$, where $\Gamma$ is a closed contour
$\Gamma \sim \Gamma_{0}$. We thus have a reflection principle for
$\Top$ expressed in preservation of the condition on flux. That is
$Flux(\Top W)=\int_{\Gamma} d \Top W$. Staring with the condition
$-H/G \sim \Top W/\overline{\Top W}$, we are assuming $Flux(\Top
W)=Flux(\overline{\Top W})=0$, reflecting symmetry with respect to
harmonic conjugates. Through the condition on parabolic
singularities, there must be a point at the boundary, where $\Top
W=0$ such that $W=0$. If the point is singular for the associated
dynamical system $(M,W)$, the condition must be symmetric, that is
must have $M=W=0$. Thus the condition on flux is necessary for
regularity for this dynamical system. Consider
$\Top(\chi)=\chi^{\vartriangle}$ such that $\gamma_{1} \ldots
\gamma_{N} \sim 1$ and $\Top(\gamma_{1} \ldots \gamma_{N}) \sim 1$,
that is preserves the closed property. The proposition is thus that
$\Top$ preserves flux, but through symplecticity, that is $\Top f=f$
implies $f$ a point, it is not considered to be a normal mapping on
a non-trivial set at the boundary.

\subsection{ Codimension one singularities }
To determine if the symbol corresponds to a hypoelliptic operator, we must
prove that every complex line, transversal to first surfaces, does not contribute to lineality. A complex line is considered as transversal, if it intersects a first surface and the origin. We are assuming that $F$
is the lifting function with $\Delta F = -W$ on a domain $\Omega$
(and the associated equation $\Delta F^{\diamondsuit}=-M$ on
$\Omega$).
The boundary condition is assumed parabolic, that is $\Gamma$ is
such that $P(\delta_{T})F_{T}=0$ implies $T=0$ for a polynomial $P$ and $F_{T}=F$ in
$T=0$.  We are
thus assuming that the parametrix is invariant at the boundary,
$E_{T}F_{0}=F_{0}$ on $\Gamma$. Let $\phi_{T}=E_{T}F$. We are
assuming in this approach that $W_{T}$ is defined by $\Delta
\phi_{T}=E_{T}W+R_{T}$. Assume the operator $\Top$ is continuous at
the boundary, with $\overline{\Top W}=\Top W^{\diamondsuit}$, we
then have $\Top W^{\diamondsuit}=\zeta \Top W$, for $\zeta \in
\mathbb{C}$ and $\mid \zeta \mid=1$. In a Puiseux-expansion, we have
that the coefficient for $t$, $\chi^{*}=a_{0}+a_{1}t+\ldots$ is
$\neq 0$ and $t=\surd \zeta$. Thus, the order of the critical
surface is one and we have singularities of codimension $1$. When
$\Top$ preserves the order of contact, we have the same conclusion
for the multiplier $\chi$.

\vsp

A co-dimensional one variety $S(p)$, is such that $S(p,x)=0$ and
$s_{x} \neq 0$, where $p=s_{x}$ is a characteristic variety, if
$g(x,s_{x})=0$, $p.dx=0$ and $g_{x}dx + g_{p}dp=0$, $dS=pdx$. As
before, we have $\ram ds=(\ram p)dx$. $S(p)$ is involutive for $g,
\gamma$, if $H_{g}(\gamma)=0$ (the Poisson bracket). For parabolic
singularities at the boundary, we are considering isolated
singularities in higher order derivatives $\psi_{N-1}=\frac{d^{N-1}
\psi}{d T^{N-1}}=0$ and $\psi_{N} \neq 0$. Assume $g$ is defined
through $dS=\frac{d g}{d x} dx + \frac{d g}{d y} d y+ \frac{d g}{d
x^{*}} d x^{*}+ \frac{d g}{d y^{*}}d y^{*}$ $=-Ydx+X dy- \widehat{Y}
d x^{*}+ \widehat{X} dy^{*}$, that is $dS=W^{\diamondsuit}+\Top
W^{\diamondsuit}$. Symplecticity gives that $d \gamma_{T}^{*}/d
\gamma_{T}=const. \Rightarrow T=0$. We have a canonical symplectic
form, $d< p,dx>=0$ on $S(p)$, where $p=d \varphi / dx$ and $d( d
\varphi)=0$.  Note that we are assuming $\Top : x \frac{d \chi}{d
x} \rightarrow x^{*} \frac{d \chi^{*}}{d x^{*}}$ and we must assume
$x \frac{d \chi}{d x} > 0$ (Poincar\'e-index) implies $x^{*} \frac{d
\chi^{*}}{d x^{*}} > 0$. If $\Top$ is algebraic, the spectral
theorem can be applied with advantage and we must assume that $\Top$
does not have zero's on the boundary, that is $x^{*}=\frac{d \chi}{d
x}$ and $< \frac{d \chi}{d x},x>=<x^{*},\frac{d \chi^{*}}{d x^{*}}>
\neq 0$. These sections correspond to $\Top : <x^{*},x> \rightarrow
<x,x^{*}>$ (normal sections). We are not assuming the trajectories
in a reflexive space, but it is sufficient to consider the Lagrange
case $\Gamma^{\bot}=\Gamma$.

\subsection{ The mapping $\Top$ and parabolicity }

Consider $\Top$ over a set where it is algebraic in $H'$, for instance
Legendre. We have seen that there can be no closed contour in the
infinity. On the other hand, if $\Top
e^{\vartheta}=e^{\vartheta^{*}+\varphi^{\vartriangle}}$, where $\vartheta^{*}$
is related to $\vartheta$ through a Legendre transform, we have that
$e^{\varphi^{\vartriangle}}$ can define a circle in the infinity. Further,
if $\vartheta^{\vartriangle}_{T}=P(\frac{1}{T})\vartheta^{*}$ and
$\vartheta_{T}=P(T)\vartheta$, taking the closure of the domain
means that
$\vartheta^{\vartriangle}_{T}=e^{\alpha(\frac{1}{T})}\vartheta^{*}$ and
that $\vartheta^{\vartriangle}_{T}/\vartheta_{T}=e^{\alpha(
\frac{1}{T})-\alpha(T)}\vartheta^{*}/\vartheta$ and we may well have
that $\alpha(\frac{1}{T})-\alpha(T) \rightarrow 0$ as $T \rightarrow
\infty$ simultaneously as $T \rightarrow 0$. 

\vsp

We are assuming that $\Top$ preserves the parabolic property for the
stratification, that is it maps exponentials on exponentials.
Note
that we have a parabolic approximation if and only if for all
functions $u$ harmonic in a neighborhood of the ideal boundary, with
finite Dirichlet integral, we have vanishing flux (cf.
\cite{AhlforsSario60}). Assume $\Top$ maps finite Dirichlet
integrals on finite Dirichlet integrals. We have defined $\Top$ such
that $\{ W=0 \} \rightarrow \{ \widehat{W}=0 \}$, why $\Delta F=-W
\rightarrow \Delta^{\vartriangle} \Top F =-\Top W$ and we see that $\Top$
preserves parabolic approximations.
Consider $\Top_{1}$ algebraic, in the sense that
$\Top_{1}(e^{-v}\chi)=\Top_{1}(e^{-v}) \Top_{1}(\chi)$. We then have
that $\Top_{1}$ maps constants on constants, why if $\chi=e^{v}$, we
have that $\Top_{1}(e^{-v} \chi)=const.$ and
$\Top_{1}(\chi)=\Top_{1}(e^{v})=e^{v^{*}}$. Assume that $\Top(e^{-v}
\chi)=e^{\varphi^{\vartriangle}}$, that is $\Top$ maps constants on to
exponentials and that $\Top_{1}(e^{-v} \chi)=e^{\varphi^{\vartriangle}}$. Thus,
if $\chi^{\vartriangle}=\Top_{1}(\chi)$, we have that
$\chi^{\vartriangle}=e^{v^{*}+\varphi^{\vartriangle}}$. If $\Top=\Top_{1}$, there is no
room for a closed contour, in the infinity. If we assume $\Top \sim
\Top_{1}$, in the sense that $e^{\varphi^{\vartriangle}}=\Top(e^{-v}
\chi)=\Top_{1}(e^{-v} \chi)$ along two different paths, the result
depends on the property of monodromy for the stratification.
That is if $\Top$ is locally injective with respect to path, we can write $\Top_{1} \Top \sim I$.

\subsection{The vanishing flux condition in phase space}
Consider the linear functional $L(\phi)=\int_{\beta} d \phi$ and consider the difference $L(e^{\phi})-e^{L(\phi)}=\big[ \widehat{I},L \big]-
\big[ \widehat{L},I\big]$. We note that the vanishing flux condition $L(\varphi^{\diamondsuit})=0$ does not imply $L(e^{\phi^{\diamondsuit}})=0$.
Assume that $(\beta)$ is a neighborhood of the origin and note
\begin{equation} \label{Green}
\int_{\beta} \tilde{W}^{\diamondsuit}=\int_{(\beta)}W dx dy
\end{equation}
where $\tilde{W}^{\diamondsuit}=P dx+Q dy$. Immediately, we note that $\widehat{L}(\phi)=L(e^{\phi})=\int_{\beta} d e^{\phi}$ bounded in the infinity
implies that $\widehat{L} \sim P(\frac{1}{T})$, as $T \rightarrow \infty$. Note that if the indicator for $L$ is $\alpha$, then the functional $L$ has 
support on a ball of radius $\alpha$. Starting with (\ref{Green}), if $\tilde{W}_{T}$ is algebraic in $T$, then the measure for $(\beta)$ is zero, that is
$L(e^{\phi})$ is of real type. Note also that $\int_{\beta^{*}} d \Top W=0$ does not imply $\int_{\beta} d W=0$, however $\int_{\beta^{*}}d \overline{\Top W}=0$.
Particularly, if $e^{\phi}=e^{e^{\alpha}}$, we can consider $\big[ L, \big[ \widehat{I},\widehat{I} \big] \big]=\big[ \widehat{I}, \big[ \widehat{I},L \big]$. For instance
if $\int_{\beta} d \phi= -\infty$ then $L(e^{\phi})=e^{L(\phi)}$ implies $L(e^{\phi})=0$. Assume $\tilde{\phi}_{T}$ algebraic in $T$, then $Flux(\phi_{T})=0$.
In the same manner if $\phi_{T}$ is harmonic, then $Flux{\phi_{T}}=0$, further if $\tilde{\phi}_{T}$ is algebraic in $T$, then $\int_{\beta} d \phi_{T}^{\diamondsuit}=0$
implies $(\beta)$ has measure zero. Thus, if $L(e^{\phi_{T}})=0$ implies that the measure of $(\beta)$ is zero.
We note the following result. Assume $\tilde{\phi}_{T}$ algebraic in $T$, then $d \tilde{\phi}_{T},d \tilde{\phi}_{T}^{\diamondsuit}$ are closed, which implies
that $\phi_{T}$ is harmonic, why we have a real type operator.
Note that if $\phi_{T}$ is harmonic on a disc, the mean is constant ($\equiv -\infty$) and the measure for the ideal boundary $(\beta)$ is zero.

\section{Boundary conditions}
\subsection{Introduction}
In the model, the singularities on the first surfaces to the symbol are mapped on to the boundary of the stratification, which is parabolic or more generally very regular. Hypoellipticity is a condition on behavior for the symbol in the infinity but the method using $\Top$ (basically the projection method) requires a discussion on the simultaneous behavior at the boundary. The boundary to the strata is defined by $\{ \M_{N}(f)=f \}$
but we also discuss the phase correspondent $\{ \M_{N}(\log f)=\log f \}$.   
\subsection{ The $\overline{\delta}$-Neumannproblem }
 We will deal with the following problem, given a regular
approximation of a singular point in the boundary, $\Gamma$,
determine $F_{T}$ such that
\begin{itemize}
\item[i)] $F_{T}(\gamma)(\zeta)=f(r_{T} \zeta)$
\item[ii)] $F$ is holomorphic in $\gamma$ and algebraic in $T$.
\end{itemize}
We assume here that the boundary is finitely generated and in
semi-algebraic neighborhood. The first part of $ii)$ is the lifting
principle for a semi-algebraic domain. For $i)$ we note that as $ f
\in (I)$, a finitely generated ideal, there is no problem to
determine $\gamma$, such that the formal series for $F_{T}$
converges. Consider now the second part of $ii)$. Given a regular
approximation $U_{T}$ with algebraic dependence of the parameter
$T$, we can use the $\overline{\delta}-$ Neumann problem to
determine a lifting function $F_{T}$ such that $\delta_{T}
F_{T}=\delta_{T} U_{T}$ and $F_{T}=U_{T}+L_{T}$. We assume the
boundary finitely generated and we can in a suitable topology assume
that given $F_{T}$, there is a domain of $\gamma$ such that the
dependence of $T$ is as prescribed and $F_{T}(\gamma)=f(r_{T}
\zeta)$. The domain in $\zeta$ is a neighborhood of the first
surface generated by one variable $T$ and it is pseudo convex.

\vsp

Note the symmetry condition that if for $T \in V$, $\big[ L \ram
\big]^{*}=\ram^{*} L^{*}$ is algebraic in $T$  corresponding to a real coefficients polynomial in $\overline{T}$, we have that $\ram - c_{T} I$
holomorphic, that is $\ram$ is holomorphic modulo monotropy.  This means that for an
algebraic dependence of $T$, $\delta_{T} F_{T}$ can have level
surfaces. 

\subsection{ Multipliers } As the mapping $x \rightarrow x^{*}$
preserves order of zero, we see that if $x$ is locally reduced, then
the same holds for $x^{*}$. Consider the system of invariant curves
$\{ C_{j} \}=\{ (x,h^{j}(x)) \}$. Assume
$\eta_{j}(x)=\frac{h(x^{j})}{x}$ such that $\eta_{j}(x)=x^{j-1}
\eta(x)=x \eta_{j-1}(x)$. Assume also that we identify using
monotropy, the curves $\{ C_{j} \} \sim_{m} \{ (x,h(x^{j})) \}$. Assume $<\eta
X,\widehat{\varphi}>=<\eta^{*}
\widehat{X},\varphi>=C<\widehat{X},\varphi>$. Thus, $\{ \eta=const
\} \rightarrow \{ \eta^{*}=const \}$. The condition
$\eta_{2}^{*}=const.$ $\Leftrightarrow \eta^{*}=const./x^{*}$ and so
on. Assume existence of an algebraic homomorphism $w$, such that
$\mid w^{-1} \eta^{*}-\frac{1}{x^{*}} \mid< \epsilon$ as $\mid x^{*}
\mid \rightarrow \infty$. Then $\eta^{*}(x^{*}) \sim_{m}
w(\frac{1}{x^{*}})$ as $\mid x^{*} \mid \rightarrow \infty$ and that
$\eta^{*}$ preserves constant value in the $x^{*}-$ infinity, why
the projection method can be applied to $\eta^{*}$. It is of no
significance what level surface we start with, that is
$\eta_{j}^{*}=const.$ $\Rightarrow \eta^{*}$ preserves a constant
value in the infinity.

\newtheorem{proj}{Lemma}[section]
\begin{proj}
Given a system of invariant curves $\{ (x),h^{j}(x) \}$ such that
$h(x^{j}) \sim_{m} h^{j}(x)$, we have that $\eta^{*}$ preserves a
constant value in the infinity and the projection method can be
applied.
\end{proj}

Assume existence of a finite $j$ such that for a lifting function
$F_{T}$, $\frac{d^{j}}{d T^{j}}F_{T}$ is algebraic in $T$, that is
$F_{T}^{j} \sim  \alpha_{j}(T)F$, in a neighborhood of $T_{0}$,
where $\alpha_{j}$ is a locally defined polynomial. Assume
$F_{T}^{j}(x,y)=G(e^{x},y) \rightarrow 0$ as $\mid x \mid
\rightarrow \infty$ for a $G \in H'$ and $y=h(x)$ finite in modulus.
Thus, $G(e^{x},h(e^{x}))=G(e^{x},e^{\tilde{h}(x)})$ is the
representation we prefer. We have also assumed that $F^{j}_{T}$
preserves a constant value in $x$ as $\mid x \mid \rightarrow
\infty$ and $y$ finite. If $\frac{d}{dT}F(\gamma_{T})=F_{1}(\frac{d
\gamma_{T}}{dT})$ and $\theta_{T} \in T \Sigma$, we see that if
$F_{1}$ is algebraic, then
$\theta_{T}=F^{-1}\frac{d}{dT}F(\gamma_{T})$ for a $\gamma_{T} \in
\Sigma$. Further, $<\gamma_{T},\theta_{T}>_{1}=0$ implies
$\gamma_{T} \in \{ B(\gamma_{T})=\mu \}$, for a constant $\mu$ and
$B$ algebraic. Assume $B$ self-transposed and such that
$B(F(\gamma_{T}))={}^{t}FB(\gamma_{T})$. A sufficient condition for
$F$ to map $\Gamma^{\bot} \rightarrow \Gamma^{\bot}$, given that
$\Gamma$ is symmetric, is that it maps constants on constants.

\vsp

For a finitely generated boundary, we have the following result.
Assume the singularity at the boundary, described by $T$, that is
$F_{T}-\mu \rightarrow 0$, for a small constant $\mu$, such that
$F_{1/T}$ is close to a polynomial as $T \rightarrow 0$.  We are assuming $\gamma_{0}$
fix at the boundary and $\gamma_{T} \in$ a neighborhood of
$\gamma_{0}$ with $F_{T}(\gamma_{0})=F(\gamma_{T})$ and where $\mid
w^{-1}F_{T}-1/T \mid < \epsilon$, for an algebraic homomorphism $w$,
as $T \rightarrow 0$. Assume $F_{T}$ invertible over regular
approximations, such that $\frac{d \gamma_{T}}{dT}=F^{-1}\phi_{T}
\neq 0$ and $\frac{d^{2}}{dT^{2}}
\gamma_{T}=\frac{d}{dT}F^{-1}\phi_{T} \neq 0$. Thus, if $F_{1}$ maps
constants on constant, we have that $F$ maps singular points on
singular points and regular points on regular points.

\subsection{ The orthogonal to the boundary }

Assume $F$ preserves a constant value as $\mid x \mid \rightarrow
\infty$ and $\mid y \mid \rightarrow \infty$. Assume $\FB(\eta
X)=\eta^{*} \FB X$. Degenerate points are then on the form $y^{*} d
x^{*}-x^{*}d y^{*}=0$. If $F$ is not algebraic, we are at least
assuming that $F^{-1}$ maps constants on constants or that $F
\sim_{m}$ an algebraic function close to the boundary. Let
$<y^{*},y>_{1}=\mbox{ Re }<y^{*},y>-1$, for $y \in \Gamma$ and $y^{*} \in
\Gamma^{0}$. Then if $<y^{*},y>_ {1}=0$, $\Gamma^{0}$ is a line if
$\Gamma$ is a line. We write $y^{*} \bot y$ for $<y^{*},y>_{1}=0$
and we define $\mbox{bd} \Gamma^{0}$ as the set where this relation
holds over $\Gamma$. If $\Gamma^{00}=\Gamma$, then $\Gamma^{0}
\subset \Gamma$. If $<y^{*}, \cdot>_{1}$ is reduced, we have
isolated singularities at the boundary of $\Gamma$. If $<y,y>_{1}=0$
for all $y \in \Gamma$, we have that $\Gamma \subset \mbox{ bd
}\Gamma^{0}$. Since $y \rightarrow y^{*}$ is a contact transform, we
have $N(\Gamma) \subset N(\Gamma^{0})$. Further, assuming $y^{*}
\rightarrow y^{**}$ is a contact transform with $y^{**} \sim_{0} y$,
we get $rad(\Gamma) \sim rad(\Gamma^{0})$ (equivalence in sense of
ideals). Consider with these conditions $\Gamma_{1}=\{y \quad
<y^{*},y>_ {1}=0 \quad y^{*} \in \Gamma^{0} \}$, then also
$rad(\Gamma_{1}) \sim rad(\Gamma_{1}^{\bot})$. Note that $(T
\Sigma^{\bot})$ does not completely describe the micro-local
contribution.

 We have a few immediate results. Let $\Gamma_{T}=\{ \gamma_{T}
 \quad \gamma_{T}=r_{T}' \gamma \quad \gamma \in \Gamma \}$ and assume the boundary condition
(\ref{bdcond})

 \newtheorem{Lagrange}{Proposition}[section]
\begin{Lagrange}
$\ram-1$ is locally algebraic if and only if $\Gamma_{T} \sim \Gamma_{T}^{\bot}$.
\end{Lagrange}
Proof:\\
We are assuming $\ram \gamma \bot \gamma$ for $\gamma \in \Gamma$,
that is $< \ram \gamma,\gamma>_{1}=0$ for $\gamma \in \Gamma$
implies $\Gamma_{T} \subset \Gamma^{\bot}$. Assume further the
ramifier symmetric, in the sense that
$<r_{T}'x,y>_{1}=<x,r_{T}'y>_{1}$. If $<\eta^{*},\gamma>_{1}=0$ for
$\gamma \in \Gamma$, we have existence of $\eta \in \Gamma$ such
that $< \ram \eta,\gamma>_{1}=0$ why $\Gamma^{\bot} \subset
\Gamma_{T}$. Thus, $\Gamma_{T} \sim \Gamma^{\bot}$. Assume
$\gamma_{T} \in \Gamma_{T} \Leftrightarrow \gamma_{2T} \in
\Gamma_{T}$, then the first implication follows. The converse
implication is obvious. $\Box$

\vsp

The micro-local contribution that is given by $\{ \ram -1=0 \}$ is
thus a subset of the contribution given by $\ram -1$ algebraic. We
claim that it is necessary for the application to
pseudo-differential operators, to assume the Lagrange case, $\Gamma \sim \Gamma^{\bot}$.

\subsection{Treves' curves}

We define $<\gamma_{T},\theta_{T}>_{1}=\mbox{ Re }
<\gamma_{T},\theta_{T}>-1$. Thus, $<F(\gamma_{T}),i \mbox{ Im }
\theta_{T}>_{1}=<\mbox{ Re } F(\gamma_{T}),\overline{\theta_{T}}>$ and if
$\theta_{T} \in T \Sigma$ has the property that $\theta_{T} \in T
\Sigma \Leftrightarrow \overline{\theta_{T}} \in T \Sigma$, we have
that $\mbox{ Re } F(\gamma_{T}) \bot \mbox{ Re } \theta_{T}$ implies $i 
F(\gamma_{T}) \bot \theta_{T}$. Conversely, if $T \Sigma$ is
symmetric with respect to the origin, we have existence of
$\gamma_{T}$ such that $\mbox{ Re } F(\gamma_{T}) \bot \theta_{T}$ implies
$F(\gamma_{T}) \bot i \mbox{ Im }  \theta_{T}$ and analyticity for $\mbox{ Im }
F(\gamma_{T})$ only means that $<F(\gamma_{T}),\theta_{T}>$
 $=<F(\gamma_{T}),\overline{\theta_{T}}>$. Assume $\frac{dF}{d \gamma}
\mbox{:} \Gamma^{\bot} \rightarrow \Gamma^{\bot}$, such that
$<\frac{dF}{dT}(\gamma_{T}),\theta_{T}>$ $=0$ implies $<\frac{d
\gamma_{T}}{dT},\theta_{T}>=0$. Over $\{ F_{T}=\frac{d}{dT}F_{T}
\}$, we thus have that if $F_{T} \bot \theta_{T}$, we have existence
of Tr\`eves curves. Conversely, given existence of $\gamma_{T}$ such
that $<\frac{d \gamma_{T}}{dT},\theta_{T}>_{1}=0$ implies $<
\frac{d}{dT} \mbox{ Im }  F(\gamma_{T}),\theta_{T}>_{1}=0$, why $<
\frac{d}{dT}F(\gamma_{T})-\frac{d}{d
\overline{T}}F(\gamma_{T}),\theta_{T}>_{1}=0$ which is always true
for real $T$. We conclude as has been noted before that the condition that $F_{T}$
is analytically hypoelliptic does not imply that the real and imaginary parts are analytically hypoelliptic.

\vsp

Note that it is possible to have $(dI)$ has a global pseudo-base,
when the pseudo-base for $(I)$ is only local. However,

\newtheorem{globalpsbase}{Proposition}[section]
\begin{globalpsbase}
If $(J)$ is a finitely generated ideal of Schwartz type topology and
with a compact ramification, such that $r_{T}' \phi/\phi \sim_{m} 0$
in the $\zeta-$ infinity. Then, $(I) \sim (r_{T}' J)$ has a global
pseudo base.
\end{globalpsbase}

Proof:\\
We are considering $(I) \sim (r_{T}' J)$, where $(J)=\mbox{ ker }h$
and as before $\eta(\phi)=h(\phi)/\phi$. Through the conditions, we
can satisfy $\mid r_{T}' h(\phi) \mid < c + \mid h(\phi) \mid$, for
a constant $c$ and for $\phi \in (J)$, so $r_{T}'$ is quasi
conformal and $(I)$ is finitely generated, if $(J)$ is. Given that
$h$ is algebraic and such that $h^{N}=1$, we can show that $h(r_{T}'
\phi) \sim_{m} 0$ and by choosing $\phi$ reduced, we have that $\eta
\sim_{m}0$ over $(I)$. If $h$ is analytic, we assume locally $h$ is
monotropic to an algebraic homomorphism. Thus, we can find an entire
function $\gamma$ such that $\gamma \sim_{m} \eta$. $\Box$

\subsection{ Analytic set theory }

Starting with the boundary condition in a higher, finite order
derivative $F_{T}^{(j)}=const.$ implies $T=0$. Consider the sets
$\Sigma_{1}=\{ \zeta_{T} \quad F_{T}=const \}$, $\Sigma_{2}=\{
\zeta_{T} \quad F_{T}=const., F_{T}^{(1)}=const. \}$ and so on. This
gives a finite sequence $\Sigma_{j} \downarrow \{ T=0 \}$. We can
form the corresponding ideals in $\gamma$, such that $N(I_{j})=V_{1}
\cup \ldots \cup V_{j}$, where $V_{1}=\{ \zeta_{T} \quad F_{T}
\mbox{ not constant} \}$. If we assume algebraic dependence of the
parameter for $F_{T}-const.I$ and that we have a neighborhood of
$\zeta_{T}$ that is a domain of holomorphy. Then the sets $V_{j}$ as
geometric complements of algebraic sets, are analytic. We have noted
examples where $V_{1}  \nsubseteq V_{2}$. We also note the following
example, assume $F_{T}=\alpha_{T}/\beta_{T}$ such that
$\alpha_{T}'=\gamma_{T} \beta_{T}'$, where $\gamma_{T}$ is assumed
non-constant and regular holomorphic (not-Fuchs equation), then if
$V_{1},V_{2}$ are analytic, then since $V_{2} \nsubseteq V_{1}$, the
inclusion $V_{1} \subset V_{1} \cup V_{2}$ is strict and we have for
the corresponding ideals $I_{2} \subset I_{1}$. For a parabolic
approximation, the set $\{ \zeta_{T} \quad \vartheta_{T}(\zeta) >0
\}$ is the geometric complement to a first surface, which is with
the conditions above an analytic set. Note that we may still have
that the sets $V_{1}^{\bot}$ are first surfaces.

\vsp

 A different argument can be given using Tr\`{e}ves-curves.
Assume $\Omega$ a domain of holomorphy with $\frac{d}{dT}\gamma_{T}
\neq 0$ on $\Omega$. If $\int_{\Omega} (\frac{d}{dT}\gamma_{T})
\theta_{T} dT=0$ and if we assume the integrand holomorphic, we must
have $\theta_{T}=0$ on $\Omega$, assuming it of positive measure.
Assume $\mbox{bd} \Omega$ on one side locally of a hyperplane, then
we have that $\gamma_{T} \neq
0$ on $\Omega$. Assume now $\Delta \subset \Omega$ where $\Delta$ is
algebraic. Then, $\gamma_{T}$ has no zero's on $\Delta$, but is not
constant. We have assumed that constant functions are not
holomorphic and we must also assume that they are not algebraic. For
instance the complement to a first surface in a domain of holomorphy
is not necessarily an analytic set. Note further that if $\Delta$ is
algebraic, we can assume $\Delta^{\bot}$ is not algebraic and with
the conditions under hand, it must be a first surface. Thus, we have
that $\Omega \backslash \Delta$ is analytic and simultaneously
$\Delta^{\bot}$ is a first surface.

 \subsection{A Tauberian problem}

 Let $V_{N}=\{\zeta_{T} \quad e^{\M_{N}(\eta_{T})}\quad \M_{N}(\eta_{T})  \quad \text{ subharmonic } \}$, which is a subset to $V_{N+1}$. Assume $V$ the set corresponding to $\eta$ subharmonic. Consider the complement set in $V_{N}$,
 $V^{c}=\{ \zeta_{T} \quad \eta > \M(\eta) \}$, then $V \subset V_{N}$,
 through the conditions and $V^{c}$ is analytic, if $V$ is analytic.
  We have that
 $\log X_{1} \sim I_{1}=\{ \eta \geq 0 \}$,  In the same manner we consider $I_{2}=\{ \eta_{2} \geq 0 \}=\{ \eta + \M(\eta) \geq 0 \}$,
 $\ldots ,I_{N}=\{ \eta_{N} \geq 0 \}=\{ \sum_{j=0}^{N-1} \M_{(j)}(\eta) \geq 0 \}$.
 Associated to these ideals, we consider $J_{1}=\{ e^{\eta} \quad \eta \geq 0 \}$
 and so on and we have if $e^{\eta_{N}}$ is analytic, that $N(J_{N})$ contains a path $\zeta_{T}$ that
 is continuous.
 The proposition is thus that given a Tauberian condition, we have existence of
 a continuous approximation of a singular point. For instance if we have that $J_{N}$
is defined by an analytic function and the set $V^{c}$ above is analytic, then we have existence of a continuous path in $N(J_{1})$.

\subsection{ The distance function }

If $e^{g}$ represents the distance to isolated (essential)
singularities, all situated on a finite distance from the origin,
then this distance function is globally reduced. For a holomorphic
function $u$, bounded in the infinity, we must have that the
distance to essential singularities is finite. It is sufficient to
consider points $P$ in a punctuated neighborhood of the origin. For
a harmonic function $u$, we have that it is bounded in the finite
plane, and we only have to apply Phragm\'en-Lindel\"of's theorem.
Consider the representation $u=e^{g+m_{1}}$ harmonic, where $m_{1}$ is
symmetric. Over a parabolic approximation where $-g-m_{1}$ is
subharmonic, we assume $m_{1} \rightarrow 0$ close to the boundary
$\Gamma$. Assume now $d$ globally reduced and that $d \rightarrow
0$, as $P \rightarrow P_{0} \in \Gamma$. Further, $\frac{1}{d}(P)
\rightarrow \infty$, as $P \rightarrow \Gamma$ and $d(\frac{1}{P})
\rightarrow \infty$, as $P \rightarrow \Gamma$. Then we can find
$\epsilon$ small such that $\mid d(\frac{1}{P})-\frac{1}{d}(P) \mid
< \epsilon$, as $P \rightarrow \Gamma$. If all singularities for $u$
are at a finite distance from the infinity, we have that $d(P)
\rightarrow 0$, when $P \rightarrow \Gamma$ over the set where $\frac{1}{d}(P)>0$ and on this set it is clear that in the
limit $P=P_{0}$. The singularities at the boundary are assumed
removable. Assume $P_{0} \neq 0$ and that, for instance $d(P) \sim
\mid P-P_{0} \mid$, then $d(\frac{1}{P}) \sim \mid
\frac{1}{P}-\frac{1}{P_{0}} \mid \rightarrow 0$ implies that
$P=P_{0}$, using reducedness for $d$ and we can conclude that the
Dirichlet problem $\Delta u=0$ on a set (it is sufficient to
assume parabolic) with boundary value $\frac{1}{d}$ is solvable,
modulo monotropy.

\newtheorem{DP}{Proposition }[section]
\begin{DP}
Assume the boundary holomorphic and only with parabolic
singularities, then there is a regular approximation of a singular
point that will reach the point.
\end{DP}

\subsection{ Localization at the boundary }
Assume $P(\delta_{T})$ is the operator used to define the boundary
condition, such that $P^{2}(\delta_{T})$ is hypoelliptic.
$C_{T}(\phi)=P(\delta_{T})\phi-\phi P(\delta_{T})$, where $\phi$ is
a real test function and where $P$ is assumed such that $P^{2}$ is
hypoelliptic. Thus, $P(\delta_{T})(\phi F_{T})=\phi P(\delta_{T})
F_{T} + C_{T}(\phi)$. If $P$ is hypoelliptic, we have that $C_{T}
\prec \prec \mbox{Re} P_{T}$. Otherwise, we will assume that
$P^{2}(\delta_{T})(\phi)-P(\delta_{T}) \phi P(\delta_{T})$ $+
P(\delta_{T}) \phi P(\delta_{T})- \phi P^{2}(\delta_{T}) \sim \mbox{Im}
P^{2}(\delta_{T}) \prec \prec \mbox{Re} P^{2}(\delta_{T})$. Thus, if
$P(\delta_{T}) C_{T}(\phi)+ C_{T}(\phi) P(\delta_{T}) \sim \mbox{Im}
P^{2}(\delta_{T})$. As $T \rightarrow 0$, we have that $C_{T}(\phi)
\rightarrow 0$ (we assume $\phi=1$ at the boundary). Using
Nullstellensatz, that is $\frac{PC_{T}+C_{T}P}{P^{2}} \rightarrow
0$, as $1/T \rightarrow 0$ implies $\mid PC_{T} + C_{T}P \mid <
\epsilon$, for large $T$ ( and real). Let
$(PC_{T})^{*}=C_{T}^{*}P^{*}$ and
if $P^{*}=P$ implies $C_{T}^{*} \sim -C_{T}$.
Symmetry with respect to ${}^{*}$ gives $C_{T} \prec
\prec \mbox{Re} P$ in the infinity, for $P$ such that the square is
hypoelliptic.

\vsp

 Assume now that the boundary condition
is given by a differential operator (reduced) $P(x,\frac{d}{dT})$
such that there is a function $g_{N} \in L^{1}$ in the parameter
close to the boundary, ($g_{N}=M_{N}(f)$) with
$P(x,\frac{d}{dT})g_{N}=I-r_{N}$, where $r_{N}$ is regularizing as a
pseudo differential operator. If we regard $g_{N}$ as an operator
$L^{1} \rightarrow \mathcal{D}_{L^{1}}'$, we can construct $g_{N}$
as an operator with kernel $G_{N} \in \mathcal{D}_{L^{1}}'$, that is
a parametrix, $g_{N}(\phi)=\int G_{N}(x,y)\phi(y)dy$. Given a
parabolic boundary condition, we can assume $P(x, \frac{d}{dT})g_{N}
\in \mathcal{D}_{L^{1}}'$ for $x$ in a neighborhood of a point
$x_{0}$ at the boundary and with $\mbox{ sing supp } g_{N} =\{ x_{0}
\}$. Assume $\phi$ a regular approximation of the singular point
such that $P(x,\frac{d}{dT})\phi=0$ implies $x=x_{0}$. For the
parametrix, we then have $< G_{N},P(x, \frac{d}{dT})\phi>=0$ implies
$x=x_{0}$ (modulo regularizing action), that is $P(x,\frac{d}{dT})G_{N}
\bot \phi$ implies $x=x_{0}$

\subsection{ Further remarks on the boundary }

If $\ram$ is an algebraic homomorphism, then $\ram e^{\phi}=e^{\ram
\phi}$, and consequently $\int_{(I)} \ram e^{\phi} dz(T)=\int_{(I)}
e^{\ram \phi} dz(T)=0$, implies $m(I)=0$ (measure zero set), using a
result by Hurwitz. If $\big[ e^{v^{\vartriangle}_{T}} -1 \big]$ is holomorphic,
we have either that $e^{v_{T}^{\vartriangle}} \equiv 1$ or $m(I)=0$.

\newtheorem{Hurwitz}[DP]{ Lemma }
\begin{Hurwitz}
Assume $e^{v_{T}^{\vartriangle}}-1$ is holomorphic in the parameter $T$ with
$v^{\vartriangle}_{T}$ algebraic in $T$. We then have that $\int_{(I)}
(e^{v_{T}^{\vartriangle}} -1)dz(T)=0$ implies $m(I)=0$.
\end{Hurwitz}

 \vsp
 
 We also note the following consequence of the condition on
 vanishing flux, $\int_{(C)} M d x d y=$  $\int_{C}
 \tilde{M}^{\diamondsuit}=0$, means that there is a trajectory
 $\gamma$ such that $M(\gamma)=0$ in $(C)$. Particularly, if $f$ is
 such that $\int_{S} \M_{N}(d \tilde{f})^{\diamondsuit}=0$ that is we have
 $\int_{(S)} \Delta \M_{N}(f) dx dy=0$, we must have that $\M_{N}(f)$
 changes sign in points inner to $(S)$.

\subsection{A very regular boundary }

The boundary is said to be very regular, if the singularities are
located in a locally finite set of isolated points or segments of
analytic curves (cf. \cite{Parreau51}). Thus, we are assuming that
if $f_{0}$ is a boundary element, then a very regular representation
of the boundary preserves the locality of singularities, but not
necessarily the order. Assume $\Gamma=\{ \Gamma_{j} \}$ is a locally
finite set of analytic curves, where the set of common points is a
discrete set. Given an element in $(\Bm)' \subset \DL$, we know for
the real Fourier transform, that $\widehat{f}=P(\xi) f_{0}$, where
$f_{0}$ is a continuous function in the real space and $P$ a
polynomial. Extend $f_{0}$ to a continuous function in a complex
neighborhood of the real space and denote $\tilde{f_{0}}$ the
function such that $\widehat{\tilde{f_{0}}}$ is the extended
function. More precisely $\widehat{\tilde{f_{0}}}
\mid_{\mathbb{R}^{n}}=f_{0}$. Assume $\tilde{f_{0}}$ has a very
regular representation at the boundary, with isolated singularities.
Then $\widehat{f}=0$ from $P(\xi)=0$, gives an extension of
singularities to $\Gamma_{j}$, locally algebraic segments. At the
boundary, in a complex neighborhood of the real space, we are
considering the symbol as $F(\gamma)=P(D) \tilde{f_{0}}$.

\vsp

Consider in $\widehat{\DL}$, $R(\zeta)f_{0}(\zeta)$, where $f_{0}$
is the Fourier transform of a very regular operator, that is
$F(\gamma) \sim R(D) \tilde{f_{0}}$, where $\widehat{\tilde{f_{0}}}
\mid_{\R}=f_{0}$ and $\tilde{f_{0}}$ is very regular. Consider
$\Gamma \rightarrow \Gamma^{*}$ through a simple Legendre transform.
If we assume $\Sigma$ discrete and that all approximations of
$\Sigma$ through $\Gamma$ are regular (transversal intersection),
then we can assume existence of a norm $\rho$, such that $\rho(z)
\leq 1 \Leftrightarrow z \in \Gamma$. In conclusion, we are
assuming a very regular boundary continued to
$\delta_{\Gamma}-\gamma_{\delta}$, that corresponds to a normal tube
in $\Omega$, thus that all singularities are situated on first
surfaces. The parabolic singularities can be given by a very regular
boundary, that $\frac{d^{k}}{d T^{K}} F_{T}=0$ implies $T=0$ and it
corresponds to a very regular representation in the right hand side.
If the differential operator is given as hypoelliptic, then $F_{T}$
is very regular.

\subsection{A very regular representation}
For a boundary operator $L$, a very regular representation is given by
$L(f)=f + \gamma_{\delta}(f)$, where $\gamma_{\delta}$ is
regularizing, for $\delta > 0$. Note that for a finite $N$, the term
$\M_{-N} \gamma_{\delta}$ is still regularizing, for $\delta > 0$. We
note that in this representation, the locality of singularities is
not affected by the means, but the order of singularity is decreased
by the mean and increased by the mean of negative order. Thus, given
singularities of finite order, we must have that application of a
finite order mean, decreases the set of singularities. Through the
result from Iversen, we see that the set of singularities in a very
regular boundary, must be assumed of measure zero. Note however,
that if $\M_{N}(f)$ is locally injective, then the corresponding
$\M_{-N}(f)$ is locally surjective. We are
assuming the neighborhood of $\Gamma$ one-sided, that is
$\gamma_{\delta}(f) \geq 0$ for $f \geq 0$ (possibly multi-valued).
The condition on vanishing flux, $\int_{bd} (d L f)^{\vartriangle}= \int_{bd}
 d f^{\vartriangle}=0$, since $d \gamma_{\delta}=0$ for $\delta=0$. That is if
 we have a "planar" reflection through the boundary, this is
 preserved by the boundary representation. 

 \subsection{The extended Dirac distribution}
 Assume $\Sigma$ a set of common points for finitely many analytic
 curve segments, a discrete set without positive measure. The boundary
 condition corresponding to a very regular boundary can now be
 formulated, $f_{T}$ is regular outside $\Gamma_{1} \cup \ldots \cup
 \Gamma_{N}$, that is at least one of the segments $\Gamma_{j}$ is
 singular for $f_{T}$. This means that if the boundary element is
 $\delta_{\Gamma_{1}}- \gamma_{\delta}$, then $\big[
 \delta_{\Gamma_{1}},\delta_{\Gamma_{2}} \big] \neq \big[
 \delta_{\Gamma_{2}},\delta_{\Gamma_{1}} \big]$. Only points in
 $\Sigma$ give raise to a commutative system. Further, the system
 will be finitely generated in the sense that finitely many
 (sufficiently many)
 iterations of $\delta_{\Gamma_{j}}$, will for different
 $\Gamma_{j}$ produce regular points. More precisely, assume
 $\Gamma_{1}$ is a singular analytic curve, for $\phi$ and
 $\Gamma_{2}$ is a regular analytic curve for $\phi$ except for a
 point in $\Sigma$ and in $\Gamma_{1} \cap \Gamma_{2}$. Then $\phi
 \mid_{\Gamma_{2}} \mid_{\Gamma_{1}}$ is the result of a regular
 approximation of a singular point, but $\phi \mid_{\Gamma_{1}}
 \mid_{\Gamma_{2}}$ gives a singular approximation of a singular
 point. Compare Nishino's concept of a normal tube (cf.
 \cite{Nishino68}).

\vsp

Consider $\delta_{x_{0}} \rightarrow \delta_{\Gamma}$, where
$\Gamma$ has an analytic parametrization. This means $\lim_{\Gamma
\ni x \rightarrow x_{0}} \phi(x)=\phi(x_{0})$. Assume the set
$\Sigma= \cap_{j} \Gamma_{j}$ discrete (compact), that is algebraic.
If $r_{\Gamma,\Sigma}$ is the restriction homomorphism
$H_{\Gamma}(V) \rightarrow H_{\Sigma}(V)$, for a domain $V$ and if
$T \in H_{I,\Sigma}'(V)$ and we have existence of $U_{T} \in
H_{I,\Gamma}'(V)$ such that $T=r_{\Gamma,\Sigma}(U_{T})$, then we
say that $T$ has a continuation to $\Gamma$. Thus, we can see
$\delta_{\Gamma}$ as a continuation of $\delta_{\Sigma}$ to a very
regular boundary. If $V$ is a domain of holomorphy and if the
definition of a normal operator $L$ at the boundary, is not
dependent on choice of $\Gamma$, we say that $\{ \Gamma \}$ is a
quasi porteur for $\delta_{\Gamma}$. When $\Gamma$ is analytic, we
say that it is a porteur.

\section{ A monodromy condition }
\subsection{Introduction}
Since we are discussing the symbol in $(\Bm)'$, there is no obvious monodromy concept that can be assumed. Assume condition $(M_{1})$ is the condition that all means are of real type.
For $f \in (\Bm)'$ we have a local representation $f \sim P(D) \tilde{f}_{0}$, where close to the boundary $\tilde{f}_{0}$ is very regular. Assume the local condition $(M_{2})$ is the condition that $P(D)$ is hypoelliptic. Finally,

Assume $pr_{1} \mbox{:} L^{1} \rightarrow R(L_{\lambda})$, for a formally self-adjoint
and hypo-elliptic differential operator. Assume $L_{\lambda} \in \phi(\DL,L^{1})$ (unbounded Fredholm-operator)
then $\DL=X_{0} \bigotimes N(L_{\lambda})$,$L^{1}=Y_{0} \bigotimes R(L_{\lambda})$ and the
inverse $B_{\lambda}$ is bounded in $(L^{1},\DL)$. Through a fundamental theorem in Fredholm
theory, we have existence of $B_{\lambda} \in \Phi(L^{1},\DL)$ such that $B_{\lambda} L_{\lambda}=I$ on $X_{0}$
and $L_{\lambda} B_{\lambda}=I$ on $R(L_{\lambda})$.

\newtheorem{vreg}{Proposition}[section]
\begin{vreg}
 Assume $E_{\lambda}$ a parametrix to $L_{\lambda}$ a hypoelliptic operator in $L^{1}$,
then $(\delta_{x}-E_{\lambda})$ is regularizing. Conversely, if $E_{\lambda}$ is regularizing,
then $\delta_{x}-E_{\lambda}$ is hypoelliptic.
\end{vreg}

 We give a short proof, if $R$ is regularizing
then $\mid \mid \psi R u \mid \mid_{s} \leq C$, for a constant $C=C(\psi,s)$ and $\psi \in C^{\infty}_{0}$
and a real $s$ and $u \in \mathcal{D}'$. We can write $\mid \mid \psi u \mid \mid_{s} \leq
C \mid \mid \psi (A-R)u \mid \mid_{s}$, where $A=E-I$ is hypoelliptic. If $E$ is a parametrix
$PE u - u \in C^{\infty}$, for a $u \in \mathcal{D}'$ and if $P$ hypoelliptic in $L^{1}$,
$PE u-P u \in C^{\infty}$ implies $E u-u \in C^{\infty}$, that is $(\delta_{x}-E)$ is regularizing.
The result can e extended to $\mathcal{D}^{F '}$, then $I_{E}-I$ is regularizing in the space
of $\mathcal{D}^{F '}$, but the result can not be extended to $\mathcal{D}'$.

\vsp

The condition $(M_{3})$ is the condition that the symbol considered as a parametrix to a boundary operator has a trivial kernel.

\subsection{The (M)-conditions and orthogonality}
Assume $T_{\gamma}$ corresponds to analytic continuation.  We will assume that $T_{\gamma}$ can be divided into translational
movement and rotational movement, not necessarily independent. For $V=(V_{1},V_{2})$ the vorticity is given
as $w=\frac{\delta V_{2}}{\delta x}-\frac{\delta V_{1}}{\delta y}$. Given that $\frac{d y}{d x}$ is bounded in the infinity,
we have that $\frac{d V_{2}}{d V_{1}} \frac{d y}{d x}=1$ in the infinity iff $w=0$ in the infinity, that is
if the limit $\frac{V_{2}}{V_{1}}=e^{-\varphi} \rightarrow 0$, in the infinity, is a ``simple zero``. If $e^{-\varphi_{T}} \sim P(1/T)$
$\sim c_{0} + c_{1} / T+ \ldots$, then we must have $c_{0}=0$ and $c_{1} \neq 0$. If we compare
with the global problem, the condition $\ram f-f=0$ implies $T=0$ is locally a condition $(M)$ and $\ram \log f-\log f=0$
implies $T=0$ is related to local parabolicity. The condition that $\ram \log f=\log \ram f$, means that 
the ramifier only contributes to the phase. If we assume $V_{1} \bot V_{2}$, we can find a polynomial $P$
such that $\frac{V_{2}}{V_{1}} \sim_{m} \frac{1}{P}$, in the infinity.

\vsp

Note in connection with the conditions $(M)$, if we assume $V_{1} \bot V_{2}$,
then for a hypoelliptic symbol both $V_{1}$ and $V_{2}$ will be unbounded in the infinity and thus the respective inverse is bounded in $\infty$.  In presence of lineality, we will later argue that the imaginary part can be assumed bounded, where we assume $\mbox{ Im }F_{1/T} \rightarrow 1/\mbox{ Im }F_{T}$ preserves conditions $(M)$. Thus, in order to discuss the conditions $(M)$ using only bounded symbols for orthogonal parts, it is necessary to assume $F_{T}^{adj}-F_{T}$ bounded, where $F_{T}^{adj}$ is the adjoint symbol in $(B_{m})'$.

\subsection{A condition $(M_{2})$ operator at the boundary}

 We now wish to define
a formal condition $(M_{2})$ operator on $(\Bm)'$. This operator can be used to define a global condition $(M_{2})$. For such an operator, we must have that ${}^{t} T_{\gamma}$ has
representation with a point support measure $\mu_{0}$. For a condition $(M_{2})$ with respect to paths, 
the limit must be independent of starting point, why we must assume $\mu_{0} \in \mathcal{E'}^{(0)}$,
that is of order $0$. If we have the parabolic property, we must further assume the measure is positive
in phase space. 
For a formal $(M_{2})$-operator, we only assume 
$F \in L^{1}$ implies ${}^{t} T_{\gamma} F \in (\Bm)'$. Consider now ${}^{t} T_{\gamma}$ modulo regularizing action, 
that is ${}^{t} T_{\gamma} \sim \delta_{0} + \nu_{\gamma}$, where $\nu_{\gamma}$ has support outside a point.
If we assume ${}^{t} T_{\gamma} F \in \Bm$, we have $\widehat{{}^{t} T_{\gamma} F}=Qf_{0}$. We write
${}^{t} T_{\gamma} F = Q(D) \tilde{f_{0}}$, where $\tilde{f_{0}}$ is a very regular operator. 
Assume $F \in \Bm$, then we have $F \sim Q_{\lambda} \tilde{f_{0}}$. We can in this context consider $\tilde{f_{0}}$
as a global representation. Given that $Q_{\lambda}$ is hypoelliptic and such that $Q_{\lambda} {}^{t} T_{\gamma}={}^{t} T_{\gamma}Q_{\lambda}$
we see that in the case with condition $(M_{2})$, ${}^{t} T_{\gamma}$ preserves parametrices to $Q_{\lambda}$.
Conversely, given that ${}^{t} T_{\gamma}$ preserves very regular parametrices to $Q_{\lambda}$ (that is we have condition $(M_{2})$)
then we can derive that $Q_{\lambda}$ is hypoelliptic. The conclusion is that with the conditions above,
we have that hypoellipticity for $Q_{\lambda}$ means that ${}^{t} T_{\gamma}$ defines condition $(M_{2})$
and conversely. Note also that ${}^{t} T_{\gamma}$ can be globally represented in $(\Bm)'$. Given $Q$,
we can define the possible continuations that preserve condition $(M_{2})$, as ${}^{t} T_{\gamma}$ such that ${}^{t} T_{\gamma} Q=Q {}^{t} T_{\gamma}$.
For instance if ${}^{t} T_{\gamma} F=F + F_{1}$, we must assume $QF_{1}=0$, why $F_{1}$ is regularizing.
If on the other hand ${}^{t} T_{\gamma} F=c F$, for a $c$ in the ${\infty}$ and ${}^{t} Q=Q^{*}$,
then for all real $c$, we have $c Q - Qc \prec \prec Q$, that is ${}^{t} T_{\gamma} Q - Q {}^{t} T_{\gamma}=0$
modulo regularizing action. Thus, localizing with a real $c$ is possible and corresponds to
localizing $\tilde{f_{0}}$. Given that $F$ is hypoelliptic and represented as $Q \tilde{f_{0}}$
at the boundary, this property can be continued using ${}^{t} T_{\gamma}$, given that this is an
algebraic action. Assume for this reason, $({}^{t} T_{\gamma} F)(\varphi)=\int_{\Omega} k_{\gamma} \varphi d \sigma=0$
implies $Q \varphi=0$ on $\Omega$ or $\sigma( \Omega)=0$, that is in the case where ${}^{t} T_{\gamma}$ is algebraic, we have condition $(M_{2})$.
Note that $\int_{\Omega} \tilde{f_{0}} \varphi d \sigma=0$ implies $\sigma (\Omega)=0$, which corresponds to ''isolated zeros`` to $\tilde{f_{0}}$.
In the case where ${}^{t} T_{\gamma}(\tilde{f}_{0})(\varphi)=0$, when ${}^{t} T_{\gamma} \tilde{f}_{0}$ is
algebraic, we see that $\sigma(\Omega)=0$, that is we have condition $(M_{2})$.Note that if ${}^{t} T_{\gamma} \tilde{f_{0}}=0$
along a line in $\Omega$, this can be compared with the extension $\delta_{0}$ to $\delta_{\Gamma}$,
where condition $(M_{2})$ is not longer with respect to a point.

\newtheorem{M2}[vreg]{ Lemma }
\begin{M2}
Assuming the symbol has a representation in $(\Bm)'$ satisfying condition $(M_{2})$, then this property is preserved if the continuation ${}^{t} T_{\gamma}$ as above is algebraic. 
\end{M2}
Condition $(M_{2})$ does not however imply that ${}^{t} T_{\gamma}$ is algebraic.
\subsection{The operator $\Top$ and the conditions (M)}
We have that $(w_{\mu})^{\diamondsuit}=(w^{\diamondsuit})_{\mu}$ but $(\Top w)_{\mu} \neq \Top(w_{\mu})$.
Assume $w=W/M \sim e^{\varphi} \rightarrow \Top e^{\varphi}=e^{\varphi^{\vartriangle}}$.The problem is when $e^{\varphi^{\vartriangle}}$
respects the conditions (M). Assume $H(w)=w^{\diamondsuit}$ and $\Top$ Legendre and that we have algebraic dependence of the parameter $T$, then $\Top H \tau_{\mu}-\tau_{\mu} \Top H \sim_{0} e^{P(\frac{1}{T})}-e^{P(T)}$, for a polynomial $P$. The 
condition $(M_{2})$ means that $e^{P(\frac{1}{T})-P(T)}$ is not $\sim 1$.
Thus, the operator $\Top$ does not necessarily preserve the condition $(M_{2})$. The same holds for the conditions $(M_{1})$,$(M_{3})$.

\vsp

 If the operator $\Top$ is studied using ${}^{t} T_{\gamma}$ acting on the Legendre transform,
we see that an algebraic continuation in the infinity implies that we do not have lineality
and further the conditions $(M_{1})$ and $(M_{2})$ in parameter infinity. If we only assume
the continuation very regular on all strata, we still have conditions $(M_{1}),(M_{2})$
and $(M_{3})$, but not necessarily an algebraic continuation. We can consider it to be
algebraic modulo monotropy locally. Finally note that an algebraic continuation does not necessarily preserve condition $(M_{1})$. However, by considering symbols modulo regularizing action, we can restrict the representation to real type symbols. For this representation the corresponding functional is an infinite sequence of constant coefficients polynomials acting on a point support measure.
 
\section{Further remarks on algebraicity}

\subsection{Introduction}
An algebraic set is geometrically equivalent to a zero set of polynomials. Characteristic for an algebraic mapping $L$ is that $L(e^{\varphi})=e^{L(\varphi)}$ and the zero set to an algebraic mapping is locally an algebraic set. The identity (evaluation) is considered as algebraic and we consider any operator that commutes with the identity in $H'$ as (topologically) algebraic.

\subsection{Clustersets for multipliers}
We will prove that given that $M \bot W$ implies $V_{1} \bot V_{2}$, we have that $\Top V^{-N}$ algebraic in $\infty$ if we have hypoellipticity.
Assume for this reason $V_{1} \bot V_{2}$ with $\frac{\widehat{V_{2}}}{\widehat{V_{1}}}=e^{-\vartheta^{*}} \rightarrow 0$, as $\mid \xi \mid \rightarrow \infty$
If $\Top V_{1}=e^{\vartheta^{\vartriangle}} \Top V_{2}$, then we must have $e^{\vartheta^{\vartriangle}} \rightarrow \infty$, as $\mid \xi \mid \rightarrow \infty$.
For instance if $\vartheta^{\vartriangle}=\vartheta^{*}+p$, where we assume $\frac{1}{p}$ is harmonic, in the sense that $p \rightarrow \infty$, as $\mid \xi \mid \rightarrow \infty$.
Otherwise, if we assume $e^{\vartheta^{\vartriangle}}$ bounded as $\mid \xi \mid \rightarrow \infty$, we have unbounded sublevel sets (cluster sets) for the multipliers. When $\Top$ is the Legendre transform, this will disrupt the concept of monodromy in the infinity.
We have noted that if $V_{T}$ is algebraic in $T$ and $Flux(V)=0$, we have that $\mbox{ Im }V_{T} \equiv 0$. Note that if $\Top V^{-N}$ is taken as
the limit over strata, in case the symbol is not hypoelliptic, we have distributional limits in symbol space for the representation of $V$.

\subsection{Example on algebraic mappings}
Consider the mapping $A \big[ a,b \big]=\widehat{b}a$. If $E$ is a (topologically)
algebraic mapping, and if $E(e^{\phi})=e^{\tilde{E}(\phi)}$, then $\tilde{E}$ is odd and also $\tilde{A}$. The
proposition that $\tilde{A}^{2}$ is odd, corresponds to
$E(e^{\frac{1}{\widehat{I}(\phi)}})=\frac{1}{E(e^{\widehat{I}(\phi)})}$,
thus if $A$ is (topologically) algebraic, then $A^{2}$ is (topologically) algebraic. Note that
$E(e^{e^{-\phi}})=\frac{1}{E(e^{e^{\phi}})}$. If $E \sim_{m} E_{1}$
and $E_{1}$ is (topologically) algebraic, then $\mid E(e^{f})-E_{1}(e^{1/f}) \mid <
\epsilon$ in $\infty$. If $E \in L^{1}$, then $\widehat{E}(f)
\rightarrow 0$, as $f \rightarrow \infty$. This means that $E
\sim_{m} E_{1}$, where $E_{1}$ is (topologically) algebraic, then $\mid
\widehat{E}(f)-\widehat{E}_{1}(\frac{1}{f}) \mid \rightarrow 0$. We
can start with a condition $\M_{N}(E) \in L^{1}$, which means that
$\M_{N}(E) \sim_{m} W$, where $\tilde{W}$ is odd.

\vsp

Further, we have
$A^{3}\big[I,E\big](\phi)=E(e^{\phi})$ and in the same manner
$A^{3}\big[E,I \big](\phi)=$ $e^{\widehat{\widehat{E}}(\phi)}$.
Thus, if $E(e^{\phi})=$ $e^{e(\phi)}$, that is $v(\phi)=\big[
\widehat{I},e \big](\phi)$, then $\widehat{\widehat{E}} \sim e$. If
$E$ is assumed (topologically) algebraic and $\widehat{\widehat{E}} \sim E$, then we
should have that $e$ is (topologically) algebraic (modulo algebraic sets). Further,
if $e$ is linear and $\widehat{\widehat{E}} \rightarrow I$ at the
boundary and if $\widehat{\widehat{E}} \sim E$, then we must have
$E-I \sim W$, where $W$ is algebraic at the boundary (with respect to
concatenation of curve segments). More precisely
$(E-I)(E+I)=E^{2}-I+R$, where $R$ is assumed to vanish. If $E-I
\rightarrow I$ at the boundary and the same condition holds for
$E+I$, we must have $E^{2}-I \rightarrow 0$, which is seen as $E$
being (topologically) algebraic at the boundary, that is $E-I \sim W(e^{-\phi})$,
where $W$ is algebraic in $T$ (compare the
Lagrange condition). We are assuming $W(e^{\phi})W(e^{-\phi})=W(1)$
Thus, if $W(1)=1$, we have $W(e^{-\phi})=e^{\tilde{W}^{-1}(\phi)}$
and the odd condition means that
$\tilde{W}^{-1}(\phi)=-\tilde{W}(\phi)$.

\subsection{The Legendre transform has removable singularities}

 Assume $\vartheta_{T}^{\vartriangle}/\vartheta_{T} \rightarrow
0$, as $T \rightarrow \infty$ and further that
$<\vartheta^{\vartriangle} \pm \vartheta,\vartheta^{\vartriangle}
\pm \vartheta> =0$ or formally $\vartheta^{\vartriangle}/\vartheta +
\vartheta/\vartheta^{\vartriangle} \sim 1$. Assume
$F_{T}^{\vartriangle}.F_{T} \sim 1$, as $T \rightarrow \infty$ and
the first condition, then the condition $ v/v^{\vartriangle} \sim 1$
is impossible. If $\vartheta^{\vartriangle}-\vartheta$ real, then
the condition contradicts symplecticity. Otherwise, we are assuming
$\mbox{ Im }(\vartheta^{\vartriangle}-\vartheta) \bot \mbox{ Re
}(\vartheta^{\vartriangle}- \vartheta)$. The conclusion is that if
$\vartheta^{\vartriangle}$ and $\vartheta$ are related through a
simple Legendre transform, then $F^{\vartriangle}_{T}.F \sim_{T} 1$,
as $T \rightarrow \infty$ is not possible, that is the singularities
are removable. 

\subsection{Sufficient conditions for orthogonality}

Assume now $M \bot W$ and consider the lift $V(M,W)=V_{1}+iV_{2}$.
The problem is now under what conditions we have that $V_{1} \bot
V_{2}$ (that is $V_{2} \prec \prec V_{1}$). Assume $V$ is the
Hamilton function, that is $M= \frac{\delta V_{2}}{\delta
y}=\frac{\delta V_{1}}{\delta x}$ and $W=\frac{\delta V_{1}}{\delta
y}=-\frac{\delta V_{2}}{\delta x}$. We write formally the condition
that $M \bot W$ as $\frac{\delta}{\delta x}(V_{1} \bot V_{2})$ and
$\frac{\delta}{\delta y}(V_{1} \bot V_{2})$. The vanishing flux
condition is $Flux(\mathcal{T} W)=0$. The corresponding condition
for $\mathcal{T} V_{2}$ is $V_{2} \prec \prec V_{1}$. Note that the
lifting function $F(X,Y)$ and the Hamilton function $V(M,W)$ may
have quite different algebraic properties, but are considered as
related by involution and the condition on equal derivatives in the
first order $x^{*},y^{*}$-variables, with respect to $(X,Y)-$
arguments. For $V_{1}$ algebraic, we always have that
$\frac{\delta}{\delta x}V_{1} \bot V_{1}$, that is
$\frac{\delta}{\delta x} \log V_{1} \rightarrow 0$, as $\mid T \mid
\rightarrow \infty$. The condition $ \frac{\delta}{\delta x} V_{1}
\bot V_{2}$ is $(\frac{\delta}{\delta x} \log V_{1})
\frac{V_{1}}{V_{2}} \rightarrow 0$, as $\mid T \mid \rightarrow
\infty$ and if $\frac{\delta}{\delta x} \log V_{2}$ is bounded in
the infinity, then we could write the condition
$\frac{\delta}{\delta x} V_{2} \bot V_{1}$ as $(\frac{\delta}{\delta
x} \log V_{2}) \frac{V_{2}}{V_{1}} \rightarrow 0$, as $\mid T \mid
\rightarrow \infty$.

\vsp

Let $W/M=e^{\varphi}$ and consider the iterated mean
$\M_{N}(\varphi)$, for $N$ large. We can then assume $\frac{d}{d x}
\M_{N}(\varphi)$ reduced at the boundary. Let
$V_{1}^{(N)},V_{2}^{(N)}$ be the corresponding Hamilton function
(locally). Consider the condition $(\frac{\delta}{\delta x} \log
V_{1}^{(N)}) \frac{V_{1}^{(N)}}{V_{2}^{(N)}} \rightarrow 0$, as $T
\rightarrow 0$. This is obviously true for large $N$. If
$\frac{\delta}{\delta x} \log V_{1}^{(N)}$ is reduced at the
boundary, then we have that $V_{2}^{(N)}/V_{1}^{(N)} \rightarrow 0$,
as $T \rightarrow \infty$. Further, $\frac{\delta}{\delta
x}V_{1}^{(N)} \bot V_{2}^{(N)}$. A sufficient condition to conclude
that $V_{1} \bot V_{2}$ is that $V$ corresponds to a reduced symbol.
Note that we have assumed that $V_{1} \mid_{\Gamma} \sim
\frac{1}{V_{1}} \mid_{\infty}$ and $\frac{1}{V_{2}} \mid_{\Gamma}
\sim V_{2} \mid_{\infty}$.

\subsection{ The orthogonal condition and the Dirichlet problem }
Assume instead of the condition $\mbox{Im} F_{T} \prec \prec \mbox{Re} F_{T}$,
that
\begin{equation}\label{Tauber}
    \mbox{Re} F_{\frac{1}{T}} \mbox{Im} F_{T} \rightarrow 0 \mbox{  as  } T
    \rightarrow \infty
\end{equation}

The condition is to be understood using $T=(\mbox{Re} T,\mbox{Im} T) \in \R$,
$n=2$ and $\mid \mbox{Im} F_{T} \mid \mid \mbox{Re} F_{\frac{1}{T}} \mid
\rightarrow 0$, as $\mid T \mid \rightarrow \infty$, where we assume
that the factors are not without support in a neighborhood of the
infinity. Given that $F_{T}$ is holomorphic in $T$, we know that
$F_{T}$ can not be reduced in the origin and in the infinity
simultaneously. Assume $\mbox{Re} F_{\frac{1}{T}}$ reduced in the $T-$
infinity, we then have that $\mbox{Im} F_{T}$ is bounded in the infinity.
Assume $F_{T}^{\vartriangle}={}^{t}F_{T}$, so that $\mbox{Im} F_{T} \sim
F^{\vartriangle}_{T}-\overline{F}_{T}$. The condition that $F_{T}$ is algebraic
in $T$ close to the boundary, does not imply that $\mbox{Im} F_{T}$ is
algebraic in $T$, close to the boundary. As we are assuming real
dominance for the operator, we will assume the real part
algebraically dependent on the parameter.

\vsp

In the same manner, given that $F_{T} \in \DL$, we do not
necessarily have that $F_{T}^{\vartriangle}-\overline{F}_{T} \in \DL$. Note
that for the restriction to the real space, the Fourier transform to
$F_{T}$ is on the form $P(\xi)f_{0}(\xi)$, why on the support of
$f_{0}$, we necessarily have finite order of zero. We are assuming
$0 \in \mbox{ supp }f_{0}$. We have that $\mbox{Re} F_{\frac{1}{T}}$ can
not be reduced in a neighborhood of parameter origin. If $\mbox{Im} F_{T}$
were bounded in a neighborhood of the origin, the condition
(\ref{Tauber}) could not be possible. The conclusion is that $\mbox{Im}
F_{T}$ must be unbounded at the boundary.

Conversely, if we assume $\mbox{Re} F_{\frac{1}{T}}$ is reduced in a
neighborhood of the parameter origin, then we have that $\mbox{Im} F_{T}$
is bounded in a neighborhood of the boundary and we have that $\mbox{Im}
F_{T}$ is necessarily unbounded in the infinity. This, does not
contradict the condition (\ref{Tauber}).

\newtheorem{bounded}{ Proposition }[section]
\begin{bounded}
 Assume condition (\ref{Tauber}). For a hypoelliptic operator, both the real part and the imaginary part are unbounded in the infinity. In presence of lineality for the real part
the imaginary part is bounded in the infinity.
\end{bounded}

The conclusion is that given the condition (\ref{Tauber}), if we
further assume that $\mbox{Re} F_{T}$ is reduced in the infinity, then we
have that $\mbox{Im} F_{T}$ can be assumed to be bounded at the boundary.
The proposition that $\mbox{Im} F_{T}$ is bounded in the infinity, means
that it can not be represented be a polynomial operator, unless the
set of unboundedness is of measure zero. In this case if the set is
normal (finite Dirichlet integral), we can give an algebraic
representation for this part of the symbol. Sato gives a well known
example of a (hyper-) function defined at the boundary, that is not
a distribution in $0$. Assume $\vartheta_{T}$ algebraic in $T$ and
that $\vartheta_{T}^{\vartriangle}-\overline{\vartheta}_{T} \sim \mbox{Im}
\vartheta_{T}$, such that $\mbox{Im} \vartheta_{T} \sim
P(\frac{1}{T})\vartheta^{\vartriangle}-\overline{P(T)}\vartheta$.

\subsection{ A necessary condition }

Assume
\begin{equation}\label{TauberII}
    \frac{e^{\alpha(\frac{1}{T})}-e^{\overline{\alpha(T)}}}{e^{\alpha(T)}} \rightarrow 0 \mbox{ as } T \rightarrow \infty
\end{equation}

We then have, $e^{\alpha(\frac{1}{T})-\alpha(T)}-e^{-2 i \mbox{Im}
\alpha(T)} \rightarrow 0$, as $T \rightarrow \infty$. A sufficient
condition for this is that $\alpha(\frac{1}{T})-\alpha(T)- 2 i \mbox{Im}
\alpha(T) \rightarrow 0$, where $\alpha(\frac{1}{T}) \sim \overline{\alpha(T)}$.
Thus, the relationship between $\alpha^{\vartriangle}(T)=\alpha(\frac{1}{T})$
and $\alpha(\frac{1}{T})=\overline{\alpha(T)}$, as a simple Legendre
type condition, means that if $e^{\alpha}$ is locally algebraic,
then the same must hold for $e^{\alpha(\frac{1}{T})}$. Formally
$$ \Top(\mbox{Im} e^{\vartheta_{T}}) \sim e^{\overline{\alpha(T)} \mbox{Im}
\vartheta}$$ where $e^{\overline{\vartheta}_{T}}$ is locally algebraic.
Note that without the condition (\ref{TauberII}), it is necessary to
consider $\Top$ as acting on hyper-functions.

\subsection{The complement sets to the first surfaces}
Solvability in this context corresponds to the regularity conditions
in dynamical systems and the corresponding conditions for first
surfaces. The neighborhood $\{ \varphi > 0 \} \sim \{ \vartheta = 0
\}$, that is $\sum \varphi_{j} > 0$ or equivalently $\sum
\vartheta_{j} \equiv 0$, gives an analytic parametrization. If the
set is not analytic, we do not have local solvability. (The Fuchs
condition). If we do not have regular approximations of a first surface,
the representation of the operator is not defined.
Thus a local complement $\{ f=c \}^{c}$ analytic, is necessary
for solvability. A locally algebraic transversal is necessary for
hypoellipticity. For example, if $e^{-\phi} \widehat{f}_{T}
\rightarrow 0$, as $T \rightarrow \infty$, given an essential
singularity in the infinity, we may have local solutions that are
not analytic.

\vsp

Assume $\{ f=c \}$ a first surface to a holomorphic function. Given
minimally defined singularities to $f$, we know that the first
surfaces are also analytic. This means that if $f \in L^{1}$, the
complement in a domain of holomorphy is analytic. If we only assume
$f \in \DL$, we do not have this strong result.
Given that $\Delta$ is analytic and $\Delta \rightarrow
\Delta^{\bot}$ is a Legendre transform, that is a contact transform,
we should be able to prove $\Omega \backslash \Delta \sim \Omega
\backslash \Delta^{\bot}$, that is if $\Delta^{\bot}$ is a first
surface, it must be analytic. If $f$ is reduced, we have that
$\frac{d}{dT} \log f \rightarrow 0$ implies $\frac{d}{dT}f
\rightarrow 0$. Otherwise, if $1/f$ is bounded, then the converse
implication holds. Algebraic dependence implies that $\frac{d}{dT}
\log f \rightarrow 0$. Since $\{ f=c \} \subset \{ \frac{d}{dT}f =0
\}$, we must have $\Omega \backslash \{ f=c \}$ analytic implies $\{
\frac{d}{dT} f \neq 0 \}$ analytic. Obviously $\frac{d}{dT} \log f
\neq 0$ implies $\frac{d}{dT} f \neq 0$ (assuming both not zero). Thus, $\{ \log f \neq 0 \}$
is analytic and $\{ f \neq c \}$ is analytic. In this case, if $\{
f=c \}$ is analytic, then $\{ \log f=c \}$ is analytic. If the
complement to the second set is analytic, the complement to the
first set is analytic. The last condition implies solvability.
Counterexamples can be given by $f=\beta e^{\varphi}$.

 \subsection{A counterexample to solvability}
 The problem that we start with is when the complement to a fist
 surface in a domain of holomorphy, is analytic? Assuming parabolicity,
 we make the approach $e^{\phi}$, with $\phi$ subharmonic.
 Assume the neighborhood of the first surface is given by $\{ \phi
 \geq 0 \}$. We are now discussing $\phi(x-y)=-\phi(y)+\psi(x,y)$.
 Through Tarski-Seidenberg's result, we have if $x,y$ in a
 semi-algebraic $\{ P(x,y) \geq 0 \}$ and related through a Legendre transform, we have
 $y$ in a semi-algebraic set $\{ Q(y) \geq 0 \}$. The problem is
 now,  if $\{ Q(y) \geq 0 \}=\{ F_{Q}(y)=0 \}$ locally, where
 $F_{Q}$ is analytic. For instance $R$ analytic with
 $R(x,y)F_{P}(x-y)=F_{Q}(y)$. Assume now $x=0$, so
 $\phi(-y)=C-\phi(y)$. If $\phi(-y) \leq 0$, then $\phi(0)-\phi(y)
 \leq 0$ and in the parabolic case $\phi(-y)=C_{1}$, a constant, as
 $y \rightarrow \infty$. Thus, $\phi(y)=C-C_{1}$ constant. If
 $\phi(-y) > 0$, then $C-\phi(y)>0$, that is $\phi(y)$ is bounded,
 as $y \rightarrow \infty$. In this case $\{ y \quad \phi(y) <C \}$
 is unbounded. The integral $\int_{\phi < C}$ now corresponds to a
 distribution. In this case $\{ \phi > 0 \}$ is not a semi-algebraic
 set. If $\phi$ is hypoelliptic, then $\phi \rightarrow \infty$, as
 $y \rightarrow \infty$. If $\phi$ is bounded, as $y \rightarrow
 \infty$, then $\phi$ must be not-hypoelliptic. In this case it is
 not possible to approximate $\infty$ through $\{ \phi > 0 \}$. The
 proposition in this case $\int_{\{ \phi > 0 \}}$ can not be
 represented by an analytic function. Thus, if $\phi$ has a
 parabolic and odd representation and if we do not have an essential
 singularity in $\infty$, then there is a $\varphi$ analytic and
 zero on $\{ \phi > 0 \}$. If however, $\phi(-y)=-\phi(y)+C$, we do
 not longer have representation with an analytic function.

 \vsp

 We have that given an ideal with Schwartz type topology, finitely
 generated and symmetric over $\Omega$, where $\Omega$ is
 pseudo-convex, if the dependence of the parameter
 is in $L^{1}$ in the symbol space and algebraic in its tangent space, we can assume
 that we have (topologically) isolated singularities. In the case where
 $\phi(-y)=-\phi(y)+C$, the approximation of the singularity only
 exists in a weak sense.

\section{Not hypoelliptic operators}

\subsection{Introduction}
To determine the class of hypoelliptic pseudo differential operators, we will first assume the operator has representation as an unbounded Fredholm operator with symbol in the radical to the ideal of hypoelliptic operators. In the last section we consider a different representation based on Cousin (cf. \cite{Cousin95}).
 
\subsection{ Pseudodifferential operators as Fredholm operators }
From the theory of linear Fredholm-operators, we know that any Fredholm operator $A: E_1 \rightarrow E_2$
between Banach spaces, has a twosided Fredholm inverse, that is there is a $B: E_2 \rightarrow
E_1$, such that $BA=I-P_1$, $AB=I-P_2$, with $P_1,P_2$ finite rank operators, $P_1$ is the
projection $E_1 \rightarrow \mbox{ Ker }A$ and $(I-P_2)$ is the projection $E_2 \rightarrow \mbox{
Im }A$. Conversely, given an operator $A$, continuous and linear $E_1 \rightarrow E_2$, such that we have
operators $B_1,B_2$ with $B_1A=I+R_1$, $AB_2=I+R_2$, $R_1,R_2$ compact operators, then $A$ is
Fredholm. Finally, the class of Fredholm operators is invariant to addition of a compact operator.

\vspace*{.5cm}

For our pseudodifferential operator $A$, given existence of left- and right parametrices, the operator $A$ is Fredholm and
we have a Fredholm inverse. We have earlier noted that our pseudodifferential operator $A$ can be compared with
polynomial operators according to $H_{\lambda}=A_{\lambda}-P_{\lambda}$, with $H_{\lambda}$
regularizing. Any parametrix to $P_{\lambda}$ can be considered as a parametrix to $A_{\lambda}$.
Assume now $B_{\lambda}$
the Fredholm inverse to $P_{\lambda}$, modified as in the preliminaries. We then know that $B_{\lambda}-I$ is
regularizing outside $\mbox{ Ker }B_{\lambda}$ as $\mid \lambda \mid \rightarrow \infty$. We shall see below, that
this is a "radical" property, which means that for $\varphi \notin \mbox{ Ker }B_{\lambda}$ we have
$\mbox{ sing supp }_{L^2}(\varphi)=$ $\mbox{ sing supp }_{L^2}(B_{\lambda}P_{\lambda}\varphi)=$ $\mbox{ sing supp
}_{L^2}(P_{\lambda}\varphi)$. Naturally, $A_{\lambda}$ has the same domain for hypoellipticity as
$P_{\lambda}$. Assume $E_{\lambda}$ a left- and $F_{\lambda}$ a right-parametrix to $A_{\lambda}$
and $B_{\lambda}$ the modified Fredholm parametrix to $P_{\lambda}$. Then
$A_{\lambda}\big( E_{\lambda}-B_{\lambda} \big)=R_1 + P_2$, $R_1$ regularizing and $P_2: E_2
\rightarrow \mbox{ Im }P_{\lambda} {}^{\bot}$ continued with regularizing action outside $\mbox{
Ker }B_{\lambda }$. As $A_{\lambda}$ is hypoelliptic outside $\mbox{ Ker }B_{\lambda}$,
$E_{\lambda}=B_{\lambda} + \Gamma_1$, with $\Gamma_1$ regularizing. In the same manner
$F_{\lambda}=B_{\lambda} + \Gamma_2$, with $\Gamma_2$ regularizing. The construction gives that $\mbox{ Ker
}B_{\lambda}$ is a finite dimensional space and on this space any parametrix to $A_{\lambda}$ is
either regularizing or $0$. We can assume $\Gamma_j=0$ on $\mbox{ Ker }B_{\lambda}$.

\newtheorem{b_he-psdo}{ Proposition}[section]
\begin{b_he-psdo}
Assume $B_{\lambda}$ the modified Fredholm inverse to $P_{\lambda}$ as above, that is
$B_{\lambda}P_{\lambda}=I-P_1$, where $P_1$ is regularizing and $P_{\lambda}B_{\lambda}=I-P_2$
with $P_2$ regularizing outside $\mbox{ Ker }B_{\lambda}$. Further $A_{\lambda}$ a
pseudodifferential operator so that $H_{\lambda}=A_{\lambda}-P_{\lambda}$ with $H_{\lambda}$
regularizing, then $A_{\lambda}$ is
hypoelliptic in $L^2$ if and only if $\mbox{ Ker }B_{\lambda}=\{ 0 \}$
\end{b_he-psdo}

The proposition can be read as follows, for a hypoelliptic pseudodifferential operator in our class $A_{\lambda}$, we
have that for $\varphi \in L^2$, $B_{\lambda}\varphi=0$ implies $\varphi = 0$. The following Lemma
is trivial

\newtheorem{b_psdo-he}[b_he-psdo]{ Lemma }
\begin{b_psdo-he} \label{b_Lemma1}
If $P_{\lambda}^N$ is hypoelliptic, then $B_{\lambda}^N - B_{N \lambda} \in C^{\infty}$.
\end{b_psdo-he}

By choosing $\lambda$ appropriately, we can assume $N=2$.

\newtheorem{b_psdo-he2}[b_he-psdo]{ Lemma }
\begin{b_psdo-he2} \label{b_Lemma2}
For $u \in L^2$
$\mbox{ sing supp }_{L^2} \big( B_{\lambda}u - u \big)=\mbox{ sing supp }_{L^2}\big( B_{\lambda}u + u \big)$
\end{b_psdo-he2}

Proof:
Assume $N=2$. We obviously have $P_{2 \lambda} - I \sim P_{2 \lambda} + I$, which means that the
singular supports for $P_{2 \lambda}(B_{2 \lambda} + I)u$ and $P_{2 \lambda}(I - B_{2 \lambda})u$
coincide. Thus, the lemma holds for $B_{2 \lambda}$. Finally, the singular support for $(I + B_{2
\lambda} + B_{2 \lambda}B_{\lambda})u$ is the same as the one for $(I + B_{2 \lambda}B_{\lambda} -
B_{2 \lambda})u$. $\Box$

\vspace*{.5cm}

According to Lemma \ref{b_Lemma1} $(B_{\lambda} - I)(B_{\lambda} + I) + I-B_{2 \lambda}$ is
regularizing and according to Lemma \ref{b_Lemma2}, this implies $(B_{\lambda} - I)^2 + I -B_{2
\lambda}$ is regularizing. So, if $I - B_{2 \lambda}$ is regularizing we have $I - B_{\lambda}$ is
regularizing outside $\mbox{ Ker }B_{\lambda}$

\newtheorem{b_psdo-he3}[b_he-psdo]{ Proposition }
\begin{b_psdo-he3}
If $P_{\lambda}^N$ is hypoelliptic, then $I - \sum^{N-1}_{j=1} B_{\lambda}^j$ is regularizing outside
$\mbox{ Ker }B_{\lambda}$
\end{b_psdo-he3}

Proof:
We have $(B_{\lambda} - I)\sum^{\infty}_{j=0} B_{\lambda}^j = I$. Thus, $( B_{N \lambda} -
I)\sum^{\infty}_{j=0} B_{\lambda}^j+$ $+ \sum_{j=1}^{N-1} B_{\lambda}^j \sim I$, where $(B_{N \lambda} - I)$
is regularizing outside $\mbox{ Ker }B_{\lambda}$. $\Box$

\vspace*{.5cm}

Assume one more time $N=2$, then $N(B_{2 \lambda}) \subset N(B_{\lambda})$ and if $\varphi \in
N(B_{\lambda})$, according to Lemma \ref{b_Lemma1}, we have $B_{2 \lambda} \varphi \in C^{\infty}$.
This gives good estimates using Hadamard's lemma as we have, $\parallel B_{N \lambda}
\parallel_{c} \leq C_N \mid \lambda \mid^{-N} \  N \geq 2$. Finally, since $B_{2 \lambda}$ is a
$L^2$-kernel of finite rank, we can find a canonical kernel $K_{2 \lambda}=\sum^p_1 k_{\mu \nu}^2
e_{\mu} \otimes e_{\nu}$ such that $\parallel B_{2 \lambda} - K_{2 \lambda} \parallel=0$ and
$e_{\mu}$ some orthonormal system. That is $K_{2 \lambda}\varphi=\sum k_{\mu \nu}
(\varphi,e_{\nu})e_{\mu}$. If $\varphi \in L^2$ is such that $B_{\lambda}\varphi=0$, using the theory of
integral equations (cf. \cite{b_Smithies}) we can make an
orthogonal decomposition of $N(B_{\lambda})$ according to
$$ S_{N \lambda}=R(B_{N \lambda}) \cap R( \overline{B}_{(N+1)\lambda} )^{\bot}$$
Thus, $\varphi=B_{N\lambda}\varphi + N \lambda B_{(N+1)\lambda} \varphi$ on $S_{N\lambda}$. If
$N>2$, say $N_0$, then $S_{M\lambda} \subset C^{\infty}$ for $M \geq N_0$. On $S_{N \lambda}$ and
$N < N_0$ we have $\varphi - B_{N \lambda}\varphi \in C^{\infty}$ with $K_{N
\lambda}\varphi=\sum_{\mu,\nu} k_{\mu \nu}^N(\varphi,e_{\nu})e_{\mu}$. These kernels $B_{N\lambda}$ can
thus be considered as hypoelliptic.

\subsection{Asymptotically hypoelliptic operators}

If $B_{\lambda}$ is a parametrix to an operator $L_{\lambda}$ in $L^{1}$ and $B_{\lambda}$
$L^{1}-$ hypoelliptic on $N(B_{\lambda}^{N_{0}})^{\bot}$, where $N_{0}$ is the minimal integer
such that the zero-space is stable. Then, if $N_{0}=1$, we have that $L_{\lambda}$ is hypoelliptic
on $L^{1}$. Assume $B_{\lambda,N}$ a $L^{1}-$ parametrix in $L^{1}$ to the iterated operator
$L_{\lambda}^{N}$, then $N(B_{\lambda,N_{0}}^{N_{1}})=N(B_{\lambda}^{B_{0}})$ implies $N_{1}=1$.
and $B_{\lambda,N_{0}} L_{\lambda}^{N_{0}}=L_{\lambda}^{N_{0}} B_{\lambda,N_{0}}=\delta_{x}-\gamma$
for some $\gamma \in C^{\infty}$ where $\gamma=0$ on $N(E_{\lambda}^{N_{0}})$. Further, we have
that $L_{\lambda}^{N_{0}}$ is hypoelliptic on $L^{1}$. If $B_{\lambda,N_{0}}$ is parametrix
to a $L^{1}-$ operator $L_{\lambda}^{N_{0}}$, by adding a solution to the homogeneous equation
$H \neq 0$ on $N(B_{\lambda,N_{0}})$, we get $N(E_{\lambda,N_{0}}+H)=\{ 0 \}$ and
$B_{\lambda,N_{0}}-B_{\lambda}^{N_{0}} \in C^{\infty}$,
that is the parametrix to $L_{\lambda}^{N_{0}}$ is regularizing on $N(E_{\lambda}^{N_{0}})$.

\subsection{Iversen's condition and hypoellipticity}
\label{Iversen}
We have seen that when $\mathcal{T}=\mathcal{L}$ (Legendre), we do
not have micro local contribution. Consider for this reason the
condition: If the phase space sequence $v_{\mu}$ in

\begin{equation} \label{HE}
 0 < \mid \mathcal{T}
v_{\mu} - \mathcal{L} v_{\mu} \mid < \epsilon \qquad \mu \rightarrow
0
\end{equation}
gives a regular approximation $e^{v_{\mu}}$ of a singular point that can be continued
analytically over the origin, then the corresponding operator does
not have micro local contribution from this sequence. More
precisely, if $V(M,W)=V_{1}+iV_{2}$ is the lifting function, let
$V^{N}$ be the localization to $X_{N}$. For
$\widehat{V^{N}}=P_{N}(T) f_{0}$, for $f_{0}$ corresponding to a
very regular operator, we assume $\widehat{V}^{-N}=(1/P_{N})g_{0}$,
where $g_{0}$ is in the same class as $f_{0}$. We then have that for
$V$ corresponding to a hypoelliptic operator, we must have that
the continuation from $\widehat{V}^{-N}$,
$\mathcal{(T-L)} V^{-N}$ is algebraic as $T \rightarrow \infty$.
Conversely, given that $\mathcal{(T-L)}V^{-N}$ is algebraic, there is no
room for lineality. In the terminology of Iversen (cf. \cite{Iversen14}), if we can find a
Jordan arc emanating from the origin, on which the limes inferior of
the modulus of the closed contours corresponding to the strata have
a limit, then there is a subsequence of $\mu_{n}$ such that
$v_{\mu_{n}} \rightarrow v_{0}$ on this arc. Assuming existence of a
point in (\ref{HE}) where $v_{\mu}$ is finite, we can find the arc
$\varrho$ such that $lim_{\mu \in \varrho} v_{\mu}=v_{0}$. The
conclusion is that we have analytic continuation in this case. Thus,
assuming the conditions on $V$ as above and that $M \bot W$, we can
use the path in Iversen's proof to derive hypoellipticity.

\newtheorem{main}[b_he-psdo]{ Proposition }
\begin{main}
 Given the representation $f=V(M,W)$ with $V_{1} \bot V_{2}$, assume $V^{-N}$ is the restriction to strata $X_{-N}$
and that the stratification of $(\Bm)$ has property $(M_{2})$. Then we have that if $(\Top-\mathcal{L}) V^{-N}$ is algebraic, as $T \rightarrow \infty$, there is no closed contour contributing microlocally to the symbol. 

\end{main}

\subsection{Algebraic continuation and orthogonality}
The information on the closed contour giving microlocal contribution
is in $\mathcal{T} V$, in the
manner that if $(\mathcal{T}-\mathcal{L}) V^{-N}$ is algebraic, there is no
possibility for presence of a closed contour that would contribute
and we can conclude that $V_{2} \prec \prec V_{1}$. Assume now that $V^{(N)} \in (\Bm)'
\subset \DL$, such that $\widehat{V}^{(N)}=P_{N}(T) f_{0}$, with
$f_{0}$ continuous (very regular). We can now assume
$\widehat{V}^{-(N)} \sim \frac{1}{P_{N}} g_{0}$, where $g_{0}$ has the same property as $f_{0}$. Our proposition is
that if $\mathcal{T}(\vartheta) \sim_{0} \alpha \mathcal{L} \vartheta$, with
$\alpha$ holomorphic, then for a hypoelliptic operator,
$\alpha/P_{N}$ will be algebraic, as $N \rightarrow \infty$.
Presence of a contributing closed contour means that $\alpha$ will be exponential.
We will argue that for $V$ corresponding to a hypoelliptic operator,
we have that $\mathcal{(T-L)} V^{-(N)}$ is algebraic in $T-$ infinity.
Assume $\alpha=q_{N}+r_{N}$, where $q_{N}$ is a polynomial. The
condition on involution means that $r_{N}/P_{N}$ behaves like the
symbol, in the tangent space. If $\mid r_{N}
\mid < \epsilon$, then according to Nullstellensatz, we have that
$\mid \alpha /P_{N} \mid < \epsilon$ as $T \rightarrow \infty$.
Thus, the condition that $P_{N}$ is reduced in the $T$-infinity and
$r_{N}$ small, means that $\alpha/P_{N}$ is algebraic in the
$T$-infinity, as $N \rightarrow \infty$. If we assume $\alpha/P_{N}$
bounded in the $T$-infinity, then we have that also $V^{-(N)}$ has
type $0$. The conclusion is given that $M \bot W$ and
$M=e^{\varphi}W$, if $F$ is algebraic, then we must have
$V_{1}=e^{\Phi}V_{2}$ and that $e^{-\Phi} \rightarrow 0$ as
$e^{-\varphi} \rightarrow 0$, as $T \rightarrow \infty$. This means
for the continuation of $F$ using $\Top$, that the orthogonality is preserved.

\subsection{Final remarks on hypoellipticity}
We have seen that under the condition that $M \bot W$ implies $V_{1} \bot V_{2}$, then we have that absence of closed contour that contributes microlocally is equivalent with the proposition $(\Top-\mathcal{L}) V^{-N}$ is algebraic.
Through the condition we have that $V$ preserves parabolic approximations.

\vsp

It is known that for a symbol  such that $f^{N} \in (I_{HE})$, for some integer $N$, we have $\Delta(f^{2}) \subset \Delta(f)$. We can give an interpolation problem for the iterates, 
$\psi^{j}+\varphi_{N}^{j}=M_{N}(\psi^{j})$ and $e^{\psi^{j}}=f^{j}$. As the kernel to the parametrix gets smaller as $j$ increases, why $\mbox{ Im }f \equiv 0$ on $V_{N}^{j}$, for all $N$ when $j$ is large. More precisely, if we assume 
$f(\xi) (\mbox{ Im }F) \rightarrow \infty$, as $\mid \xi \mid \rightarrow \infty$. Further, if we have $f^{\lambda} \in (I_{HE})$, as $\mid \xi \mid \rightarrow \infty$. Assume for simplicity, $\mid \xi \mid^{\delta} \leq C \mid f(\xi)^{N} \mid$, for $N$ positive and $\mid \xi \mid \rightarrow \infty$. The problem is now if we can find $\delta$ such that $\mid \xi \mid^{\delta} \leq C' \mid f(\xi) \mid$, as $\mid \xi \mid \rightarrow \infty$?. In this case we can choose $\delta=1/N$. In this manner we can prove that given that the real part has lower bound with exponent $\sigma$, then we can select $\delta=\sigma/N$ as exponent for the lower bound to the imaginary part. 

\vsp

We have noted that presence of lineality for the symbol, may result
in $\mbox{ Im }F$ in the space of hyperfunctions. We now note

\newtheorem{im}[b_he-psdo]{ Proposition }
\begin{im}
If $F$ is symmetric, entire and of finite type in $Exp$, then the
condition that $f$ represents a hypoelliptic operator, means that
for some $\lambda$, $(\mbox{ Im})^{\lambda} =\sum A_{j} F_{j}$ on a
domain of holomorphy, for constant coefficients and a global
pseudo-base $F_{j}$ representing the ideal of hypoelliptic operators.
\end{im}
Thus, symbols to hypoelliptic operators do not have imaginary part outside
the space of distributions and if hyper-function representation is
necessary, we must have contribution of lineality in the infinity.

\section{Examples}
\label{Ex}
\subsection{Introduction}
There is a big number of examples published in the literature (cf.
\cite{Treves99}) and we will deal with only some of them briefly
here. We are assuming the pseudo differential operator $P$ is
defined as $\lim_{\lambda \rightarrow 0}P_{\lambda}$, where the dependence of
$\lambda$ is locally algebraic in the symbol space. We are assuming a
dependence $\xi(x)$ through reciprocal polars, in this context, that
is $x_{T} \rightarrow \xi_{1/T}$. We are assuming the limit unique,
in the sense that $\frac{dP_{T}}{dT}=0$ implies $T=0$, for small $T$
(regular approximations). However, we may have
$\frac{d^{2}}{dT^{2}}P_{T}=0$, even locally. An operator that is the
regular limit of analytically hypoelliptic operators, with the
conditions that we have, is analytically hypoelliptic. We can use Proposition (\ref{Yodovich})
to construct an approximating sequence of symbols. The $L^{1}-$ dependence for parameter, means that the limit of $P_{\lambda}$ in operator space is continuous.

\subsection{Example 1}

Consider $P(x,y,\xi,\eta)=\xi^{2}+x^{2m} \eta^{2}$. This corresponds
to an operator analytically hypoelliptic in $\R$ with $n=2,m \geq
1$. If, for a constant $c_{T}$, $\xi_{1/T}^{2}=x_{T}^{2m}
\eta^{2}-c_{T}$, then $P_{T}$ is not analytically hypoelliptic, when
$c_{T} = 0$. But, when $c_{T} \neq 0$, we have that $P_{T}$ is
analytically hypoelliptic and the limit $P$ is a limit of
analytically hypoelliptic operators.

\subsection{ Example 2}

Consider $P(x,y,z,\xi,\eta,\nu)=\xi^{2}+x^{2}\eta^{2}+\nu^{2}$. This
corresponds to an operator not analytically hypoelliptic in $0$ in
$\R$, for $n=3$. For this reason we consider
$\xi_{1/T}^{2}+x^{2}_{T}\eta^{2}+\nu_{1/T}^{2}=c_{T}$, for a
(possibly zero) constant $c_{T}$ such that $c_{T} \rightarrow 0$ as
$T \rightarrow 0$ and such that for a constant $c_{T}'$, as $T
\rightarrow 0$,
$x_{T}^{2}\eta^{2}=c_{T}'(\xi_{1/T}^{2}+\nu_{1/T}^{2})$. We now have
constant surfaces through a suitable choice of $\eta$, which means
that the operator is not analytically hypoelliptic in $0$.

\subsection{ Example 3}

Consider $P(x,y,z,\xi,\eta,\nu)=\xi^{2}+ \big[ \eta +
(\frac{1}{3}x^{2}+xy^{2})v \big]^{2}=\xi^{2}+ \big[ \eta +
x(\frac{1}{3}x^{2}+y^{2})v \big]^{2}$. As we have assumed real
arguments, it can be proved that $P$ is a regular limit of
analytically hypoelliptic operators in $(x,
\xi)$,$(y,\eta)$,$(z,\nu)$ and in combinations of these.

\section{ The mapping $(I) \rightarrow Op(I)$ }
Assume $T=\{ y=F \}$ a transversal manifold, that is for a
submersion $p$, $\mbox{ ker }p=T_{x} \bigoplus L_{x}$, for a
manifold $L$. For instance, assume $h$ such that $dh(f)=0$ implies $f=const.$, then if $dh$ is locally algebraic (in the parameter), we must have where $f=const.$, that $dh(f)=const.$. Assuming on an irreducible component in $\{ f=const \}$, there is at least one point where $dh(f)=0$ we can conclude that $f=const$ implies $dh(f)=0$. A sufficient condition for this is that the tangent set $(I_{\mu})=\{\zeta \quad  h(f)=\mu f \}$ exists and has irreducible.components. Note that under these conditions, the first surfaces have a locally algebraic definition and the complement sets are assumed locally analytic. We can thus assume
that in a neighborhood of $df=0$, we have that $dh(f)=0$ gives a regular set. If we consider transversals on the form $d h(f)=\rho d f$, for a locally regular function $\rho$, we can form the extended transversal as a Baire function. 

\vsp

Assume $\Gamma_{A}$ is the boundary given by $d F=0$ and $\Gamma_{\mu}$ is given by 
$\{ \zeta \qquad \log  d F=0 \}$. If $\gamma$ is transversal to
$\Gamma_{A}^{U}(=\Gamma_{A} \cap U$), we have
$$ \Phi_{dF}=a_{\gamma}(z_{-},\zeta_{k})=\int_{\Gamma_{A}^{U}}
\frac{d_{z} F(z_{-}+ \zeta)}{\zeta - \zeta_{k}} d \zeta$$ We then
have, if $\widetilde{\Phi_{dF}}$ is the analytic continuation over
transversals, that $\widetilde{\Phi_{dF}}-\Phi_{dF} \rightarrow
dF=0$ on $\Gamma_{A}$. Further,
$d_{\zeta}(\widetilde{\Phi_{dF}}-\Phi_{dF})=0$ at isolated points.

\vsp

Now consider $F(\zeta,z)$ analytic for $\zeta$ bounded and $z$
large, such that $F(\zeta,z) \rightarrow 0$ as $z \rightarrow
\infty$. We assume
$$ F(\zeta,z)=\frac{d_{z} F(z+\zeta)}{d_{z} F(z)}-1$$
Using Cousin (cf. \cite{Cousin95}, Ch.4), there is a polynomial $Q$
in $1/t$ such that $\mid F(\zeta,z) \mid < \epsilon +Q(1/t)$ on
$\Delta_{t,\epsilon}$, a conical neighborhood of the lineality. We
say that $\zeta$ preserves a constant value for $F(\zeta,z)$. Thus,
$$ \mid d_{z}F(\zeta + z)-d_{z}F(z) \mid <(\epsilon_{t} +Q(1/t))
\mid d_{z}F(z) \mid$$ with $Q(1/t) \rightarrow 0$ as $t \rightarrow
\infty$. Assume further $\mid dF \mid < \mid dF + dL \mid$, where
$\mbox{ ker }p=\{ dF + dL \}$. Then formally $d_{z}F(\zeta + z)
\sim_{m}(dF + dL)$ and $dL \sim_{m} 0$. Within monotropy, we thus
have slow oscillation in the limit $z \rightarrow \infty$. We have
that $d_{z}F$ is holomorphic with respect to $\{ d_{z}F=0 \}$
according to $d_{z}F=\mbox{ tr }\widetilde{\Phi_{dF}}$, where
$\widetilde{\Phi_{dF}}$ corresponds to to $d_{z}L=L_{z}$ in the
transversal decomposition. If $F$ is a minimally defining function of
$N(J_{h})$, we have $dF \neq 0$ on this set. If $\mu$ is such that
$h(F)=\mu F$, we have that $\widetilde{\Phi_{F}}-\Phi_{F}
\rightarrow 0$ on $N(J_{h})$, why $\Phi_{F}$ is analytic over the
characteristic set $\{ F=0 \}$. Note that
$$ \mid \frac{d_{z} F(z+\zeta)}{dF(z)+dL(\zeta)}-1 \mid < \epsilon
\Rightarrow \mid \frac{dL(z)}{dF(z)} \mid <\epsilon$$ as $z
\rightarrow \infty$ along a transversal $\gamma$ emanating from the
origin. Consider first $L_{z}$ and the corresponding
$\widetilde{\Phi_{dF}}$. This corresponds through the inverse
Fourier-Borel transform, to an analytic functional that allows real
support, why it is sufficient to consider the real space and
$b_{\Gamma}$. Assume $u$ an entire function on a univalent domain
and zero on points for multi-valentness and on $\Gamma_{\mu}$. We
shall see that the corresponding form defines a good contour for the
associated pseudo differential operator. This is according to (cf.
\cite{Sjostrand82}) given by a realization with regularizing action
on $\mathcal{D}'$ why we have no loss of generality from the
conditions on the zero-set in the approach. Consider now the form
$\theta \text{:} \frac{2}{i} \frac{d \varphi}{d x} + iR
\overline{(x-y)}$ with $\mid x-y \mid \leq r$. Assume $y=x-\zeta$
and $\mid \zeta \mid \leq r$. We then have
$$ e^{\lambda (\varphi(x+\zeta)-\varphi(x))-2 \lambda \mbox{ Re }
\frac{d \varphi}{dx} \cdot \zeta}e^{- \mid \zeta \mid^{2}}=\big[
\frac{u^{\lambda}(x+\zeta)}{u^{\lambda}(x)} e^{-2 \mbox{  Re } \big[
\big(  \sum_{j} \frac{\delta_{j} u^{\lambda}}{u^{\lambda}} \cdot
\zeta \big) \big]} \big] e^{-R \mid \zeta \mid^{2}}$$ The slow
oscillation that we already established (within monotropy) implies
that particularly $\mid \delta_{j} u^{\lambda}/u^{\lambda} \mid$
bounded, as $\mid x \mid \rightarrow \infty$

\vsp

We conclude that for $\mid \zeta \mid$ bounded, as $\mid x
\mid \rightarrow \infty$, have that the bracket tends to $1$, as
$\mid x \mid \rightarrow \infty$ and we have a good contour
$\Gamma(x)$ for the form $\theta$. The pseudo differential operator
can be realized through
$$\widetilde{H_{\mu}}u_{\lambda}(x)=C_{\lambda} \int_{\Gamma(x)}
e^{i \lambda (x-y) \cdot \xi} L_{z}^{\mu}(x,y,\xi) u_{\lambda}(y) dy
d \xi$$ where $L_{z}^{\mu}$ has compact support and the operator
acts $\mathcal{D}' \rightarrow C^{\infty}$. Finally, we have the
case with $T_{z}$. Using Weierstrass theorem, we can find a
polynomial $P_{\mu,c}$ which include the foliation in its zero's.
For any polynomial, we have that $\delta_{j}P_{\mu} \prec P_{\mu}$,
why the second term is bounded, for $\lambda$ finite. The first term
is bounded by slow oscillation as before. We have apparently a good
contour also in this case.

 \vsp

We conclude that given the foliation for the symbol $\{ f=c \}$
and the tangent spaces $(I_{\mu})$, we can realize the pseudo
differential operator as a locally polynomial operator where the
polynomial part of the operator has zero's on the foliation. This is
a Levi decomposition of the operator
$$ A u_{\lambda}=\sum_{j,\mu}(P_{\mu,c_{j}} +
\widetilde{H_{\mu,c_{j}}}) u_{\lambda}$$

\bibliographystyle{amsplain}
\bibliography{sing}
\end{document}